\documentclass{amsart}
\usepackage{amssymb}
\usepackage{mathrsfs}
\usepackage{stmaryrd}
\usepackage{imakeidx}
\usepackage[colorlinks = true]{hyperref}
\usepackage{enumitem}
\usepackage{comment}
\input xy \xyoption {all}
\usepackage{bbm}

\usepackage{tikz}
\usetikzlibrary{cd}

\usepackage{dsfont}

\usepackage{rotating}
\usepackage{lscape}

\newcommand{\Z}{\mathbb{Z}}

\newcommand{\F}{\mathbb{F}}
\newcommand{\Q}{\mathbb{Q}}
\newcommand{\K}{\mathbb{K}}

\newcommand{\bk}{\Bbbk}
\newcommand{\Ga}{\mathbb{G}_{\mathrm{a}}}
\newcommand{\Gm}{\mathbb{G}_{\mathrm{m}}}

\newcommand{\Db}{D^{\mathrm{b}}}

\newcommand{\Perv}{\mathrm{Perv}}

\newcommand{\hatstar}{\mathbin{\widehat{\star}}}

\newcommand{\Av}{\mathsf{Av}}
\newcommand{\For}{\mathsf{For}}
\newcommand{\IC}{\mathscr{I}\hspace{-1pt}\mathscr{C}}

\newcommand{\id}{\mathrm{id}}

\newcommand{\simto}{\xrightarrow{\sim}}

\DeclareMathOperator{\Hom}{Hom}
\DeclareMathOperator{\Ext}{Ext}
\DeclareMathOperator{\End}{End}

\newcommand{\Iw}{\mathrm{I}}
\newcommand{\Iwu}{\mathrm{I}_{\mathrm{u}}}
\newcommand{\Fl}{\mathrm{Fl}}
\newcommand{\tFl}{\widetilde{\mathrm{Fl}}}
\newcommand{\Gr}{\mathrm{Gr}}
\newcommand{\Wf}{W_{\mathrm{f}}}
\newcommand{\fW}{{}^{\mathrm{f}} W}

\newcommand{\IW}{\mathcal{IW}}

\newcommand{\bG}{\mathbf{G}}
\newcommand{\bB}{\mathbf{B}}
\newcommand{\bU}{\mathbf{U}}
\newcommand{\bT}{\mathbf{T}}
\newcommand{\bg}{\mathbf{g}}
\newcommand{\bb}{\mathbf{b}}
\newcommand{\bu}{\mathbf{u}}
\newcommand{\bt}{\mathbf{t}}
\newcommand{\bR}{\mathbf{R}}

\renewcommand{\tilde}{\widetilde}

\makeatletter
\def\lotimes{\@ifnextchar_{\@lotimessub}{\@lotimesnosub}}
\def\@lotimessub_#1{\mathchoice{\mathbin{\mathop{\otimes}^L}_{#1}}%
  {\otimes^L_{#1}}{\otimes^L_{#1}}{\otimes^L_{#1}}}
\def\@lotimesnosub{\mathbin{\mathop{\otimes}^L}}
\makeatother

\makeatletter
\def\lboxtimes{\@ifnextchar_{\@lboxtimessub}{\@lboxtimesnosub}}
\def\@lboxtimessub_#1{\mathchoice{\mathbin{\mathop{\boxtimes}^L}_{#1}}%
  {\boxtimes^L_{#1}}{\boxtimes^L_{#1}}{\boxtimes^L_{#1}}}
\def\@lboxtimesnosub{\mathbin{\mathop{\boxtimes}^L}}
\makeatother

\newcommand{\scO}{\mathscr{O}}

\newcommand{\pH}{{}^{\mathrm{p}} \hspace{-1.5pt} \mathscr{H}}
\newcommand{\Rep}{\mathrm{Rep}}

\newcommand{\rInd}{{}^* \hspace{-0.5pt} \mathsf{Ind}}
\newcommand{\lInd}{{}^! \hspace{-0.5pt} \mathsf{Ind}}

\newcommand{\Groth}{\widetilde{\bg}}
\newcommand{\Coh}{\mathrm{Coh}}
\newcommand{\QCoh}{\mathrm{QCoh}}

\newcommand{\sZ}{\mathscr{Z}}

\newcommand{\sm}{\mathsf{m}}
\newcommand{\st}{\mathsf{t}}
\newcommand{\su}{\mathsf{u}}
\newcommand{\Wak}{\mathscr{W}}
\newcommand{\Wakop}{\mathscr{W}^{\mathrm{op}}}
\newcommand{\asp}{\mathrm{asp}}

\newcommand{\wt}{\mathsf{W}}
\newcommand{\sfD}{\mathsf{D}}
\newcommand{\sfP}{\mathsf{P}}

\newcommand{\Ind}{\mathrm{Ind}}

\newcommand{\sfZ}{\mathsf{Z}}

\newcommand{\scA}{\mathscr{A}}
\newcommand{\scB}{\mathscr{B}}

\newcommand{\scF}{\mathscr{F}}
\newcommand{\scG}{\mathscr{G}}

\newcommand{\scL}{\mathscr{L}}

\newcommand{\bX}{\mathbf{X}}
\newcommand{\bY}{\mathbf{Y}}

\newcommand{\Loop}{\mathrm{L}}
\newcommand{\Sat}{\mathsf{Sat}}
\newcommand{\rmZ}{\mathrm{Z}}

\numberwithin{equation}{section}
\numberwithin{figure}{section}
\newtheorem{thm}{Theorem}[section]
\newtheorem{lem}[thm]{Lemma}
\newtheorem{prop}[thm]{Proposition}
\newtheorem{cor}[thm]{Corollary}

\theoremstyle{definition}

\theoremstyle{remark}
\newtheorem{rmk}[thm]{Remark}

\title[Affine Hecke category and regular unipotent centralizer]{Modular affine Hecke category and regular unipotent centralizer}

\author[R. Bezrukavnikov]{Roman Bezrukavnikov}
\address{Department of Mathematics \\ Massachusetts Institute of Technology \\ Cambridge, MA \\ 02139 \\ USA.}
\email{bezrukav@math.mit.edu}

 \author[S.~Riche]{Simon Riche}
 \address{Universit\'e Clermont Auvergne, CNRS, LMBP, F-63000 Clermont-Ferrand, France.}
 \email{simon.riche@uca.fr}
 
 \author[L.~Rider]{Laura Rider}
 \address{Department of Mathematics, University of Georgia, Athens Georgia 30602, USA.}
 \email{laurajoy@uga.edu}
 
\begin{document}

\begin{abstract}
In this paper we provide, under some mild explicit assumptions, a geometric description of the category of representations of the centralizer of a regular unipotent element in a reductive algebraic group in terms of perverse sheaves on the Langlands dual affine flag variety. This equivalence is suggested and motivated by the ``geometric Langlands'' philosophy, and is used in later work~\cite{reg-quotient-pt2, reg-quotient-pt3} to construct equivalences of categories relating various geometric incarnations of the affine Hecke algebra of the given reductive group.
\end{abstract}

\maketitle


\section{Introduction}

\subsection{Towards an equivalence between the two modular categorifications
of the affine Hecke algebra}

The pre\-sent paper is the first one in a series dedicated to proving a modular variant (i.e.~a variant for coefficients of positive characteristic) of the equivalence between the two categorifications
of the affine Hecke algebra attached to a reductive algebraic group constructed by the first author for $\overline{\mathbb{Q}}_\ell$-coefficients~\cite{be}.
Such a variant is motivated by applications in modular representation theory; in particular, in~\cite{reg-quotient-pt3} we derive a proof of the Finkelberg--Mirkovi\'c conjecture~\cite{fm} describing the principal block of representations of a connected reductive algebraic group over a field of positive characteristic in terms of perverse sheaves on the Langlands dual affine Grassmannian. 

The origin of this construction is the fact that the Hecke algebra of the affine Weyl group attached to a reductive algebraic group admits two natural (and useful) geometric realizations: one (due to Iwahori--Matsumoto~\cite{im}) in terms of Iwa\-hori-bi-invariant functions on the loop group of the given group, and one (due to Kazhdan--Lusztig~\cite{kl} and Ginzburg~\cite{cg}) in terms of the equivariant K-theory of the Steinberg variety of the Langlands dual group. The equivalence of~\cite{be} ``lifts'' the isomorphism between these algebras to an equivalence of monoidal categories relating the Iwahori-equivariant derived category of the affine flag variety of the given group to a certain derived category of equivariant coherent sheaves on the Langlands dual Steinberg variety, both of the categories being defined over $\overline{\mathbb{Q}}_\ell$. (The precise statement of the equivalence requires some care, and will not be reviewed here.)

The construction of this equivalence in~\cite{be} was based on a previous work of S.~Arkhipov and the first author~\cite{ab} which provided an equivalence of categories ``lifting'' an isomorphism between two geometric realizations of the \emph{antispherical module} for the affine Hecke algebra; see the introduction of~\cite{ab} for more details. It turns out that adapting the constructions of~\cite{ab} to the modular setting presents several technical issues that we cannot overcome at the moment in general.\footnote{In~\cite[\S 8]{ar-book}, P.~Achar and the second author explain how to generalize this construction in case $G=\mathrm{GL}(n)$ and the characteristic $\ell$ of the coefficient field satisfies $\ell > \frac{1}{2} \binom{n}{\lfloor n/2 \rfloor}$.} Although it will be based on similar ingredients, our strategy will therefore be different: in addition to the constructions of local geometric Langlands theory, we
will employ Soergel's method of describing the categories we want to compare in algebro-combinatorial
terms.

\subsection{The regular quotient and representations of the regular unipotent centralizer}

The goal of the present paper and its sequel~\cite{reg-quotient-pt2} is to lay the foundations for this algebro-combinatorial description on the constructible side. 
The role of Soergel bimodules in the present paper is played by the category
of representations of the centralizer of a regular unipotent element; this is motivated
by~\cite{ab} where a similar construction is carried out for $\overline{\Q}_\ell$-coefficients, see also~\cite{dodd} where affine Soergel bimodules are related
to centralizers of regular elements. Thus, following the pattern of~\cite{soergel-Kat}, we
establish an equivalence between a Serre quotient of the category of Iwahori-equivariant $\bk$-perverse sheaves on the affine flag variety
with one (if $G$ is semisimple and simply connected) simple object and representations of the regular unipotent centralizer in the dual group defined over $\bk$;
this can be thought of as an analogue of Soergel's Endomorphismensatz, in its interpretation from~\cite{bbm}. (Here $\bk$ is an algebraic closure of $\mathbb{F}_\ell$ for some prime $\ell$.)
In~\cite{reg-quotient-pt2} we bootstrap this construction to an equivalence between appropriate categories of
constructible $\bk$-sheaves on the affine flag variety and coherent
sheaves on the Steinberg variety, establishing an analogue of Soergel's Struktursatz.

To describe more explicitly the equivalence that we prove here, let us introduce more notation. We will denote by $G$ a connected reductive algebraic group over an algebraically closed field of characteristic $p>0$. Let us choose a Borel subgroup $B \subset G$ and a maximal torus $T \subset B$. The choice of $B$ determines an Iwahori subgroup $\Iw$ in the arc group $\Loop^+ G$ of $G$, and we denote by $\Fl_G$ the associated affine flag variety. This ind-scheme admits a natural action of $\Iw$, and we can therefore consider, for $\bk$ an algebraic closure of $\mathbb{F}_\ell$ with $\ell \neq p$, the $\Iw$-equivariant derived category
\[
\sfD_{\Iw,\Iw} := \Db_{\Iw}(\Fl_G,\bk).
\]
This category is endowed with the perverse t-structure, whose heart will be denoted $\sfP_{\Iw,\Iw}$, and it admits a natural monoidal structure with product given by convolution, denoted $\star_\Iw$. 
Note that $\star_\Iw$ does {\em not} restrict to a monoidal product on $\sfP_{\Iw,\Iw}$; however we can construct an abelian monoidal category out of these data as follows. If we denote by $\sfP_{\Iw,\Iw}^0$ the quotient of the category $\sfP_{\Iw,\Iw}$ by the Serre subcategory generated by the simple objects whose support has positive dimension, then it is not difficult to check (see~\S\ref{ss:regular-quotient}) that the assignment $(\scF,\scG) \mapsto \pH^0(\scF \star_\Iw \scG)$ descends to an exact monoidal product on $\sfP_{\Iw,\Iw}^0$, which we denote by $\star^0_\Iw$.

On the dual side we consider the Langlands dual reductive group $G^\vee_\bk$ over $\bk$, a regular unipotent element $\su$, and the corresponding centralizer subgroup $\rmZ_{G^\vee_\bk}(\su)$. (Under the assumptions we will impose, this group scheme is smooth.) We denote by $\Rep(\rmZ_{G^\vee_\bk}(\su))$ the abelian category of finite-dimensional representations of $\rmZ_{G^\vee_\bk}(\su)$; this category has a natural monoidal structure given by the tensor product of representations.

Let us denote by $\mathfrak{R}$ the root system of $(G,T)$, and by $\mathfrak{R}^\vee$ the corresponding coroot system. We will assume that: 
\begin{itemize}
\item either $\ell$ is very good, or
$X^*(T)/\Z\mathfrak{R}$ is free (i.e.~$G^\vee_\bk$ has simply connected derived subgroup) and $X_*(T)/\Z\mathfrak{R}^\vee$ has no $\ell$-torsion;
\item $\ell$ is bigger than some explicit bound depending on $\mathfrak{R}$ and given in Figure~\ref{fig:bounds} below.
\end{itemize}
Under these assumptions,
the main result of this paper (Theorem~\ref{thm:main}) states that there exists an equivalence of abelian monoidal categories
\begin{equation}
\label{eqn:main-intro}
(\sfP_{\Iw,\Iw}^0, \star^0_\Iw) \cong (\Rep(\rmZ_{G^\vee_\bk}(\su)), \otimes_\bk)
\end{equation}
which satisfies a natural compatibility property with the geometric Satake equivalence and Gaitsgory's central functor (see~\S\ref{ss:intro-strategy}).

\begin{rmk}
In case $\bk=\overline{\Q}_\ell$, the equivalence~\eqref{eqn:main-intro} (or in fact a slightly less explicit variant) is a special case of the main result of~\cite{bez}, and it plays a technical role in a proof in~\cite{ab}. The equivalence itself is stated in a remark in~\cite{ab}, and is proved explicitly in~\cite[\S 7.2]{ar-book}.
\end{rmk}

In~\cite{reg-quotient-pt2} we will use this theorem to construct a fully faithful ``functor $\mathbb{V}$'' for certain tilting monodromic perverse sheaves on the natural $T$-torsor $\widetilde{\Fl}_G \to \Fl_G$, with codomain the category of representations of the regular centralizer group scheme (whose fiber over the base point is $\rmZ_{G^\vee_\bk}(\su)$). The latter category is related to the category of equivariant coherent sheaves on the Steinberg variety of $G^\vee_\bk$ (via restriction to the regular locus) in~\cite{reg-quotient-pt3}.


\subsection{Strategy of proof}
\label{ss:intro-strategy}

Our construction of the equivalence~\eqref{eqn:main-intro} will rely in a crucial way on the geometric Satake equivalence and Gaitsgory's theory of \emph{central sheaves}. Let $\Gr_G$ denote the affine Grassmannian of $G$. This ind-scheme admits a natural action of $\Loop^+ G$, and we can consider the category $\sfP_{\Loop^+G, \Loop^+G}$ of $\Loop^+ G$-equivariant \'etale $\bk$-perverse sheaves on $\Gr_G$. This category admits a natural convolution product $\star_{\Loop^+G}$, and the geometric Satake equivalence~\cite{mv} provides an equivalence of abelian monoidal categories
\begin{equation}
\label{eqn:Satake-intro}
(\sfP_{\Loop^+G, \Loop^+G},\star_{\Loop^+G}) \simto (\Rep(G^\vee_\bk), \otimes_\bk).
\end{equation}

The relation between perverse sheaves on $\Gr_G$ and on $\Fl_G$ is provided by Gaitsgory's functor
\[
\mathsf{Z} : \sfP_{\Loop^+G, \Loop^+G} \to \sfP_{\Iw,\Iw}.
\]
This functor, defined in terms of nearby cycles, provides a ``categorical lift'' of Bernstein's description of the center of the affine Hecke algebra, and the perverse sheaves in its image have a number of very favorable properties, studied in particular in~\cite{gaitsgory,gaitsgory-app,ab} and reviewed in~\cite{ar-book}. Using a general lemma from~\cite{bez} (or, more precisely, a slight extension of this lemma), the properties of this functor provide a unipotent element $\su \in G^\vee_\bk$, a subgroup scheme $H \subset \rmZ_{G^\vee_\bk}(\su)$, and an equivalence of monoidal categories between a certain subcategory $\widetilde{\sfP}_{\Iw,\Iw}^0$ of $\sfP_{\Iw,\Iw}^0$ and the category of finite-dimensional representations of $H$, under which the composition of $\mathsf{Z}$ with the quotient functor $\sfP_{\Iw,\Iw} \to \sfP_{\Iw,\Iw}^0$ corresponds to the composition of the Satake equivalence~\eqref{eqn:Satake-intro} with the restriction functor $\Rep(G^\vee_\bk) \to \Rep(H)$. Most of the new material in this paper will then be used to prove the following claims:
\begin{enumerate}
\item 
\label{it:u-reg}
$\su$ is regular;
\item 
\label{it:intro-tP-P}
the subcategory $\widetilde{\sfP}_{\Iw,\Iw}^0$ is the whole of $\sfP_{\Iw,\Iw}^0$;
\item 
\label{it:intro-H-Zu}
the embedding $H \subset \rmZ_{G^\vee_\bk}(\su)$ is an equality.
\end{enumerate}

These proofs will use modular variants of some technical constructions from~\cite{ab}. In particular we will prove that, under appropriate assumptions, the images in the category of Iwahori--Whittaker perverse sheaves of the objects $\mathsf{Z}(\scF)$ with $\scF$ tilting in $\sfP_{\Loop^+G, \Loop^+G}$ are tilting.

\subsection{Contents}

In Section~\ref{sec:preliminaries} we 
recall properties of centralizers of regular unipotent elements in reductive groups,
and provide a general criterion on subgroup schemes of such centralizers that will be used to prove statement~\eqref{it:intro-H-Zu} from~\S\ref{ss:intro-strategy}. In Section~\ref{sec:reconstructing} we prove a slight extension of the general lemma from~\cite{bez} alluded to in~\S\ref{ss:intro-strategy}. Section~\ref{sec:Satake} provides a reminder on the geometric Satake equivalence and Gaitsgory's central functor. In Section~\ref{sec:regular-quotient} we state our main result more precisely, and explain the construction of $\su$ and $H$. In Section~\ref{sec:regularity} we prove statement~\eqref{it:u-reg} from~\S\ref{ss:intro-strategy}. In Section~\ref{sec:asp-IW} we introduce the category of Iwahori--Whittaker perverse sheaves on $\Fl_G$, and explain the relation with $\sfP_{\Iw,\Iw}$ and $\sfP_{\Iw,\Iw}^0$. In Section~\ref{sec:central-tilting} we prove that the images of certain central sheaves in the category of Iwahori--Whittaker perverse sheaves are tilting. Finally, in Section~\ref{sec:proof-main} we use this to prove statements~\eqref{it:intro-tP-P} and~\eqref{it:intro-H-Zu} from~\S\ref{ss:intro-strategy}, and thus complete the proof of our main result.

\subsection{Acknowledgements}

R.B.~was supported by NSF Grants No.~DMS-1601953 and~DMS-2101507.
This project has received funding from the European Research Council (ERC) under the European Union's Horizon 2020 research and innovation programme (S.R., grant agreements No.~677147 and~101002592). L.R.~was supported by NSF Grant No.~DMS-1802378.

We thank M.~Korhonen for explaining the proof of Lemma~\ref{lem:fixed-points-qmin} to one of us and for a helpful discussion on centralizers of regular unipotent elements, A.~Bouthier and S.~Cotner for other useful discussions on centralizers of regular unipotent elements, and M.~Finkelberg for spotting some typos.

\section{Preliminaries on reductive groups}
\label{sec:preliminaries}


In this section we collect a number of general results on the geometry or representations of reductive algebraic groups which will be needed later in our arguments. We fix an algebraically closed field $\K$, whose characteristic we denote by $\ell$.

\subsection{Quotients of group schemes}
\label{ss:quotients}

In our considerations below we will have to consider some quotients of affine $\K$-group schemes of finite type by closed subgroup schemes, without any reducedness assumptions. For the reader's convenience, in this subsection we briefly recall the properties of this construction that we will use.

Let $H$ be an affine $\K$-group scheme of finite type.
\begin{itemize}
\item For any subgroup scheme $H' \subset H$, the fppf-sheafification of the functor $R \mapsto H(R)/H'(R)$ is represented by a $\K$-scheme of finite type, denoted $H/H'$; see~\cite[III, \S3, Th\'eor\`eme~5.4]{dg}.
\item We have $\dim(H/H')=\dim(H)-\dim(H')$; see~\cite[III, \S 3, Remarque 5.5(a)]{dg}.
\item If $H'$ is smooth then the natural map $H \to H/H'$ is smooth; see~\cite[III, \S 3, Corollaire~2.6]{dg}.
\item If $H$ is smooth, then the scheme $H/H'$ is smooth; see~\cite[III, \S 3, Proposition~2.7]{dg}.
\item If $H'$ is a normal subgroup, then the quotient $H/H'$ is an affine group scheme (of finite type); see~\cite[III, \S 3, Th\'eor\`eme~5.6]{dg}.
\end{itemize}

For $h \in H(\K)$ 
we will denote by $\rmZ_{H}(h)$ the scheme-theoretic centralizer of $h$ in $H$, i.e.~the scheme-theoretic fiber over $h$ of the morphism $H \to H$ defined by $g \mapsto ghg^{-1}$.

\subsection{Centralizers of regular unipotent elements in reductive groups}
\label{ss:centralizer-ureg}

From now on we assume that $\bG$ is a connected reductive algebraic group over $\K$, and denote by $u \in \bG(\K)$ a regular unipotent element. 
We denote by $\bB \subset \bG$ the unique Borel subgroup containing $u$, by $\bU$ its unipotent radical (so that $u \in \bU(\K)$) and choose a maximal torus $\bT \subset \bB$, so that multiplication induces an isomorphism $\bT \ltimes \bU \simto \bB$.
We will also denote by $\bR \subset X^*(\bT)$ the root  system of $(\bG,\bT)$, by $\bR^\vee \subset X_*(\bT)$ the associated coroots, 
and by $\Z\bR \subset X^*(\bT)$ and $\Z\bR^\vee \subset X_*(\bT)$ the root and coroot lattices. 
We will denote by $\rmZ(\bG)$ the (scheme-theoretic) center of $\bG$. We start with the following observation.

\begin{lem}
\label{lem:center-smooth}
The group scheme $\rmZ(\bG)$ is smooth iff the quotient $X^*(\bT)/\Z\bR$ has no $\ell$-torsion.
\end{lem}

\begin{proof}
It is well known that $\rmZ(\bG)$ is the scheme-theoretic intersection of the kernels of the roots of $(\bG,\bT)$; in other words it is the diagonalizable $\K$-group scheme whose group of characters is $X^*(\bT)/\Z\bR$. This implies the desired claim.
\end{proof}

The following statement gathers some statements regarding the centralizers of $u$ in $\bU$, $\bB$ and $\bG$ established in~\cite[Lemma~3.4 and Corollary~3.11]{cotner} (building on classical results of Steinberg) in a more general setting.

\begin{prop}
\label{prop:Z-ureg}
The group scheme $\rmZ_{\bU}(u)$ is smooth, and multiplication induces an isomorphism of group schemes
\[
\rmZ_{\bU}(u) \times \rmZ(\bG) \simto \rmZ_{\bB}(u).
\]
If moreover $X_*(\bT) / \Z\bR^\vee$ has no $\ell$-torsion, then the natural closed immersion
\[
\rmZ_{\bB}(u) \to \rmZ_{\bG}(u)
\]
is an isomorphism, and these group schemes are commutative. 

In particular, if neither $X^*(\bT)/\Z\bR$ nor $X_*(\bT) / \Z\bR^\vee$ has $\ell$-torsion, then $\rmZ_{\bG}(u)$ is smooth.
\end{prop}

\begin{rmk}
At the time when the first version of this paper was written, the rerence~\cite{cotner} was not yet available. That version contained an independent proof of a slightly weaker form of Proposition~\ref{prop:Z-ureg}.
\end{rmk}

\subsection{The unipotent cone and the multiplicative Springer resolution}
\label{ss:multiplicative-Springer}

We continue with the notation of~\S\ref{ss:centralizer-ureg}.
The \emph{multiplicative Springer resolution} is the induced scheme
\[
 \widetilde{\mathcal{U}} := \bG \times^{\bB} \bU,
\]
where $\bB$ acts on $\bU$ by conjugation. This scheme is a locally trivial bundle over the flag variety $\bG/\bB$, with fibers isomorphic to $\bU$; in particular it is smooth. Let us denote by $\mathcal{U} \subset \bG$ the unipotent cone, i.e.~the closed subvariety whose $\K$-points are the unipotent elements in $\bG$. (This variety is denoted $\mathscr{U}_{\bG}^{\mathrm{var}}$ in~\cite[\S 4.2]{cotner}.) The adjoint action of $\bG$ on itself induces a morphism $\widetilde{\mathcal{U}} \to \bG$ which factors through a morphism
\[
 \pi : \widetilde{\mathcal{U}} \to \mathcal{U}.
\]
This morphism is surjective (see~\cite[Theorem~11.10]{borel}), hence induces an embedding of $\bG$-modules
\begin{equation}
\label{eqn:embedding-U-tU}
\scO(\mathcal{U}) \hookrightarrow \scO(\widetilde{\mathcal{U}})
\end{equation}
where, for a scheme $X$, we write $\scO(X)$ for $\Gamma(X, \scO_X)$.

We will also denote by $\mathcal{U}_{\mathrm{reg}} \subset \mathcal{U}$ the regular locus; then, as explained in~\cite[Proof of Lemma~6.5]{cotner}, if $u \in \bG(\K)$ is as in~\S\ref{ss:centralizer-ureg},
the conjugation action on $u$ induces an isomorphism
\begin{equation}
\label{eqn:Ureg-quotient}
\bG/\rmZ_{\bG}(u) \simto \mathcal{U}_{\mathrm{reg}};
\end{equation}
in particular, using~\cite[\S I.5.10 and~Proposition~I.5.12]{jantzen} we deduce that we have an isomorphism of $\bG$-modules
\begin{equation}
\label{eqn:O-Ureg}
\scO(\mathcal{U}_{\mathrm{reg}}) \cong \mathrm{Ind}^{\bG}_{\rmZ_{\bG}(u)}(\K).
\end{equation}


The choice of $\bB$ determines a system of positive roots for $(\bG,\bT)$, chosen as the opposites of the roots appearing in the Lie algebra of $\bB$, and therefore a notion of dominant weights. For $\lambda \in X^*(\bT)$ we denote by $\scO_{\bG/\bB}(\lambda)$ the line bundle on $\bG/\bB$ naturally associated with $\lambda$, and by
$\scO_{\widetilde{\mathcal{U}}}(\lambda)$ its pullback under the projection $\widetilde{\mathcal{U}} \to \bG/\bB$.

\begin{lem}
\label{lem:line-bundle-nonzero}
If $\lambda \in X^*(\bT)$ is dominant, then there exists an embedding of $\bG$-modules $\Gamma(\bG/\bB, \scO_{\bG/\bB}(\lambda)) \hookrightarrow \Gamma(\widetilde{\mathcal{U}}, \scO_{\widetilde{\mathcal{U}}}(\lambda))$. In particular, $\Gamma(\widetilde{\mathcal{U}}, \scO_{\widetilde{\mathcal{U}}}(\lambda)) \neq 0$.
\end{lem}

\begin{proof}
The pushforward of the structure sheaf under the projection $\widetilde{\mathcal{U}} \to \bG/\bB$ is the quasi-coherent sheaf on $\bG/\bB$ associated with the $\bB$-module $\scO(\bU)$ (in the sense of~\cite[Chap.~I.5]{jantzen}). In particular, the natural morphism from $\scO_{\bG/\bB}$ is induced by the embedding $\K \hookrightarrow \scO(\bU)$, hence is injective. We deduce an embedding $\Gamma(\bG/\bB, \scO_{\bG/\bB}(\lambda)) \hookrightarrow \Gamma(\widetilde{\mathcal{U}}, \scO_{\widetilde{\mathcal{U}}}(\lambda))$.

The claim follows, since the domain of this map is nonzero if $\lambda$ is dominant.
\end{proof}

We set
\[
\widetilde{\mathcal{U}}_{\mathrm{reg}} = \pi^{-1}(\mathcal{U}_{\mathrm{reg}}).
\]
For each simple root $\alpha$ we denote by $\bU_{(\alpha)} \subset \bU$ the subgroup generated by the root subgroups associated with roots different from $-\alpha$ (in other words, the unipotent radical of the parabolic subgroup of $\bG$ containing $\bB$ associated with the subset $\{\alpha\}$ of the set of simple roots). Then $\bU_{(\alpha)}$ is stable under the $\bB$-action, so that we can define
\[
 \mathcal{D}_\alpha := \bG \times^{\bB} \bU_{(\alpha)}.
\]
The subvariety $\mathcal{D}_\alpha$ is a divisor in $\widetilde{\mathcal{U}}$, and we have
\[
 \scO_{\widetilde{\mathcal{U}}}(\mathcal{D}_\alpha)=\scO_{\widetilde{\mathcal{U}}}(-\alpha).
\]

For any sum of positive roots $\lambda$, written as $\lambda=\sum_\alpha n_\alpha \cdot \alpha$ (where $\alpha$ runs over the simple roots), we set
\[
 \mathcal{D}_\lambda := \sum_\alpha n_\alpha \cdot \mathcal{D}_\alpha.
\]
Then $\mathcal{D}_\lambda$ is an effective divisor in $\widetilde{\mathcal{U}}$ with
\begin{equation}
\label{eqn:line-bundle-Dlambda}
 \scO_{\widetilde{\mathcal{U}}}(\mathcal{D}_\lambda)=\scO_{\widetilde{\mathcal{U}}}(-\lambda),
\end{equation}
and if each coefficient $n_\alpha$ is positive the associated closed subscheme $\mathcal{Y}_\lambda \subset \widetilde{\mathcal{U}}$ satisfies
\[
 \widetilde{\mathcal{U}} \smallsetminus \mathcal{Y}_\lambda = \widetilde{\mathcal{U}}_{\mathrm{reg}},
\]
see e.g.~\cite[\S 4.1]{humphreys}. It follows that if $j_{\mathrm{reg}} : \widetilde{\mathcal{U}}_{\mathrm{reg}} \hookrightarrow \widetilde{\mathcal{U}}$ is the embedding, then we have
\begin{equation}
\label{eqn:pushforward-jreg}
 (j_{\mathrm{reg}})_* \scO_{\widetilde{\mathcal{U}}_{\mathrm{reg}}} = \varinjlim_{m \in \Z_{\geq 0}} \scO_{\widetilde{\mathcal{U}}}(m \cdot \mathcal{D}_\lambda)
\end{equation}
(where the maps $\scO_{\widetilde{\mathcal{U}}}(m \cdot \mathcal{D}_\lambda) \to \scO_{\widetilde{\mathcal{U}}}((m+1) \cdot \mathcal{D}_\lambda)$ are induced by the natural embedding $\scO_{\widetilde{\mathcal{U}}} \hookrightarrow \scO_{\widetilde{\mathcal{U}}}(\mathcal{D}_\lambda)$).

In the following statement we denote by $\Coh^{\bG}(\widetilde{\mathcal{U}})$, resp.~$\Coh^{\bG}(\widetilde{\mathcal{U}}_{\mathrm{reg}})$, the abelian category of $\bG$-equivariant coherent sheaves on $\widetilde{\mathcal{U}}$, resp.~$\widetilde{\mathcal{U}}_{\mathrm{reg}}$.

\begin{lem}
\label{lem:Ext1-jreg}
 Let $\lambda=\sum_\alpha n_\alpha \cdot \alpha$, where $n_\alpha>0$ for all $\alpha$. Then restriction to $\widetilde{\mathcal{U}}_{\mathrm{reg}}$ induces an isomorphism
 \[
  \varinjlim_{m \in \Z_{\geq 0}} \Ext^1_{\Coh^{\bG}(\widetilde{\mathcal{U}})}(\scO_{\widetilde{\mathcal{U}}}(m \cdot \lambda), \scO_{\widetilde{\mathcal{U}}}) \simto \Ext^1_{\Coh^{\bG}(\widetilde{\mathcal{U}}_{\mathrm{reg}})}(\scO_{\widetilde{\mathcal{U}}_{\mathrm{reg}}}, \scO_{\widetilde{\mathcal{U}}_{\mathrm{reg}}}).
 \]
\end{lem}

\begin{proof}
 The morphism $j_{\mathrm{reg}}$ is affine. Therefore, using adjunction we obtain an isomorphism
 \begin{multline*}
  \Ext^1_{\Coh^{\bG}(\widetilde{\mathcal{U}}_{\mathrm{reg}})}(\scO_{\widetilde{\mathcal{U}}_{\mathrm{reg}}}, \scO_{\widetilde{\mathcal{U}}_{\mathrm{reg}}})=\Ext^1_{\Coh^{\bG}(\widetilde{\mathcal{U}}_{\mathrm{reg}})}(j_{\mathrm{reg}}^*\scO_{\widetilde{\mathcal{U}}}, \scO_{\widetilde{\mathcal{U}}_{\mathrm{reg}}}) \\
  \cong \Ext^1_{\QCoh^{\bG}(\widetilde{\mathcal{U}})}(\scO_{\widetilde{\mathcal{U}}}, (j_{\mathrm{reg}})_* \scO_{\widetilde{\mathcal{U}}_{\mathrm{reg}}}),
 \end{multline*}
 where $\QCoh^{\bG}(\widetilde{\mathcal{U}})$ denotes the category of $\bG$-equivariant \emph{quasi}-coherent sheaves on $\widetilde{\mathcal{U}}$.
The derived functor $R\Gamma(\widetilde{\mathcal{U}},-)$ commutes with uniformly bounded below filtered direct limits (see e.g.~\cite[Lemma~3.9.3.1]{lipman}), and so does the functor of derived $\bG$-invariants, see~\cite[Lemma~I.4.17]{jantzen}. Since the functor $R\Hom_{\QCoh^{\bG}(\widetilde{\mathcal{U}})}(\scO_{\widetilde{\mathcal{U}}},-)$ is the composition of these two functors (see~\cite[Proposition~A.6]{mr1}), it also satisfies this property, which implies the desired claim in view of~\eqref{eqn:pushforward-jreg} and~\eqref{eqn:line-bundle-Dlambda}.
\end{proof}

The following claims are standard.

\begin{lem}
\label{lem:U-sc}
Assume that $\bG$ has simply connected derived subgroup.

\begin{enumerate}
\item
\label{it:functions-U}
 The variety $\mathcal{U}$ is normal, and the morphism
 \[
  \scO(\mathcal{U}) \to \scO(\mathcal{U}_{\mathrm{reg}})
 \]
induced by restriction is an isomorphism.
\item
\label{it:U-tU-reg}
The restriction of $\pi$ to $\widetilde{\mathcal{U}}_{\mathrm{reg}}$ is an isomorphism of varieties $\widetilde{\mathcal{U}}_{\mathrm{reg}} \simto \mathcal{U}_{\mathrm{reg}}$.
\item
\label{it:U-tU-Zar-Main}
The embedding~\eqref{eqn:embedding-U-tU} is an isomorphism.
\item
\label{it:functions-tU}
We have an isomorphism of $\bG$-modules
\[
\scO(\widetilde{\mathcal{U}}) \cong \mathrm{Ind}_{\rmZ_{\bG}(u)}^{\bG}(\K).
\]
\end{enumerate}
\end{lem}

\begin{proof}
Note that $\mathcal{U}$ and $\widetilde{\mathcal{U}}$ are unchanged if $\bG$ is replaced by its derived subgroup, so we can assume $\bG$ is semisimple (and simply connected).

 \eqref{it:functions-U} 
For the normality of $\mathcal{U}$, see~\cite[Theorem~4.24(iii)]{humphreys} or~\cite[Lemma~4.4(4)]{cotner}. Then the second claim follows using the algebraic Hartogs' lemma 
and the fact that the complement of $\mathcal{U}_{\mathrm{reg}}$ in $\mathcal{U}$ has codimension at least $2$ (which itself follows from the facts that $\mathcal{U}$ has finitely many orbits and that non-regular orbits have dimension at most $\dim(\mathcal{U})-2$, see~\cite[Remark in~\S 4.1]{humphreys}, or~\cite[Proof of Lemma~4.4(4)]{cotner}.) 
 
 \eqref{it:U-tU-reg} See~\cite[\S 6.3--6.4]{humphreys} or~\cite[Proof of Lemma~6.17]{cotner}.
 
 \eqref{it:U-tU-Zar-Main}
 The claim follows from Zariski's Main Theorem (see~\cite[Chap.~III, Proof of Corollary~11.4]{hartshorne}), whose assumptions are satisfied by~\eqref{it:functions-U} and~\eqref{it:U-tU-reg}.
 
 \eqref{it:functions-tU}
 This follows from~\eqref{it:functions-U},~\eqref{it:U-tU-Zar-Main} and~\eqref{eqn:O-Ureg}.
\end{proof}

\subsection{Tilting modules and subgroups of the centralizer of a regular unipotent element}
\label{ss:proof-criterion}

%
%
%
%



From now on we assume that:
\begin{enumerate}
\item $\bG$ has simply connected derived subgroup (i.e.~$X_*(\bT)/\Z\bR^\vee$ is free);
\item $\ell$ is good for $\bG$;
\item 
\label{it:ass-3}
$X^*(\bT)/\Z\bR$ has no $\ell$-torsion.
\end{enumerate}

\begin{lem}
\label{lem:klt}
We have
\[
\mathsf{H}^i(\widetilde{\mathcal{U}}, \scO_{\widetilde{\mathcal{U}}})=0 \quad \text{for any $i>0$,}
\]
and moreover the $\bG$-module $\scO(\widetilde{\mathcal{U}})$ admits a good filtration.
\end{lem}

\begin{proof}
Assume first that $\bG$ is semisimple. Under our assumption on $\ell$, $\widetilde{\mathcal{U}}$ identifies with the usual (additive) Springer resolution considered in~\S\ref{ss:fixed-points-Zu} below (see~\cite[Proposition~3]{klt} or Lemma~\ref{lem:Spr-resolution-isom} below), so that our claims are special cases of~\cite[Theorem~2]{klt} and~\cite[Theorem~7]{klt}. The general case follows, using the fact that $\widetilde{\mathcal{U}}$ is unchanged if $\bG$ is replaced by its derived subgroup $\mathscr{D}(\bG)$, and that a $\bG$-module $V$ admits a good filtration iff $V_{|\mathscr{D}(\bG)}$ admits a good filtration (see~\cite[Proposition~3.2.7]{donkin}).
\end{proof}


As in~\S\ref{ss:centralizer-ureg} we consider a regular unipotent element $u$ which belongs to $\bU(\K)$. Under our assumptions, by Lemma~\ref{lem:center-smooth} and Proposition~\ref{prop:Z-ureg} the group schemes $\rmZ(\bG)$, $\rmZ_{\bU}(u)$ and $\rmZ_{\bG}(u)$ are smooth.

\begin{lem}
\label{lem:criterion-equality}
 Let $\mathbf{H} \subset \rmZ_{\bG}(u)$ be a subgroup scheme containing $\rmZ(\bG)$. If for any finite-dimensional tilting $\bG$-module $V$ the embedding
 \[
   V^{\rmZ_{\bG}(u)}\hookrightarrow V^{\mathbf{H}} 
 \]
is an equality, then $\mathbf{H} = \rmZ_{\bG}(u)$.
\end{lem}

\begin{proof}
We will show that 
if $\mathbf{H} \subset \rmZ_{\bG}(u)$ is a proper subgroup scheme containing $\rmZ(\bG)$, then there exists a finite-dimensional tilting $\bG$-module $V$ such that the embedding $V^{\rmZ_{\bG}(u)}  \hookrightarrow V^{\mathbf{H}}$ is strict. 

By Proposition~\ref{prop:Z-ureg},
multiplication induces a group scheme isomorphism
\[
 \rmZ_{\bU}(u) \times \rmZ(\bG) \simto \rmZ_{\bG}(u),
\]
and these group schemes are commutative.
Since $\mathbf{H}$ is contained in $\rmZ_{\bG}(u)$ and contains $\rmZ(\bG)$, we deduce that we also have
\[
 (\rmZ_{\bU}(u) \cap \mathbf{H}) \times \rmZ(\bG) \simto \mathbf{H},
\]
and $\mathbf{H}$ is a normal subgroup in $\rmZ_{\bG}(u)$. 
Hence the quotient
 $\rmZ_{\bG}(u) / \mathbf{H}$
is a smooth affine group scheme of finite type (see~\S\ref{ss:quotients}).
Since this group scheme is nontrivial by assumption, it follows that $\scO(\rmZ_{\bG}(u) / \mathbf{H}) \neq \K$.

Let $L \subset \scO(\rmZ_{\bG}(u) / \mathbf{H})/\K$ be a simple $\rmZ_{\bG}(u) / \mathbf{H}$-submodule (where the action on $\scO(\rmZ_{\bG}(u) / \mathbf{H})$ is induced by right multiplication), and denote by $E$ its inverse image in $\scO(\rmZ_{\bG}(u) / \mathbf{H})$.
Since $\rmZ_{\bU}(u)$ is unipotent, 
$L$ must be the trivial module, and $E$ fits in an exact sequence of $\rmZ_{\bG}(u) / \mathbf{H}$-modules 
\begin{equation}
\label{eqn:ses-E}
 \K \hookrightarrow E \twoheadrightarrow \K.
\end{equation}

Recall that any $\rmZ_{\bG}(u)$-module $M$ determines a $\bG$-equivariant quasi-coherent sheaf $\mathscr{L}_{\bG/\rmZ_{\bG}(u)}(M)$ on $\bG/\rmZ_{\bG}(u)$, see~\cite[\S I.5.8]{jantzen}. In particular, if $f : \bG/\mathbf{H} \to \bG/\rmZ_{\bG}(u)$ is the projection, then we have
\[
 f_* \scO_{\bG/\mathbf{H}} \cong \mathscr{L}_{\bG/\rmZ_{\bG}(u)}(\scO(\rmZ_{\bG}(u) / \mathbf{H})),
\]
see~\cite[\S I.5.12, Remark 3]{jantzen}. 
The functor $\mathscr{L}_{\bG/\rmZ_{\bG}(u)}(-)$ is exact; in particular, the short exact sequence~\eqref{eqn:ses-E} provides a short exact sequence of $\bG$-equivariant coherent sheaves
\begin{equation}
\label{eqn:ses-E-L}
 \scO_{\bG/\rmZ_{\bG}(u)} \hookrightarrow \mathscr{L}_{\bG/\rmZ_{\bG}(u)}(E) \twoheadrightarrow \scO_{\bG/\rmZ_{\bG}(u)},
\end{equation}
and the embedding $E \hookrightarrow \scO(\rmZ_{\bG}(u) / \mathbf{H})$ provides an embedding
\begin{equation}
\label{eqn:embedding-E-O}
 \mathscr{L}_{\bG/\rmZ_{\bG}(u)}(E) \hookrightarrow f_* \scO_{\bG/\mathbf{H}}.
\end{equation}

In view of~\eqref{eqn:Ureg-quotient},
the exact sequence~\eqref{eqn:ses-E-L} can be considered as an exact sequence in $\Coh^{\bG}(\mathcal{U}_{\mathrm{reg}})$, or equivalently (see Lemma~\ref{lem:U-sc}\eqref{it:U-tU-reg}) in $\Coh^{\bG}(\widetilde{\mathcal{U}}_{\mathrm{reg}})$. After choosing a weight $\lambda$ as in Lemma~\ref{lem:Ext1-jreg}, which we furthermore assume to be dominant, 
this lemma ensures that there exists $m \in \mathbb{Z}_{\geq 0}$ and a short exact sequence
\begin{equation}
\label{eqn:ses-E-Springer}
 \mathscr{O}_{\widetilde{\mathcal{U}}} \hookrightarrow \mathscr{E} \twoheadrightarrow \scO_{\widetilde{\mathcal{U}}}(m \cdot \lambda)
\end{equation}
for some $m \in \Z_{\geq 0}$ whose image under $j_{\mathrm{reg}}^*$ is~\eqref{eqn:ses-E-L}. In particular, $\mathscr{E}$ is a vector bundle on $\widetilde{\mathcal{U}}$ such that $j_{\mathrm{reg}}^* \mathscr{E} = \mathscr{L}_{\bG/\rmZ_{\bG}(u)}(E)$. Since $\widetilde{\mathcal{U}}_{\mathrm{reg}}$ is open (hence dense) in $\widetilde{\mathcal{U}}$, the morphism
\[
 \mathscr{E} \to (j_{\mathrm{reg}})_* \mathscr{L}_{\bG/\rmZ_{\bG}(u)}(E)
\]
induced by adjunction is injective.
Taking global sections we deduce an embedding of $\bG$-modules
\[
 \Gamma(\widetilde{\mathcal{U}}, \mathscr{E}) \hookrightarrow \Gamma(\bG/\rmZ_{\bG}(u), \mathscr{L}_{\bG/\rmZ_{\bG}(u)}(E)).
\]
Composing this map with the morphism obtained from~\eqref{eqn:embedding-E-O} by taking global sections, we deduce an embedding 
of $\bG$-modules
\begin{equation}
\label{eqn:embedding-global-sec}
 \Gamma(\widetilde{\mathcal{U}}, \mathscr{E}) \hookrightarrow \scO(\bG/\mathbf{H}).
 \end{equation}
 
 Let us now come back to~\eqref{eqn:ses-E-Springer}, and denote by $\Rep^\infty(\bG)$ the category of all  (not necessarily finite-dimensional) algebraic $\bG$-modules. 
By Lemma~\ref{lem:klt} we have $\mathsf{H}^1(\widetilde{\mathcal{U}}, \scO_{\widetilde{\mathcal{U}}})=0$. Therefore, by taking global sections we obtain an exact sequence of $\bG$-modules
 \begin{equation}
 \label{eqn:ses-global-sec}
  \scO(\widetilde{\mathcal{U}}) \hookrightarrow \Gamma(\widetilde{\mathcal{U}}, \mathscr{E}) \twoheadrightarrow \Gamma(\widetilde{\mathcal{U}}, \scO_{\widetilde{\mathcal{U}}}(m \cdot \lambda)).
 \end{equation}
Now the second statement in Lemma~\ref{lem:klt} implies that
for any finite-dimensional tilting $\bG$-module $V$ we have
\[
\Ext^1_{\Rep^\infty(\bG)} \bigl( V,\scO(\widetilde{\mathcal{U}}) \bigr)=0.
\]
On the other hand, using Lemma~\ref{lem:line-bundle-nonzero} we see that there exists a (finite-dimensional) tilting $\bG$-module $V$ (e.g.~the indecomposable tilting $\bG$-module with highest weight $m\lambda$) such that
\[
 \Hom_{\Rep^\infty(\bG)} \bigl( V,\Gamma(\widetilde{\mathcal{U}}, \scO_{\widetilde{\mathcal{U}}}(m \cdot \lambda)) \bigr) \neq 0.
 \]
Applying the functor $\Hom_{\Rep^\infty(\bG)}(V,-)$ to~\eqref{eqn:ses-global-sec} we deduce a strict embedding
\[
 \Hom_{\Rep^\infty(\bG)}(V,\scO(\widetilde{\mathcal{U}})) \hookrightarrow \Hom_{\Rep^\infty(\bG)}(V,\Gamma(\widetilde{\mathcal{U}}, \mathscr{E}))
\]
and then, using~\eqref{eqn:embedding-U-tU} and~\eqref{eqn:embedding-global-sec}, a strict embedding
\[
\Hom_{\Rep^\infty(\bG)}(V, \scO(\mathcal{U})) \hookrightarrow \Hom_{\Rep^\infty(\bG)}(V,\scO(\bG/\mathbf{H})).
\]
By Lemma~\ref{lem:U-sc}\eqref{it:functions-U} and~\eqref{eqn:Ureg-quotient}, this map can be considered as a strict embedding
\[
\Hom_{\Rep^\infty(\bG)}(V, \scO(\bG/\rmZ_{\bG}(u))) \hookrightarrow \Hom_{\Rep^\infty(\bG)}(V,\scO(\bG/\mathbf{H})).
\]

Now, by Frobenius reciprocity we have isomorphisms
\begin{gather*}
\Hom_{\Rep^\infty(\bG)}(V, \scO(\bG/\rmZ_{\bG}(u))) \cong (V^*)^{\rmZ_{\bG}(u)}, \\
\Hom_{\Rep^\infty(\bG)}(V,\scO(\bG/\mathbf{H})) \cong (V^*)^{\mathbf{H}}.
\end{gather*}
Since $V^*$ is a tilting module, this finishes the proof.
\end{proof}


\subsection{Fixed points of the centralizer of a regular unipotent element}
\label{ss:fixed-points-Zu}

We continue with the notation and assumptions of~\S\ref{ss:proof-criterion}. Our goal in this subsection is to prove Proposition~\ref{prop:Vu-V0} below, which relates the dimensions of fixed points under $\rmZ_{\bG}(u)$ and $\bT$ for finite-dimensional $\bG$-modules admittting a good filtration. 
We denote by $\bg$, $\bb$, $\bu$, $\bt$ the Lie algebras of $\bG$, $\bB$, $\bU$, $\bT$.
The proof will involve the geometry of the ``additive'' version of the Springer resolution
%
%
\[
 \widetilde{\mathcal{N}} := \bG \times^{\bB} \bu,
\]
and of the Grothendieck resolution
\[
 \Groth:=\bG \times^{\bB} \bb.
\]
Both $\widetilde{\mathcal{N}}$ and $\Groth$ are vector bundles over the flag variety $\bG/\bB$; in particular they are smooth schemes. 

Under our present assumptions, the multiplicative and additive versions of the Springer resolution are isomorphic.

\begin{lem}
\label{lem:Spr-resolution-isom}
There exists an isomorphism of $\bG$-schemes
\[
\widetilde{\mathcal{U}} \simto \widetilde{\mathcal{N}}.
\]
\end{lem}

\begin{proof}
Let $\mathcal{N} \subset \bg$ be the nilpotent cone, i.e.~the closed subvariety whose $\K$-points are the nilpotent elements in $\bg$. (This scheme is denoted $\mathscr{N}_{\bG}^{\mathrm{var}}$ in~\cite[\S 4.3]{cotner}.)
By~\cite[Theorem~5.1]{cotner} there exists a $\bG$-equivariant ``Springer isomorphism''
\[
\mathcal{U} \simto \mathcal{N}.
\]
(Note that the unipotent and nilpotent schemes involved in~\cite[Theorem~5.1]{cotner} might not be reduced; but under our present assumptions they are by~\cite[Theorem~4.6 and~Theorem~4.12]{cotner}, so that they coincide with our varieties $\mathcal{U}$ and $\mathcal{N}$.)
By~\cite[Remark~10]{mcninch}, such an isomorphism necessarily restricts to a $\bB$-equivariant isomorphism $\bU \simto \bu$; inducing from $\bB$ to $\bG$ we deduce the desired isomorphism.
\end{proof}

Using Lemma~\ref{lem:U-sc}\eqref{it:functions-tU} and Lemma~\ref{lem:klt} we deduce the following properties.

\begin{cor}
 \phantomsection
\label{cor:properties-Spr}
 \begin{enumerate}
  \item 
  \label{it:global-sec-Spr-1}
  For any $i>0$, we have $\mathsf{H}^i(\widetilde{\mathcal{N}}, \mathscr{O}_{\widetilde{\mathcal{N}}})=0$.
  \item 
  \label{it:global-sec-Spr-2}
  We have an isomorphism of $\bG$-modules
  \[
   \mathscr{O}(\widetilde{\mathcal{N}}) \cong \mathrm{Ind}_{\rmZ_{\bG}(u)}^{\bG}(\K),
  \]
and moreover this module admits a good filtration.
 \end{enumerate}
\end{cor}

There exists a natural morphism $\Groth \to \bt$ induced by the quotient morphism $\bb \to \bb/\bu \cong \bt$, and $\widetilde{\mathcal{N}}$ identifies with the fiber of this map over $0$. We will also denote by $\bt_{\mathrm{reg}} \subset \bt$ the open subset of regular elements, i.e.~the complement of the union of the kernels of the differentials of the roots of $(\bG,\bT)$, see~\cite[Equation~(1) in~\S 13.3]{jantzen-nilp}. (Note that our assumption~\eqref{it:ass-3} implies that the differentials of roots are nonzero, so that this open subset is nonempty.)

\begin{lem}
\phantomsection
\label{lem:Groth-trs}
\begin{enumerate}
\item
\label{it:Groth-trs}
There exists a canonical isomorphism of $\bG$-varieties
\[
\bG/\bT \times \bt_{\mathrm{reg}} \simto \Groth \times_{\bt} \bt_{\mathrm{reg}}.
\]
\item
\label{it:Groth-rs}
We have a canonical isomorphism of $\bG$- and $\mathscr{O}(\bt_{\mathrm{reg}})$-modules
\[
\mathscr{O}(\bt_{\mathrm{reg}}) \otimes_{\mathscr{O}(\bt)} \mathscr{O}(\Groth) \cong \mathscr{O}(\bt_{\mathrm{reg}}) \otimes_\K \mathrm{Ind}_{\bT}^{\bG}(\K).
\]
\end{enumerate}
\end{lem}

\begin{proof}
\eqref{it:Groth-trs}
As explained in~\cite[p.~188]{jantzen-nilp}, the adjoint action induces an isomorphism
\[
\bU \times \bt_{\mathrm{reg}} \simto \bb \times_{\bt} \bt_{\mathrm{reg}},
\]
for the morphism $\bb \to \bt$ considered above. Inducing from $\bB$ to $\bG$ we deduce the desired isomorphism.

\eqref{it:Groth-rs} The claim follows from~\eqref{it:Groth-trs} by taking global sections.
\end{proof}

\begin{prop}
\label{prop:Vu-V0}
For any finite-dimensional $\bG$-module $V$ which admits a good filtration and any regular unipotent element $u \in \bG$, we have
\[
\dim(V^{\rmZ_{\bG}(u)})=\dim(V^{\bT}).
\]
\end{prop}

\begin{proof}
Corollary~\ref{cor:properties-Spr}\eqref{it:global-sec-Spr-1} and standard arguments imply that the $\bG$-module $\mathscr{O}(\Groth)$ admits a filtration with associated graded isomorphic to $\mathscr{O}(\bt) \otimes \mathscr{O}(\widetilde{\mathcal{N}})$; see~\cite[Proposition~3.4.1]{bmr} or~\cite[Lemma~4.12]{mr1} for similar considerations. 
Since the higher derived functors of the functor of $\bG$-invariants vanish on modules admitting a good filtration (see~\cite[Proposition~II.4.13]{jantzen}), and since $V \otimes \mathscr{O}(\widetilde{\mathcal{N}})$ admits a good filtration (by~\cite[Proposition~II.4.21]{jantzen}), it follows that $(V \otimes \mathscr{O}(\Groth))^{\bG}$ admits a filtration with associated graded
\[
\mathscr{O}(\bt) \otimes \bigl( V \otimes \mathscr{O}(\widetilde{\mathcal{N}}) \bigr)^{\bG}.
\]
In particular this module is free as an $\mathscr{O}(\bt)$-module, and by Corollary~\ref{cor:properties-Spr}\eqref{it:global-sec-Spr-2}, Frobenius reciprocity and the tensor identity, its rank is
\[
\dim \bigl( (V \otimes \mathscr{O}(\widetilde{\mathcal{N}}) )^{\bG} \bigr) = \dim(V^{\rmZ_{\bG}(u)}).
\]

%

On the other hand, using Lemma~\ref{lem:Groth-trs}\eqref{it:Groth-rs} we see that
\[
\mathscr{O}(\bt_{\mathrm{reg}}) \otimes_{\mathscr{O}(\bt)} (V \otimes \mathscr{O}(\Groth))^{\bG} \cong \bigl( V \otimes \mathscr{O}(\bt_{\mathrm{reg}}) \otimes \mathrm{Ind}_{\bT}^{\bG}(\K) \bigr)^{\bG}
\cong V^{\bT} \otimes \mathscr{O}(\bt_{\mathrm{reg}}),
\]
so that this $\mathscr{O}(\bt_{\mathrm{reg}})$-module is free of rank $\dim(V^{\bT})$. The desired equality follows.
\end{proof}

\subsection{Fixed points of regular unipotent elements in tilting modules with quasi-minuscule highest weight}
\label{ss:fixed-points-qmin}

In this subsection we drop our previous assumptions, but assume instead that $\bG$ is quasi-simple and simply connected, not of type $\mathbf{A}$. We will denote by $\alpha_0$ its highest short root (i.e.~its unique quasi-minuscule dominant weight).

\begin{figure}
 \begin{tabular}{|c|c|c|c|c|c|c|c|}
  \hline
  $\mathbf{B}_n$ ($n \geq 2$) & $\mathbf{C}_n$ ($n \geq 3$) & $\mathbf{D}_n$ ($n \geq 4$) & $\mathbf{E}_6$ & $\mathbf{E}_7$ & $\mathbf{E}_8$ & $\mathbf{F}_4$ & $\mathbf{G}_2$ \\
  \hline
  $\ell \neq 2$ & $\ell \nmid n$ & $\ell \neq 2$ & $\ell \neq 3$ & $\ell \neq 2$ & $\varnothing$ & $\ell \neq 3$ & $\ell \neq 2$\\
  \hline
 \end{tabular}
\caption{Conditions on $\ell$}
\label{fig:conditions}
\end{figure}


\begin{lem}
\label{lem:fixed-points-qmin}
 Assume that 
 $\ell$ satisfies the conditions in Figure~\ref{fig:conditions}. Then if $V$ is the indecomposable tilting module of highest weight $\alpha_0$ and $u$ is as in~\S\ref{ss:centralizer-ureg}
 we have
 \[
  \dim(V^u)=\dim(V^{\bT}).
 \]
\end{lem}

\begin{proof}\footnote{This proof was explained to one of us by M. Korhonen. An earlier proof of ours required stronger assumptions on $\ell$.}
 The conditions in Figure~\ref{fig:conditions} are exactly those which guarantee that the Weyl module of highest weight $\alpha_0$ is simple, hence tilting (see~\cite{lubeck}, see also~\cite[Table~1]{korhonen}).
 
 Assume first that $\bG$ is of exceptional type. Then the Jordan blocks of $u$ acting on $V$ (in particular their number, which coincides with $\dim(V^u)$) are described in~\cite{lawther}. From these tables we obtain that the number of Jordan blocks is in each case equal to the number of short simple roots in $\bR$, i.e.~to $\dim(V^{\bT})$.
 
 We assume now that $\bG$ is of classical type. If $\bG$ is of type $\mathbf{B}_n$ then $\bG \cong \mathrm{Spin}(2n+1,\K)$ and $V$ is the natural $(2n+1)$-dimensional representation of $\mathrm{SO}(2n+1,\K)$, seen as a $\bG$-module through the canonical surjection $\mathrm{Spin}(2n+1,\K) \twoheadrightarrow \mathrm{SO}(2n+1,\K)$; the desired claim is then clear. If $\bG$ is of type $\mathbf{D}_n$ then $V$ is the adjoint representation $\bg$. Since $\rmZ_{\bG}(u)$ is smooth (see Proposition~\ref{prop:Z-ureg}) we have
 \[
 \dim(\bg^u) = \dim(\rmZ_{\bG}(u)),
 \]  
 and since $u$ is regular this dimension is the rank of $\bG$,
 which of course coincides with $\dim(\bg^{\bT})$.
 
 Finally, we consider the case when $\bG$ is of type $\mathbf{C}_n$, so that $\bG \cong \mathrm{Sp}(2n,\K)$. If $E$ is the natural $2n$-dimensional representation of this group, then $V$ is the direct summand of $\wedge^2 E$ (as a $\bG$-module) given by the kernel of the natural morphism $\wedge^2 E \to \K$ sending $v \otimes v'$ to $\chi(v,v')$, where $\chi$ is the alternating form used to define the symplectic group (see e.g.~\cite{ps}); from this we see that to conclude it suffices to prove that $\dim((\wedge^2 E)^u)=n$. Now $u$ acting on $E$ has a single Jordan block, so that the desired claim follows from~\cite{lindsey} or~\cite[Theorem~2.1(2)]{barry}.
\end{proof}

\begin{rmk}
M.~Korhonen informed us that Lemma~\ref{lem:fixed-points-qmin} holds for all characteristics, except maybe in type $\mathbf{D}$ for $\ell=2$. For instance in type $\mathbf{B}$ for $\ell=2$, the tilting module is the restriction of the natural module for the group of type $\mathbf{D}_{n+1}$, and a regular unipotent element for the group of type $\mathbf{B}_n$ is regular in the group of type $\mathbf{D}_{n+1}$, hence acts on this module with two Jordan blocks.
\end{rmk}

\section{Reconstructing a subgroup from the restriction functor}
\label{sec:reconstructing}

In this section we prove a slight extension of~\cite[Proposition~1]{bez} which will be required for our arguments.

\subsection{Statement}
\label{ss:statement-reconstruction}

Let $(\mathsf{A}, \star)$ be a monoidal category.
Recall that the \emph{Drinfeld center} 
$\mathrm{DrC}(\mathsf{A})$ 
of $(\mathsf{A}, \star)$ is the category whose objects are pairs $(X,\iota)$ with $X$ an object in $\mathsf{A}$ and
\[
 \iota : X \star (-) \simto (-) \star X
\]
an isomorphism of functors such that
\[
 \iota_{Y \star Z} = (\id \star \iota_Z) \circ (\iota_Y \star \id)
\]
for any $Y,Z$ in $\mathsf{A}$ (where we omit the associativity isomorphism for simplicity), and whose morphisms from $(X,\iota)$ to $(X',\iota')$ are morphisms $X \to X'$ in $\mathsf{A}$ which are compatible with the morphisms $\iota$ and $\iota'$ in the obvious way. The Drinfeld center of a monoidal category admits a natural structure of braided monoidal category, with monoidal product induced by that of $\mathsf{A}$, and the braiding $(X,\iota) \star (X',\iota') \simto (X',\iota') \star (X,\iota)$ given by $\iota$.

If $(\mathsf{B},\otimes)$ is a \emph{symmetric} monoidal category, then as explained in~\cite[\S 2.1]{bez} the datum of a braided monoidal functor from $\mathsf{B}$ to $\mathrm{DrC}(\mathsf{A})$ is equivalent to that of a \emph{central functor} from $\mathsf{B}$ to $\mathsf{A}$, i.e.~a pair consisting of a monoidal functor $F : \mathsf{B} \to \mathsf{A}$ (with ``monoidality'' isomorphisms $\phi_{-,-}$) together with bifunctorial isomorphisms
\[
 \sigma_{X,Y} : F(X) \star Y \simto Y \star F(X)
\]
for $X$ in $\mathsf{B}$ and $Y$ in $\mathsf{A}$, which satisfy the following properties.
\begin{enumerate}
 \item 
 \label{it:central-functor-1}
 For $X,X'$ in $\mathsf{B}$, the isomorphism $\sigma_{X,F(X')}$ coincides with the composition
 \[
  F(X) \star F(X') \xrightarrow[\sim]{\phi_{X,X'}^{-1}} F(X \otimes X') \simto F(X' \otimes X) \xrightarrow[\sim]{\phi_{X',X}} F(X') \star F(X)
 \]
where the middle isomorphism is the image under $F$ of the commutativity constraint of $\mathsf{B}$ (applied to $(X,X')$).
\item
\label{it:central-functor-2}
For $Y_1,Y_2$ in $\mathsf{A}$ and $X$ in $\mathsf{B}$, the composition
\[
F(X) \star Y_1 \star Y_2 \xrightarrow[\sim]{\sigma_{X,Y_1} \star \id_{Y_2}} Y_1 \star F(X) \star Y_2 \xrightarrow[\sim]{\id_{Y_1} \star \sigma_{X,Y_2}} Y_1 \star Y_2 \star F(X)
 \]
 coincides with $\sigma_{X,Y_1 \star Y_2}$ (where we omit the associativity constraint of $\mathsf{A}$).
 \item
 \label{it:central-functor-3}
 For $Y$ in $\mathsf{A}$ and $X_1, X_2$ in $\mathsf{B}$, the composition
 \begin{multline*}
  F(X_1 \otimes X_2) \star Y \xrightarrow[\sim]{\phi_{X_1,X_2}\star Y} F(X_1) \star F(X_2) \star Y \xrightarrow[\sim]{\id_{F(X_1)} \star \sigma_{X_2,Y}} F(X_1) \star Y \star F(X_2) \\
  \xrightarrow[\sim]{\sigma_{X_1,Y} \star \id_{F(X_2)}} Y \star F(X_1) \star F(X_2) \xrightarrow[\sim]{\id_Y \star \phi_{X_1,X_2}^{-1}} Y \star F(X_1 \otimes X_2)
 \end{multline*}
 coincides with $\sigma_{X_1 \otimes X_2,Y}$.
\end{enumerate}

Our goal in this section is to prove the following statement. Here, we let $\mathbb{K}$ be an algebraically closed field. For $H$ an affine $\K$-group scheme of finite type we denote by $\Rep(H)$ the category of finite-dimensional $H$-modules (in other words, finite dimensional $\scO(H)$-comodules) and, for $H' \subset H$ a subgroup scheme we denote by
\[
\mathrm{Res}^{H}_{H'} : \Rep(H) \to \Rep(H')
\]
the restriction functor. The category $\Rep(H)$ admits a natural symmetric monoidal structure, with product given by tensor product over $\K$.

\begin{prop}
\label{prop:subgroup}
 Let $H$ be an affine $\K$-group scheme of finite type.
 Let $(\mathsf{A},\otimes)$ be a $\mathbb{K}$-linear abelian monoidal category with unit object $\mathbbm{1}$ which satisfies $\End(\mathbbm{1})=\mathbb{K}$, and such that $\otimes$ is exact in each variable. We assume furthermore that $\mathbbm{1}$ is a simple object of $\mathsf{A}$, and that $\Hom_{\mathsf{A}}(\mathbbm{1},X)$ is finite-dimensional for all $X$ in $\mathsf{A}$.
 Assume we are given an exact, central functor
 \[
  F : \Rep(H) \to \mathsf{A}
 \]
such that each object in $\mathsf{A}$ is isomorphic to a subquotient of an object $F(V)$ with $V$ in $\Rep(H)$. Then there exists a subgroup scheme $H' \subset H$ and an equivalence of monoidal categories
\[
 \Phi : \Rep(H') \simto \mathsf{A}
\]
which satisfies $F \cong \Phi \circ \mathrm{Res}^{H}_{H'}$ as monoidal functor.
\end{prop}

\begin{rmk}
\begin{enumerate}
\item
The assumption that $\mathbbm{1}$ is simple follows from the condition that $\End(\mathbbm{1})=\mathbb{K}$ in case $\mathsf{A}$ is rigid, see~\cite[Proposition~1.17]{dm}.
\item
The subgroup $H' \subset H$ is unique only up to conjugation; using the notation introduced in~\S\ref{ss:proof-prop-subgroup} below,
it depends on the subobject $\mathcal{J}$, which cannot be chosen in any canonical way.
\end{enumerate}
\end{rmk}

This proposition appears as~\cite[Proposition~1]{bez} with slightly different assumptions. (Namely, it is not assumed in~\cite{bez} that $\mathbbm{1}$ is simple, but one imposes the extra assumption that $\mathbb{K}$ is uncountable.) Below we show how to adapt the proof in~\cite{bez} in order to fit with the present assumptions. (In later sections we will apply this proposition in case $\mathbb{K}$ is an algebraic closure of a finite field, hence in particular \emph{is} countable.) No detail of this proof will be used in later sections.

\subsection{Preliminaries}

The main ingredient needed to extend the proof of~\cite{bez} will be the following easy commutative algebra lemma.

\begin{lem}
\label{lem:comm-alg}
 Let $\mathbb{K}$ be an algebraically closed field, and let $\mathbb{L}/\mathbb{K}$ be a field extension. Assume there exists a field extension $\mathbb{K'}/\mathbb{K}$, a finitely generated commutative $\mathbb{K}'$-algebra $A$, and an embedding of $\mathbb{K}'$-algebras $\mathbb{K}' \otimes_{\mathbb{K}} \mathbb{L} \subset A$. Then $\mathbb{L}=\mathbb{K}$.
\end{lem}

\begin{proof}
One can of course assume that $\mathbb{K}'$ is algebraically closed.
 Assume for a contradiction that $\mathbb{L} \neq \mathbb{K}$, and let $a \in \mathbb{L} \smallsetminus \mathbb{K}$. Let $\tilde{a}=1 \otimes a \in \mathbb{K}' \otimes_{\mathbb{K}} \mathbb{L} \subset A$. One can consider $\tilde{a}$ as a scheme morphism $\mathrm{Spec}(A) \to \mathbb{A}^1_{\mathbb{K}'}$. The image of this morphism is constructible by Chevalley's theorem (see~\cite[\href{https://stacks.math.columbia.edu/tag/00FE}{Tag 00FE}]{stacks-project}), hence is either a finite subset of $\mathbb{K}'$ or the complement of such a subset. Now $\tilde{a}-t$ is invertible in $A$ for any $t \in \mathbb{K}$, so that this image does not intersect $\mathbb{K}$. We deduce that it is finite, so that there exist $a_1, \cdots, a_r \in \mathbb{K}'$ such that the morphism $\mathrm{Spec}(A) \to \mathbb{A}^1_{\mathbb{K}'}$ associated with the element $\prod_{i=1}^r (\tilde{a}-a_i)$ has image $\{0\}$, i.e.~is nilpotent in $A$. Hence for some $N \geq 1$ the elements $1,\tilde{a},\cdots, \tilde{a}^N$ are linearly dependent (over $\mathbb{K}'$) in $\mathbb{K}' \otimes_{\mathbb{K}} \mathbb{L}$. Then the elements $1,a,\cdots, a^N$ must be linearly dependent (over $\mathbb{K}$) in $\mathbb{L}$, so that $a$ is algebraic over $\mathbb{K}$. This is impossible since $\mathbb{K}$ is algebraically closed and $a \notin \mathbb{K}$.
\end{proof}

We will also require the following lemma.

\begin{lem}
\label{lem:monoidal-cat-embedding}
Let $\K$ be a field, and let $(\mathsf{A},\otimes)$ be a $\mathbb{K}$-linear abelian monoidal category which satisfies the following assumptions:
\begin{enumerate}
\item
the unit object $\mathbbm{1}$ is simple and satisfies $\End(\mathbbm{1})=\mathbb{K}$;
\item
for any $X$ in $\mathsf{A}$, the $\K$-vector space
$\Hom_{\mathsf{A}}(\mathbbm{1},X)$ is finite-dimensional;
\item
$\otimes$ is exact in each variable.
\end{enumerate}
Then for all objects $M$ and $N$ in $\mathsf{A}$, the morphism
\[
\Hom_{\mathsf{A}}(\mathbbm{1}, M)\otimes_{\mathbb{K}}\Hom_{\mathsf{A}}(\mathbbm{1}, N) \to \Hom_{\mathsf{A}}(\mathbbm{1}, M\otimes N)
\]
which, for $f \in \Hom_{\mathsf{A}}(\mathbbm{1}, M)$ and $g \in \Hom_{\mathsf{A}}(\mathbbm{1}, N)$, sends $f \otimes g$ to the morphism
\[
\mathbbm{1} = \mathbbm{1} \otimes \mathbbm{1} \xrightarrow{f \otimes g} M \otimes N,
\]
is injective.
\end{lem}

\begin{proof}
For any $M$ in $\mathsf{A}$, since $\dim \Hom_{\mathsf{A}}(\mathbbm{1}, M) < \infty$ we can consider the object $\Hom_{\mathsf{A}}(\mathbbm{1}, M)\otimes_{\mathbb{K}} \mathbbm{1}$, and we have an obvious morphism
\begin{equation}
\label{eqn:morph-unit}
\Hom_{\mathsf{A}}(\mathbbm{1}, M)\otimes_{\mathbb{K}} \mathbbm{1} \to M.
\end{equation}
We claim that this morphism
is injective. In fact, its kernel is a subobject of the semisimple object $\Hom_{\mathsf{A}}(\mathbbm{1}, M)\otimes_{\mathbb{K}} \mathbbm{1}$, hence is itself a direct sum of copies of $\mathbbm{1}$. If it were nonzero, then there would exist a morphism $\mathbbm{1} \to \Hom_{\mathsf{A}}(\mathbbm{1}, M)\otimes_{\mathbb{K}} \mathbbm{1}$ whose composition with~\eqref{eqn:morph-unit} vanishes. Now we have
\[
\Hom_{\mathsf{A}} \bigl( \mathbbm{1}, \Hom_{\mathsf{A}}(\mathbbm{1}, M)\otimes_{\mathbb{K}} \mathbbm{1} \bigr) \cong \Hom_{\mathsf{A}}(\mathbbm{1}, M),
\]
and for any nonzero $f$ in the right-hand side, the composition of its image in the left-hand side with~\eqref{eqn:morph-unit} is $f$, hence is nonzero. The kernel under consideration must therefore vanish.

Now, consider a second object $N$ in $\mathsf{A}$.
Tensoring~\eqref{eqn:morph-unit} on the right with $N$ (which is an exact functor by assumption) we deduce an embedding
\[
\Hom_{\mathsf{A}}(\mathbbm{1}, M)\otimes_{\mathbb{K}} N \hookrightarrow M \otimes N.
\]
Finally, applying the left-exact functor $\Hom_{\mathsf{A}}(\mathbbm{1},-)$ we deduce an embedding
\[
\Hom_{\mathsf{A}}(\mathbbm{1}, M)\otimes_{\mathbb{K}}\Hom_{\mathsf{A}}(\mathbbm{1}, N) \hookrightarrow \Hom_{\mathsf{A}}(\mathbbm{1}, M\otimes N),
\]
which is easily seen to be given by the formula in the statement.
\end{proof}

\subsection{Proof of Proposition~\ref{prop:subgroup}}
\label{ss:proof-prop-subgroup}

In this subsection we assume that the conditions in Proposition~\ref{prop:subgroup} are satisfied, and we explain how to prove~\cite[Lemma~1]{bez} without assuming that $\mathbb{K}$ is uncountable.
Once that is in place, the argument for Proposition~\ref{prop:subgroup} follows in the same manner as for~\cite[Proposition~1]{bez}. 

We recall the necessary notation from \emph{loc}.~\emph{cit}. Let $\underline{\mathscr{O}}_{H}$ denote the ring of regular functions on $H$ regarded as a commutative ring object in the symmetric monoidal category $\Ind(\Rep(H))$ of ind-objects\footnote{For generalities on categories of ind-objects we refer to~\cite[Chap.~6]{ks}. We will use in particular the fact that the category of ind-objects in an abelian category is abelian, see~\cite[Theorem~8.6.5(i)]{ks}, and that extensions of exact functors to ind-objects are exact, see~\cite[Corollary~8.6.8]{ks}.} in $\Rep(H)$ with $H$-action coming from left multiplication of $H$ on itself. Then, still denoting by $F$ the natural extension of this functor to ind-objects, $F(\underline{\mathscr{O}}_{H})$ inherits a structure of ring object in the category $\Ind(\mathsf{A})$ of ind-objects in $\mathsf{A}$. Let $\mathcal{J}\subset F(\underline{\mathscr{O}}_{H})$ be a maximal left ideal subobject; then as explained in~\cite[p.~73]{bez}, $\mathcal{J}$ is automatically a right ideal subobject also, so that $\underline{\mathscr{O}}_{H'} := F(\underline{\mathscr{O}}_{H})/\mathcal{J}$ acquires a natural ring object (in $\Ind(\mathsf{A})$) structure. (For details justifying the existence of a maximal left ideal subobject, see~\cite{ar-book}.) The multiplication map for this algebra will be denoted $m_{\underline{\mathscr{O}}_{H'}}: \underline{\mathscr{O}}_{H'} \otimes \underline{\mathscr{O}}_{H'} \rightarrow \underline{\mathscr{O}}_{H'}$, and the unit (inherited from the unital structure of $\underline{\mathscr{O}}_{H}$) will be denoted $\iota: \mathbbm{1}\to \underline{\mathscr{O}}_{H'}$. 

We may now consider the category of (left) $\underline{\mathscr{O}}_{H'}$-module ind-objects in $\mathsf{A}$; we denote morphisms in this category by $\Hom_{\underline{\mathscr{O}}_{H'}}(-,-)$. Since $\underline{\mathscr{O}}_{H'}$ is simple when regarded as an $\underline{\mathscr{O}}_{H'}$-module, the $\K$-algebra $\mathbb{L} := \End_{\underline{\mathscr{O}}_{H'}}(\underline{\mathscr{O}}_{H'} )$ is a division algebra. The content of~\cite[Lemma~1]{bez} is then the following.

\begin{lem}
\label{lem:L=K}
We have $\mathbb{L}=\K$.
\end{lem}

The proof of this lemma will be given at the end of the subsection, after some preliminary results that we now consider.

\begin{lem}
\label{lem:varsigma}
There exists a canonical isomorphism
\[
\varsigma :
\underline{\mathscr{O}}_{H'} \otimes \underline{\mathscr{O}}_{H'} \simto \underline{\mathscr{O}}_{H'} \otimes \underline{\mathscr{O}}_{H'}
\]
which satisfies the following equalities:
\begin{enumerate}
\item
\label{it:prop-varsigma-1}
for any morphism $f : \mathbbm{1} \to \underline{\mathscr{O}}_{H'}$ we have
$\varsigma \circ (\id_{\underline{\mathscr{O}}_{H'}} \otimes f) = f \otimes \id_{\underline{\mathscr{O}}_{H'}}$
(where we identify $\mathbbm{1} \otimes  \underline{\mathscr{O}}_{H'}$ and $ \underline{\mathscr{O}}_{H'} \otimes \mathbbm{1}$ with $ \underline{\mathscr{O}}_{H'}$ using the unit constraint);
\item
\label{it:prop-varsigma-2}
$m_{\underline{\mathscr{O}}_{H'}} \circ \varsigma = m_{\underline{\mathscr{O}}_{H'}}$;
\item
\label{it:prop-varsigma-3}
$(\id \otimes m_{\underline{\mathscr{O}}_{H'}}) \circ (\varsigma \otimes \id) \circ (\id \otimes \varsigma) = \varsigma \circ (m_{\underline{\mathscr{O}}_{H'}} \otimes \id)$;
\item
\label{it:prop-varsigma-4}
$(m_{\underline{\mathscr{O}}_{H'}} \otimes \id) \circ (\id \otimes \varsigma) \circ (\varsigma  \otimes \id) = \varsigma \circ (\id \otimes m_{\underline{\mathscr{O}}_{H'}})$.
\end{enumerate}
\end{lem}

\begin{proof}
As part of the central structure on $F$, we have a canonical isomorphism
\[
\sigma_{\underline{\mathscr{O}}_{H},F(\underline{\mathscr{O}}_{H})} : F(\underline{\mathscr{O}}_{H}) \otimes F(\underline{\mathscr{O}}_{H}) \simto F(\underline{\mathscr{O}}_{H}) \otimes F(\underline{\mathscr{O}}_{H}),
\]
which by functoriality this isomorphism sends $F(\underline{\mathscr{O}}_{H}) \otimes \mathcal{J}$ to $\mathcal{J} \otimes F(\underline{\mathscr{O}}_{H})$. Moreover, this isomorphism is also (via the appropriate identifications) the image under $F$ of the commutativity constraint in $\Rep(H)$, which is symmetric; therefore it satisfies $\sigma_{\underline{\mathscr{O}}_{H},F(\underline{\mathscr{O}}_{H})} \circ \sigma_{\underline{\mathscr{O}}_{H},F(\underline{\mathscr{O}}_{H})}=\id$, hence also sends $\mathcal{J} \otimes F(\underline{\mathscr{O}}_{H})$ to $F(\underline{\mathscr{O}}_{H}) \otimes \mathcal{J}$. It follows that this isomorphism sends $\mathcal{J} \otimes F(\underline{\mathscr{O}}_{H}) + F(\underline{\mathscr{O}}_{H}) \otimes \mathcal{J}$ to itself, so that it induces the wished-for isomorphism
\[
\varsigma :
\underline{\mathscr{O}}_{H'} \otimes \underline{\mathscr{O}}_{H'} \simto \underline{\mathscr{O}}_{H'} \otimes \underline{\mathscr{O}}_{H'}.
\]

All the desired equalities hold for the corresponding structures on $\underline{\mathscr{O}}_{H}$ in the category $\Ind(\Rep(H))$; hence they hold for $F(\underline{\mathscr{O}}_{H})$ in $\Ind(\mathsf{A})$, and finally for its quotient $\underline{\mathscr{O}}_{H'}$.
%
\end{proof}

\begin{cor}
The division algebra $\mathbb{L}$ is commutative, hence a field.
\end{cor}

\begin{proof}
Composition with $\iota$ defines an isomorphism
\[
\mathbb{L} \simto \Hom_{\Ind(\mathsf{A})}(\mathbbm{1}, \underline{\mathscr{O}}_{H'}),
\] 
with inverse given by $f\mapsto m_{\underline{\mathscr{O}}_{H'}}\circ(\id_{\underline{\mathscr{O}}_{H'}} \otimes f)$.
In these terms multiplication can be described as follows: given $f,g : \mathbbm{1} \to \underline{\mathscr{O}}_{H'}$, the product $f \cdot g$ is the composition
\[
\mathbbm{1} = \mathbbm{1} \otimes \mathbbm{1} \xrightarrow{f \otimes g} \underline{\mathscr{O}}_{H'} \otimes \underline{\mathscr{O}}_{H'} \xrightarrow{m_{\underline{\mathscr{O}}_{H'}}} \underline{\mathscr{O}}_{H'}.
\]
Using the isomorphism $\varsigma_{\underline{\mathscr{O}}_{H'}}$ and properties~\eqref{it:prop-varsigma-1}--\eqref{it:prop-varsigma-2} in Lemma~\ref{lem:varsigma}, from this description
it is not difficult to check that $\mathbb{L}$ is commutative, hence a field.
\end{proof}

\begin{lem}
\phantomsection
\label{lem:K=L-prelim}
\begin{enumerate}
\item
\label{it:K=L-prelim-2}
For any $V$ in $\Rep(H)$ we have a canonical isomorphism of $\underline{\mathscr{O}}_{H'}$-modules
\[
\underline{\mathscr{O}}_{H'} \otimes F(V) \cong \underline{\mathscr{O}}_{H'} \otimes_{\K} V.
\]
\item
\label{it:K=L-prelim-3}
For any $M$ in $\mathsf{A}$, the $\underline{\mathscr{O}}_{H'}$-module $\underline{\mathscr{O}}_{H'} \otimes M$ is isomorphic to a finite direct sum of copies of $\underline{\mathscr{O}}_{H'}$.
\item
\label{it:K=L-prelim-4}
The functor $M \mapsto \Hom_{\Ind(\mathsf{A})}(\mathbbm{1},\underline{\mathscr{O}}_{H'} \otimes M)$ is exact.
\end{enumerate}
\end{lem}

\begin{proof}
\eqref{it:K=L-prelim-2}
We have a canonical isomorphism of $H$- and $\underline{\mathscr{O}}_H$-modules
\[
\underline{\mathscr{O}}_H \otimes V \simto \underline{\mathscr{O}}_H \otimes_\K \For^{H}(V),
\]
where in the left-hand side we consider the diagonal $H$-action and on the right-hand side we consider the action on $\underline{\mathscr{O}}_H$ only. 
Applying $F$ we deduce an isomorphism of $F(\underline{\mathscr{O}}_H)$-modules
\[
F(\underline{\mathscr{O}}_H) \otimes F(V) \simto F(\underline{\mathscr{O}}_H) \otimes_\K V.
\]
Being an isomorphism of $F(\underline{\mathscr{O}}_H)$-modules, this morphism must send
\[
\mathcal{J} \cdot \bigl( F(\underline{\mathscr{O}}_H) \otimes F(V) \bigr) = \mathcal{J} \otimes F(V)
\]
to
\[
\mathcal{J} \cdot \bigl( F(\underline{\mathscr{O}}_H) \otimes_\K V \bigr) = \mathcal{J} \otimes_\K V;
\]
we deduce the desired isomorphism by passing to the quotients by these submodules.


\eqref{it:K=L-prelim-3}
By assumption, $M$ is a subquotient of $F(V)$ for some $V$ in $\Rep(H)$. Then, in view of~\eqref{it:K=L-prelim-2} the $\underline{\mathscr{O}}_{H'}$-module $\underline{\mathscr{O}}_{H'} \otimes M$ is a subquotient of $\underline{\mathscr{O}}_{H'} \otimes_\K V$. The latter object being a direct sum of copies of the simple module $\underline{\mathscr{O}}_{H'}$, $\underline{\mathscr{O}}_{H'} \otimes M$ must also be a direct sum of copies of this module.

\eqref{it:K=L-prelim-4}
Consider an exact sequence $0 \to M_1 \to M_2 \to M_3 \to 0$ in $\mathsf{A}$. By exactness of $\otimes$, tensoring with $\underline{\mathscr{O}}_{H'}$ we deduce an exact sequence
\begin{equation}
\label{eqn:ses-OH}
0 \to \underline{\mathscr{O}}_{H'} \otimes M_1 \to \underline{\mathscr{O}}_{H'} \otimes M_2 \to \underline{\mathscr{O}}_{H'} \otimes M_3 \to 0
\end{equation}
in $\Ind(\mathsf{A})$. Then, applying the left-exact functor $\Hom_{\Ind(\mathsf{A})}(\mathbbm{1},-)$ we deduce an exact sequence
\begin{multline*}
0 \to \Hom_{\Ind(\mathsf{A})}(\mathbbm{1},\underline{\mathscr{O}}_{H'} \otimes M_1) \to \Hom_{\Ind(\mathsf{A})}(\mathbbm{1},\underline{\mathscr{O}}_{H'} \otimes M_2) \\
\to \Hom_{\Ind(\mathsf{A})}(\mathbbm{1},\underline{\mathscr{O}}_{H'} \otimes M_3).
\end{multline*}
Now, each $\Hom_{\Ind(\mathsf{A})}(\mathbbm{1},\underline{\mathscr{O}}_{H'} \otimes M_i)$ has a natural structure of $\mathbb{L}$-vector space, and each map in this sequence is $\mathbb{L}$-linear. The dimension of $\Hom_{\Ind(\mathsf{A})}(\mathbbm{1},\underline{\mathscr{O}}_{H'} \otimes M_i)$ is the number of copies of $\underline{\mathscr{O}}_{H'}$ appearing in $\underline{\mathscr{O}}_{H'} \otimes M_i$. The exact sequence~\eqref{eqn:ses-OH} shows that the number of copies in $\underline{\mathscr{O}}_{H'} \otimes M_2$ is the sum of the number of copies in $\underline{\mathscr{O}}_{H'} \otimes M_1$ and $\underline{\mathscr{O}}_{H'} \otimes M_3$. Therefore the morphism $\Hom_{\Ind(\mathsf{A})}(\mathbbm{1},\underline{\mathscr{O}}_{H'} \otimes M_2) \to \Hom_{\Ind(\mathsf{A})}(\mathbbm{1},\underline{\mathscr{O}}_{H'} \otimes M_3)$ must be surjective, showing that the above sequence is exact.
\end{proof}

\begin{proof}[Proof of Lemma~\ref{lem:L=K}]

%
We set
\[
R:=\Hom_{\Ind(\mathsf{A})}(\mathbbm{1}, \underline{\mathscr{O}}_{H'}\otimes\underline{\mathscr{O}}_{H'}).
\]
Properties~\eqref{it:prop-varsigma-3}--\eqref{it:prop-varsigma-4} in Lemma~\ref{lem:varsigma} ensure that the composition
\begin{multline*}
\bigl( \underline{\mathscr{O}}_{H'}\otimes\underline{\mathscr{O}}_{H'} \bigr) \otimes \bigl( \underline{\mathscr{O}}_{H'}\otimes\underline{\mathscr{O}}_{H'} \bigr) \xrightarrow{\id \otimes \varsigma \otimes \id} 
\underline{\mathscr{O}}_{H'}\otimes\underline{\mathscr{O}}_{H'} \otimes \underline{\mathscr{O}}_{H'} \otimes\underline{\mathscr{O}}_{H'} \\
\xrightarrow{m_{\underline{\mathscr{O}}_{H'}} \otimes m_{\underline{\mathscr{O}}_{H'}}}
\underline{\mathscr{O}}_{H'}\otimes\underline{\mathscr{O}}_{H'}
\end{multline*}
defines a structure of ring-object on $\underline{\mathscr{O}}_{H'}\otimes\underline{\mathscr{O}}_{H'}$; hence, copying the description of the product in $\Hom_{\Ind(\mathsf{A})}(\mathbbm{1}, \underline{\mathscr{O}}_{H'})$, we see that $R$ admits a natural ring structure. As for $\mathbb{L}$ one sees that this ring is commutative. 

The embeddings
\[
\underline{\mathscr{O}}_{H'} = \underline{\mathscr{O}}_{H'} \otimes \mathbbm{1} \xrightarrow{\id \otimes \iota} \underline{\mathscr{O}}_{H'}\otimes\underline{\mathscr{O}}_{H'}
\quad \text{and} \quad
\underline{\mathscr{O}}_{H'} = \mathbbm{1} \otimes \underline{\mathscr{O}}_{H'} \xrightarrow{\iota \otimes \id} \underline{\mathscr{O}}_{H'}\otimes\underline{\mathscr{O}}_{H'} 
\]
induce two $\K$-algebra morphisms $\mathbb{L} \to R$. We use the first one to turn $R$ into an $\mathbb{L}$-algebra. Considering also the second one we obtain an $\mathbb{L}$-linear morphism
\[
\mathbb{L}\otimes_{\mathbb{K}}\mathbb{L}\to \Hom_{\mathsf{A}}(\mathbbm{1}, \underline{\mathscr{O}}_{H'}\otimes\underline{\mathscr{O}}_{H'})
\]
which, for $f,g : \mathbbm{1} \to \underline{\mathscr{O}}_{H'}$, sends $f \otimes g$ to the morphism
\[
\mathbbm{1} = \mathbbm{1} \otimes \mathbbm{1} \xrightarrow{f \otimes g} \underline{\mathscr{O}}_{H'} \otimes \underline{\mathscr{O}}_{H'}.
\]
By Lemma~\ref{lem:monoidal-cat-embedding} this morphism is injective. It is also easily seen to be an $\mathbb{L}$-algebra morphism.

Now by exactness of $\otimes$ we have a surjective morphism
\[
\underline{\mathscr{O}}_{H'} \otimes F(\underline{\mathscr{O}}_{H}) \twoheadrightarrow \underline{\mathscr{O}}_{H'} \otimes \underline{\mathscr{O}}_{H'},
\]
which by Lemma~\ref{lem:K=L-prelim}\eqref{it:K=L-prelim-2} can be interpreted as a surjective morphism
\[
\underline{\mathscr{O}}_{H'} \otimes_\K \scO(H) \twoheadrightarrow \underline{\mathscr{O}}_{H'} \otimes \underline{\mathscr{O}}_{H'}.
\]
This morphism is an algebra-object morphism, where the structure on the left-hand side is the tensor product of the ring structures on $\underline{\mathscr{O}}_{H'}$ and $\scO(H)$, and the structure on the right-hand side is as above. Applying the functor $\Hom_{\Ind(\mathsf{A})}(\mathbbm{1},-)$, in view of Lemma~\ref{lem:K=L-prelim}\eqref{it:K=L-prelim-4} we deduce a surjective $\mathbb{L}$-algebra morphism
\[
\mathbb{L}\otimes_{\mathbb{K}} \scO(H) \twoheadrightarrow R,
\]
which shows that $R$ is finitely generated over $\mathbb{L}$.

%
%
We have now checked the assumptions of Lemma~\ref{lem:comm-alg} for $A =  R$ and $\mathbb{K}' = \mathbb{L}$; this lemma implies that $\mathbb{L}=\K$, which finishes the proof.
\end{proof}

\section{Geometric Satake equivalence, central sheaves and Wakimoto sheaves}
\label{sec:Satake}

In this section we provide a reminder on the geometric Satake equivalence and the construction of ``central sheaves'' on the affine flag variety.

\subsection{The affine Grassmannian and the geometric Satake equivalence}
\label{ss:Satake}

We let $G$ be a connected reductive algebraic group over an algebraically closed field $\F$ of characteristic $p>0$, and $\bk$ be an algebraic closure of a finite field of characteristic $\ell \neq p$. 

We will denote by $\Loop G$, resp.~$\Loop^+ G$, the functor from $\F$-algebras to groups defined by
\[
 R \mapsto G \bigl( R(\hspace{-1pt}(z)\hspace{-1pt}) \bigr), \quad \text{resp.} \quad G \bigl( R[\hspace{-1pt}[z]\hspace{-1pt}] \bigr),
\]
where $z$ is an indeterminate. It is well known that the functor $\Loop G$ is represented by an ind-affine group ind-scheme (which will be denoted similarly), and that $\Loop^+ G$ is represented by an affine group scheme (which will also be denoted similarly). One can then define the affine Grassmannian $\Gr_G$ as the fppf quotient
\[
 \Gr_G = (\Loop G/\Loop^+ G)_{\mathrm{fppf}}.
\]
It is well known also that $\Gr_G$ is represented by an ind-projective ind-scheme,
 which will once again be denoted similarly.

If $\bk'$ is a finite subfield of $\bk$, then we can consider the $\Loop^+ G$-equivariant derived category $\Db_{\Loop^+ G}(\Gr_G,\bk')$ of $\bk'$-sheaves on $\Gr_G$. (The definition of this category requires a little bit of care, but we will ignore these subtleties here; see~\cite[\S 6]{gaitsgory} or~\cite[\S 1.16.4]{br} for discussions on this topic. Similar comments apply for various categories of sheaves on ind-schemes considered below.) Taking the direct limit of these categories over all finite subfields of $\bk$ one obtains a category which will be denoted $\sfD_{\Loop^+G, \Loop^+G}$, and which we will consider informally as the $\Loop^+ G$-equivariant derived category of $\bk$-sheaves on $\Gr_G$. This category admits a natural convolution product $\star_{\Loop^+ G}$, and an associativity constraint $\varphi$ for this product. The skyscraper sheaf $\delta_{\Gr}$ at the base point of $\Gr_G$ is a unit for $\star$, so that we obtain a monoidal category $(\sfD_{\Loop^+G, \Loop^+G}, \star, \varphi, \delta_{\Gr})$.

On $\sfD_{\Loop^+G, \Loop^+G}$ we also have a natural perverse t-structure, whose heart will be denoted
$\sfP_{\Loop^+G, \Loop^+G}$. The bifunctor $\star_{\Loop^+G}$ restricts to an exact bifunctor
\[
 \sfP_{\Loop^+G, \Loop^+G} \times \sfP_{\Loop^+G, \Loop^+G} \to \sfP_{\Loop^+G, \Loop^+G},
\]
see~\cite[\S 1.6.3 and \S 1.10.3]{br} for details and references; moreover the restriction of $\star_{\Loop^+G}$ to the category $\sfP_{\Loop^+G, \Loop^+G}$ admits a commutativity constraint which can be constructed (following an idea of Drinfeld) using a description of this product as a ``fusion product;'' see~\cite[\S 1.7]{br} for details. In this way, $\sfP_{\Loop^+G, \Loop^+G}$ comes equipped with the structure of an abelian symmetric monoidal category.

Using a detailed study of the so-called ``weight functors,'' Mirkovi{\'c}--Vilonen construct in~\cite{mv} an affine $\bk$-group scheme $\widetilde{G}_\bk$ and an equivalence of abelian symmetric monoidal categories
\begin{equation}
\label{eqn:Satake-equiv}
\Sat : \sfP_{\Loop^+G, \Loop^+G} \simto \Rep(\widetilde{G}_\bk),
\end{equation}
where the right-hand side denotes the category of finite-dimensional $\widetilde{G}_\bk$-modules.
Following earlier work of Lusztig~\cite{lusztig}, Ginzburg~\cite{ginzburg} and Be{\u\i}linson--Drinfeld~\cite{bd}, they then prove the following.

\begin{thm}[Geometric Satake equivalence, \cite{mv}]
\label{thm:Satake}
 The group scheme $\widetilde{G}_\bk$ is a connected reductive group, and its root datum is dual to that of $G$.
\end{thm}

See~\cite{br} for a detailed study of the proof of Theorem~\ref{thm:Satake}. In view of this theorem, in the rest of the paper we will denote the reductive group $\widetilde{G}_\bk$ by $G^\vee_\bk$. In the course of the proof of Theorem~\ref{thm:Satake}, Mirkovi{\'c} and Vilonen construct a canonical maximal torus $T^\vee_\bk$ of $G^\vee_\bk$.
(More precisely, given a choice of Borel subgroup $B \subset G$ and of maximal torus $T \subset B$, one obtains ``weight functors'' which allow to define a maximal torus $T^\vee_\bk \subset G^\vee_\bk$ with a canonical identification $X^*(T^\vee_\bk) = X_*(T)$. The authors then prove that the resulting subgroup $T^\vee_\bk$ does not depend on the choice of $B$ and $T$, see~\cite[\S 1.5.5]{br} for details.)

\subsection{The affine flag variety and Gaitsgory's central functor}
\label{ss:Fl}

We now choose a Borel subgroup $B \subset G$. We have a natural group scheme morphism
\[
 \mathrm{ev} : \Loop^+ G \to G,
\]
induced by the evaluation of $z$ at $0$. The preimage of $B$ under this morphism will be denoted $\Iw$. One can then define the \emph{affine flag variety} as the fppf quotient
\[
 \Fl_G = (\Loop G/\Iw)_{\mathrm{fppf}}.
\]
As for $\Gr_G$, this functor is represented by an ind-projective ind-scheme.
 Moreover the natural projection $\pi : \Fl_G \to \Gr_G$ is a smooth projective morphism, all of whose fibers are isomorphic to the flag variety $G/B$.

As for the category $\sfD_{\Loop^+G, \Loop^+G}$, there exists a natural structure of monoidal category on the $\Iw$-equivariant derived category $\sfD_{\Iw,\Iw}=\Db_{\Iw}(\Fl_G,\bk)$, whose product will be denoted $\star_{\Iw}$, and whose unit object is the skyscraper sheaf $\delta_\Fl$ at the base point of $\Fl_G$. This category also admits a perverse t-structure, whose heart will be denoted $\sfP_{\Iw,\Iw}$; however in this setting it is \emph{not} true that a convolution of perverse sheaves is perverse.

Following an earlier construction of Gaitsgory~\cite{gaitsgory} together with a remark of Heinloth~\cite{heinloth}, Zhu constructs in~\cite{zhu-conj} an ind-scheme
\begin{equation}
\label{eqn:central-Gr}
 \Gr_G^{\mathrm{Cen}} \to \mathbb{A}^1_\F
\end{equation}
whose fiber over $0$ identifies canonically with $\Fl_G$, and whose restriction to $\mathbb{A}^1 \smallsetminus \{0\}$ identifies with $\Gr_G \times (\mathbb{A}^1 \smallsetminus \{0\})$. Using this ind-scheme one can then define the ``central functor''
\[
 \sfZ : \sfP_{\Loop^+G, \Loop^+G} \to \Perv(\Fl_G,\bk)
\]
by setting
\[
 \sfZ(\scA) := \Psi_{\Gr_G^{\mathrm{Cen}}} \bigl( \scA \lboxtimes_\bk \underline{\bk}_{\mathbb{A}^1 \smallsetminus \{0\}}[1] \bigr),
\]
where $\Psi_{\Gr_G^{\mathrm{Cen}}}$ is the nearby cycles functor associated with the morphism~\eqref{eqn:central-Gr}. (Here $\Perv(\Fl_G,\bk)$ denotes the category of all $\bk$-perverse sheaves on $\Fl_G$.)

This functor will be fundamental for our constructions. Its main properties are summarized in the following statement. (In~\eqref{it:gaitsgory-exact} we consider the natural ``extension'' of $\star_\Iw$ to a bifunctor from $\Db_{\mathrm{c}}(\Fl_G,\bk) \times \sfD_{\Iw,\Iw}$ to $\Db_{\mathrm{c}}(\Fl_G,\bk)$.)

\begin{thm}[Gaitsgory]
\phantomsection
\label{thm:gaitsgory}
\begin{enumerate} 
\item
For any $\scA\in  \sfP_{\Loop^+G, \Loop^+G}$ the perverse sheaf $\sfZ(\scA)$ is $\Iw$-equivariant; in other words, the functor $\sfZ$ factors through a functor
\[
 \sfP_{\Loop^+G, \Loop^+G} \to \sfP_{\Iw,\Iw},
\]
which will be denoted similarly.
\item 
\label{it:gaitsgory-exact}
For any $\scA\in \sfP_{\Loop^+G, \Loop^+G}$ the perverse sheaf $\sfZ(\scA)$ is convolution-exact, in the sense that for any $\scF\in \Perv(\Fl_G,\bk)$ the convolution product $\scF \star_{\Iw} \sfZ(\scA)$ is a perverse sheaf.
\item 
\label{it:gaitsgory-central}
For any $\scA \in \sfP_{\Loop^+G, \Loop^+G}$ and $\scF\in \sfD_{\Iw,\Iw}$, there exists a canonical (in particular, bifunctorial) isomorphism
\[
\sigma_{\scA,\scF} : \sfZ(\scA)\star_{\Iw} \scF \simto \scF \star_{\Iw} \sfZ(\scA).
\]
\item \label{it:gaitsgory-monoidal}
The composition
\[
\sfP_{\Loop^+G, \Loop^+G}\xrightarrow{\sfZ} \sfP_{\Iw,\Iw} \to \sfD_{\Iw,\Iw}
\]
has a canonical monoidal structure; in other words
we have an identification $\sfZ(\delta_{{\Gr}}) = \delta_{{\Fl}}$, for $\scA, \scB\in  \sfP_{\Loop^+G, \Loop^+G}$ there is a canonical isomorphism
\[
\phi_{\scA,\scB} : \sfZ(\scA \star_{\Loop^+G} \scB) \simto \sfZ(\scA)\star_{\Iw}\sfZ(\scB),
\]
and these isomorphisms intertwine the associativity constraints of the mono\-idal categories $\sfP_{\Loop^+G, \Loop^+G}$ and $\sfD_{\Iw,\Iw}$ in the obvious way.
\item 
\label{it:Z-pi!}
For $\scA\in \sfP_{\Loop^+G, \Loop^+G}$, we have a canonical isomorphism
\[
\pi_*(\sfZ(\scA))\simto\scA.
\]
\end{enumerate}
\end{thm}

For the original proof of these properties (in the setting where $\bk$ is replaced by $\overline{\mathbb{Q}}_\ell$, and using a different version of the ind-scheme $\Gr_G^{\mathrm{Cen}}$), see~\cite{gaitsgory}.\footnote{In~\cite{gaitsgory} the isomorphism $\sigma_{\scA,\scF}$ is constructed only in the special case when $\scF$ is perverse, which causes complications for our later considerations.} For a modification of these proofs using the present definition of $\Gr_G^{\mathrm{Cen}}$ (but still in the setting where $\bk$ is replaced by $\overline{\mathbb{Q}}_\ell$), see~\cite[\S 7]{zhu-conj}. For a detailed review of this proof, which allows more general coefficient rings, see~\cite{ar-book}.

%
%

Recall the notion of a central functor (see~\S\ref{ss:statement-reconstruction}). By Theorem~\ref{thm:gaitsgory} the functor $\sfZ$ admits a natural monoidal structure, and we have the natural isomorphisms $\sigma_{\scA,\scF}$ for $\scA$ in $\sfP_{\Loop^+G, \Loop^+G}$ and $\scF$ in $\sfD_{\Iw,\Iw}$. These data are exactly the ingredients that enter the definition of central functors. The following statement says that this triple indeed is a central functor.

\begin{thm}[Gaitsgory]
\label{thm:Z-central}
The triple $(\sfZ,\phi_{-,-},\sigma_{-,-})$ is a central functor from $\sfP_{\Loop^+G, \Loop^+G}$ to $\sfD_{\Iw,\Iw}$.
\end{thm}

See~\cite{gaitsgory-app} for the original proof of Theorem~\ref{thm:Z-central}, written in the setting where $\bk$ is replaced by $\overline{\mathbb{Q}}_\ell$. See~\cite{ar-book} for a detailed review of this proof, adapted to the present setting.

\subsection{Monodromy}

For any $\scA$ in $\sfP_{\Loop^+G, \Loop^+G}$, since the perverse sheaf $\sfZ(\scA)$ is defined in terms of nearby cycles, it comes equipped with a canonical \emph{monodromy automorphism}\footnote{More precisely, $\sfZ(\scA)$ comes equipped with a canonical action of $\Z_\ell(1)$, the inverse limit of the groups of $\ell^n$-th roots of unity in $\F$ for all $n$. Below we fix once and for all a topological generator of this group, and consider the action of this element.}
\[
 \sm_{\scA} : \sfZ(\scA) \simto \sfZ(\scA).
\]
Below we will need the following properties of these automorphisms.

\begin{prop}
\phantomsection
\label{prop:monodromy}
 \begin{enumerate}
 \item
 \label{it:monodromy-0}
 For any $\scA,\scB$ in $\sfP_{\Loop^+G, \Loop^+G}$,
 under the isomorphism $\phi_{\scA,\scB}$ we have
 \[
 \sm_{\scA \star_{\Loop^+G} \scB} = \sm_{\scA} \star_\Iw \sm_{\scB}.
 \]
  \item 
  \label{it:monodromy-1}
  For any $\scA$ in $\sfP_{\Loop^+G, \Loop^+G}$, the automorphism $\sm_{\scA}$ is unipotent.
  \item 
  \label{it:monodromy-2}
  For any $\scA,\scB \in \sfP_{\Loop^+G, \Loop^+G}$ and $f \in \Hom_{\sfP_{\Iw,\Iw}}(\sfZ(\scA), \sfZ(\scB))$ we have
  \[
  \sm_{\scB} \circ f = f \circ \sm_{\scA}.
  \]
 \end{enumerate}
\end{prop}

Here the proof of~\eqref{it:monodromy-0} is easy, see~\cite[Proposition~3.4.2]{ar-book}.
For the original proof of~\eqref{it:monodromy-1}, in the setting where $\bk$ is replaced by $\overline{\mathbb{Q}}_\ell$, see~\cite[Proposition~7]{gaitsgory}. This proof does not easily adapt to our present setting; however an alternative proof, which also implies property~\eqref{it:monodromy-2}, is provided in~\cite{bm}. This proof does work in our present setting, see~\cite[Proposition~2.4.6]{ar-book} for details.

The proof of Proposition~\ref{prop:monodromy}\eqref{it:monodromy-1}--\eqref{it:monodromy-2} requires a different (but equivalent) description of the automorphisms $\sm_{\scA}$, which we now explain for later use. 
Consider the ``loop rotation'' action of $\Gm$ on $\Fl_G$, induced by the action on $\F(\hspace{-1pt}(z)\hspace{-1pt})$ which ``rescales'' $z$.
Since the $\Iw$-orbits on $\Fl_G$ are stable under this action, every object $\scF$ of $\sfP_{\Iw,\Iw}$ is monodromic for this action in the sense of~\cite{verdier}; therefore it admits a canonical automorphism $\mathfrak{M}_\scF$\footnote{More precisely, $\scF$ admits a canonical action of $\Z_\ell(1)$; we deduce a canonical automorphism by taking the image of our fixed topological generator of this group.} (see~\cite{bm} or~\cite[\S 9.5.3]{ar-book} for details and more precise references), and $\mathfrak{M}_{(-)}$ defines an endomorphism of the identity functor $\id_{\sfP_{\Iw,\Iw}}$. As explained in~\cite[\S 5.2]{bm}, for any $\scA$ in $\sfP_{\Loop^+G, \Loop^+G}$ we have
\begin{equation}
\label{eqn:sm-monodromy}
\sm_{\scA} = (\mathfrak{M}_{\sfZ(\scA)})^{-1}.
\end{equation}

\subsection{Iwahori orbits}
\label{ss:Iw-orbits}

From now on we fix a choice of maximal torus $T \subset B$.
We then obtain a canonical identification $X_*(T) = X^*(T^\vee_\bk)$, see~\S\ref{ss:Satake}. We will denote by $\mathfrak{R}$ the root system of $(G,T)$.
The choice of $B$ determines a choice of positive roots $\mathfrak{R}_+ \subset \mathfrak{R}$ (such that the $T$-weights in the Lie algebra of $B$ are the \emph{negative} roots), hence a subset 
$X_*^+(T) \subset X_*(T)$
of dominant coweights. We will also denote by $B^\vee_\bk \subset G^\vee_\bk$ the Borel subgroup containing $T^\vee_\bk$ such that the $T^\vee_\bk$-weights on the Lie algebra of $B^\vee_\bk$ are the negative coroots of $G$.


We denote by $\Wf$ the Weyl group of $(G,T)$. Then $\Wf$ acts canonically on 
$X_*(T)$,
and the semi-direct product
\[
 W:= \Wf \ltimes X_*(T)
\]
is called the (extended) affine Weyl group. This group parametrizes the $\Iw$-orbits on $\Fl_G$ in the following way. For $\lambda \in X_*(T)$,
we will denote by $\st(\lambda)$ the image of $\lambda$ in $W$. Then if $w=\st(\lambda)x$ for some $\lambda \in X_*(T)$
and $x \in \Wf$, we denote by $\Fl_{G,w} \subset \Fl_G$ the $\Iw$-orbit of the coset $z^\lambda \dot{x} \Iw$, where:
\begin{itemize}
 \item $z^\lambda$ is the image of $z$ under the morphism $\F(\hspace{-1pt}(z)\hspace{-1pt})^\times \to \Loop G$ induced by $\lambda$;
 \item $\dot{x}$ is any lift of $x$ in the normalizer $\mathrm{N}_G(T)$.
\end{itemize}
Then it is well known that we have
\[
 (\Fl_G)_{\mathrm{red}} = \bigsqcup_{w \in W} \Fl_{G,w}.
\]

For $w \in W$ we will denote by $j_w : \Fl_{G,w} \to \Fl_G$ the embedding, and set
\[
 \ell(w):=\dim(\Fl_{G,w}).
\]
 A formula due to Iwahori--Matsumoto~\cite{im} states that if $w=x \st(\lambda)$ with $x \in \Wf$ and $\lambda \in X_*(T)$,
 we have
\begin{equation}
\label{eqn:formula-length}
 \ell(w)=\sum_{\substack{\alpha \in \mathfrak{R}_+ \\ x(\alpha) \in \mathfrak{R}_+}} |\langle \lambda, \alpha \rangle | + \sum_{\substack{\alpha \in \mathfrak{R}_+ \\ x(\alpha) \in -\mathfrak{R}_+}} |\langle \lambda, \alpha \rangle + 1 |.
\end{equation}
Set
\[
 W^{\mathrm{Cox}}:= \Wf \ltimes \Z\mathfrak{R}^\vee,
\]
where $\mathfrak{R}^\vee \subset X_*(T)$
is the coroot system of $(G,T)$ and $\Z\mathfrak{R}^\vee \subset X_*(T)$
is the coroot lattice. If we set $S=\{w \in W^{\mathrm{Cox}} \mid \ell(w)=1\}$, then it is well known that $(W^{\mathrm{Cox}},S)$ is a Coxeter system, whose length function is the restriction of $\ell$. Moreover, if we set $\Omega:=\{w \in W \mid \ell(w)=0\}$, then for any $\omega \in \Omega$ conjugation by $\omega$ stabilizes $S$, hence acts on $W^{\mathrm{Cox}}$ as a Coxeter group automorphism. Finally, the product induces a group isomorphism
\[
 \Omega \ltimes W^{\mathrm{Cox}} \simto W,
\]
such that $\ell(\omega w)=\ell(w)$ for $\omega \in \Omega$ and $w \in W^{\mathrm{Cox}}$.

For any $w \in W$,
we define the \emph{standard} and \emph{costandard} perverse sheaves associated with $w$ as
\[
 \Delta^{\Iw}_w := j_{w!} \underline{\bk}_{\Fl_{G,w}}[\ell(w)], \quad \nabla^{\Iw}_w := j_{w*} \underline{\bk}_{\Fl_{G,w}}[\ell(w)].
\]
These complexes are indeed perverse sheaves because $j_w$ is an affine embedding, and are $\Iw$-equivariant. It is well known that 
if $w,y \in W$ are such that $\ell(wy)=\ell(w)+\ell(y)$ there exist canonical isomorphisms
\begin{equation}
\label{eqn:Delta-nabla-conv}
 \Delta^{\Iw}_w \star_{\Iw} \Delta^{\Iw}_y \cong \Delta^{\Iw}_{wy}, \quad \nabla^{\Iw}_w \star_{\Iw} \nabla^{\Iw}_y \cong \nabla^{\Iw}_{wy},
\end{equation}
and that
for any $w \in W$ there exist isomorphisms
\begin{equation}
\label{eqn:Delta-nabla-inverse}
 \Delta^{\Iw}_w \star_{\Iw} \nabla^{\Iw}_{w^{-1}} \cong \delta_{\Fl} \cong \nabla^{\Iw}_{w^{-1}} \star_{\Iw} \Delta^{\Iw}_w.
\end{equation}
In particular, the objects $\Delta^{\Iw}_w$ and $\nabla^{\Iw}_w$ are invertible (in the monoidal category $\sfD_{\Iw,\Iw}$) for any $w \in W$.

For any $w \in W$, the intersection cohomology complex of $\Fl_{G,w}$ (i.e.~the image of the unique---up to scalar---nonzero morphism $\Delta^{\Iw}_w \to \nabla^{\Iw}_w$) will be denoted $\IC_w$; then the assignment $w \mapsto \IC_w$ induces a bijection between $W$ and the set of isomorphism classes of simple objects in $\sfP_{\Iw,\Iw}$.

The following lemma is an analogue of~\cite[Lemma~2.1]{bbm} for $\Fl_G$. The same arguments as in~\cite{bbm} apply.

\begin{lem}
\label{lem:socle-Deltas}
Let $w \in W$, and assume that $w=\omega x$ for some $\omega \in \Omega$ and $x \in W^{\mathrm{Cox}}$. Then the socle of $\Delta^{\Iw}_w$ is $\IC_\omega$, and the the cokernel of the embedding $\IC_\omega \hookrightarrow \Delta^{\Iw}_w$ does not contain any object of the form $\IC_{\omega'}$ with $\omega' \in \Omega$ as a composition factor. Dually, the top of $\nabla_w^{\Iw}$ is $\IC_\omega$, and the kernel of the surjection $\nabla^{\Iw}_w \twoheadrightarrow \IC_\omega$ does not contain any object of the form $\IC_{\omega'}$ with $\omega' \in \Omega$ as a composition factor.
\end{lem}

\subsection{Wakimoto sheaves}
\label{ss:Wakimoto-sheaves}

Fix $\lambda \in X_*(T)$.
Given $\mu,\nu \in X_*^+(T)$ such that $\lambda=\mu-\nu$ and $\scF$ in $\sfD_{\Iw,\Iw}$, we can consider the $\bk$-vector space 
\[
\Hom_{\sfD_{\Iw,\Iw}}(\nabla^{\Iw}_{\st(\mu)}, \nabla^{\Iw}_{\st(\nu)} \star_{\Iw} \scF).
\]
If $\mu',\nu' \in X_*^+(T)$ are such that $\lambda=\mu'-\nu'$ and $\mu'-\mu \in X_*^+(T)$, then convolution with $\nabla^{\Iw}_{\st(\mu'-\mu)}$ induces a canonical isomorphism
\[
 \Hom_{\sfD_{\Iw,\Iw}}(\nabla^{\Iw}_{\st(\mu)}, \nabla^{\Iw}_{\st(\nu)} \star_{\Iw} \scF) \simto \Hom_{\sfD_{\Iw,\Iw}}(\nabla^{\Iw}_{\st(\mu')}, \nabla^{\Iw}_{\st(\nu')} \star_{\Iw} \scF).
\]
One can check that if $(\mu'',\nu'')$ is another pair such that $\lambda=\mu''-\nu''$ and $\mu''-\mu' \in X_*^+(T)$, the two isomorphisms
\[
 \Hom_{\sfD_{\Iw,\Iw}}(\nabla^{\Iw}_{\st(\mu)}, \nabla^{\Iw}_{\st(\nu)} \star_{\Iw} \scF) \simto \Hom_{\sfD_{\Iw,\Iw}}(\nabla^{\Iw}_{\st(\mu'')}, \nabla^{\Iw}_{\st(\nu'')} \star_{\Iw} \scF)
\]
obtained by first convolving with $\nabla^{\Iw}_{\st(\mu'-\mu)}$ and then with $\nabla^{\Iw}_{\st(\mu''-\mu')}$, or directly with $\nabla^{\Iw}_{\st(\mu''-\mu)}$, coincide; we can therefore consider the (filtrant) direct limit
\[
 \varinjlim_{\substack{(\mu,\nu) \in (X_*^+(T))^2 \\ \lambda=\mu-\nu}} \Hom_{\sfD_{\Iw,\Iw}}(\nabla^{\Iw}_{\st(\mu)}, \nabla^{\Iw}_{\st(\nu)} \star_{\Iw} \scF)
\]
with respect to the order such that $(\mu,\nu) \unlhd (\mu',\nu')$ if $\mu'-\mu \in X_*^+(T)$. We denote by $\Wak_\lambda \in \sfD_{\Iw,\Iw}$ the object which represents the functor
\[
 \scF \mapsto \varinjlim_{\substack{(\mu,\nu) \in (X_*^+(T))^2 \\ \lambda=\mu-\nu}} \Hom_{\sfD_{\Iw,\Iw}}(\nabla^{\Iw}_{\st(\mu)}, \nabla^{\Iw}_{\st(\nu)} \star_{\Iw} \scF).
\]
(This functor is indeed representable since $\nabla^{\Iw}_{\st(\nu)}$ is invertible, see~\eqref{eqn:Delta-nabla-inverse}.)


Following Mirkovi\'c (to whom this definition is due),
the objects $(\Wak_\lambda : \lambda \in X_*(T))$ are called \emph{Wakimoto sheaves}. Their properties are investigated in~\cite{bm} (in the setting where $\bk$ is replaced by $\overline{\mathbb{Q}}_\ell$) and in~\cite[Chap.~4]{ar-book} (in the general setting). Those which will be useful for us are the following. (Here, we denote by $\preceq$ the order on $X_*(T)$ such that $\lambda \preceq \mu$ iff $\mu-\lambda$ is a sum of positive roots.)
\begin{enumerate}
\item 
\label{it:Wakimoto-Delta-nabla}
If $\lambda=\mu-\nu$ with $\mu,\nu \in X_*^+(T)$ we have \emph{noncanonical} isomorphisms
\[
\Wak_\lambda \cong \Delta_{\st(-\nu)}^{\Iw} \star_{\Iw} \nabla^{\Iw}_{\st(\mu)} \cong \nabla^{\Iw}_{\st(\mu)} \star_{\Iw} \Delta_{\st(-\nu)}^{\Iw}.
\]
\item
For any $\lambda,\lambda' \in X_*(T)$ there exists a canonical isomorphism
\[
\Wak_\lambda \star_{\Iw} \Wak_{\lambda'} \simto \Wak_{\lambda + \lambda'}.
\]
\item
\label{it:Wak-support}
For any $\lambda \in X_*(T)$, the object $\Wak_\lambda$ 
is a perverse sheaf supported on $\overline{\Fl_{G,\st(\lambda)}}$, and its restriction to $\Fl_{G,\st(\lambda)}$ is isomorphic to $\underline{\bk}_{\Fl_{G,\st(\lambda)}}[\ell(\st(\lambda))]$.
 \item
 \label{it:Hom-Wak}
 For $\lambda,\lambda' \in X_*(T)$ and $n \in \Z$ we have
 \begin{equation*}
 \Hom_{\sfD_{\Iw,\Iw}}(\Wak_\lambda, \Wak_{\lambda'}[n])=0 \quad \text{unless $\lambda' \preceq \lambda$.}
 \end{equation*}
 \item
 \label{it:Wak-End}
 For $\lambda \in X_*(T)$ we have
 \begin{equation*}
 \Hom_{\sfP_{\Iw,\Iw}}(\Wak_\lambda,\Wak_\lambda)=\bk \quad \text{and} \quad \Ext^1_{\sfP_{\Iw,\Iw}}(\Wak_\lambda,\Wak_\lambda)=0.
 \end{equation*}
\end{enumerate}

\begin{rmk}
\label{rmk:Wakop}
Following~\cite{zhu-conj}, in~\cite{ar-book} Wakimoto sheaves are defined in a more general setting, where $X_*^+(T)$ is replaced by the intersection of $X_*(T)$ with any fixed choice of Weyl chamber in $\mathbb{R} \otimes_{\Z} X_*(T)$. 
All the properties listed above remain true in this generality, with the appropriate replacement for~\eqref{it:Wakimoto-Delta-nabla} and for the order in~\eqref{it:Hom-Wak}.
Below we will also use the Wakimoto sheaves associated with the \emph{anti-dominant} chamber, which will be denoted $\Wakop_\lambda$ ($\lambda \in X_*(T)$). So, if $\lambda = \mu-\nu$ with $\mu,\nu \in X_*^+(T)$ then we have
\[
\Wakop_\lambda \cong \nabla^{\Iw}_{\st(-\nu)} \star_{\Iw} \Delta^{\Iw}_{\st(\mu)} \cong \Delta^{\Iw}_{\st(\mu)} \star_{\Iw} \nabla^{\Iw}_{\st(-\nu)}.
\]
We also have $\Hom_{\sfD_{\Iw,\Iw}}(\Wakop_\lambda, \Wakop_{\lambda'}[n])=0$ unless $\lambda' \succeq \lambda$. 
\end{rmk}

\subsection{Wakimoto filtrations}
\label{ss:filtration-central-sheaves}

We will denote by
$\sfP_{\Iw,\Iw}^{\Wak}$ the full subcategory of $\sfP_{\Iw,\Iw}$ whose objects are the perverse sheaves which admit a filtration by Wakimoto sheaves, i.e.~a filtration whose subquotients are of the form $\Wak_\lambda$ with $\lambda \in X_*(T)$. This category admits a natural structure of exact category, inherited from the ambiant abelian category $\sfP_{\Iw,\Iw}$. Since a convolution of Wakimoto sheaves is a Wakimoto sheaf (in particular, is perverse), $\sfP_{\Iw,\Iw}^{\Wak}$ is a monoidal subcategory of $(\sfD_{\Iw,\Iw},\star_{\Iw})$.

It follows from properties~\eqref{it:Hom-Wak}--\eqref{it:Wak-End} in~\S\ref{ss:Wakimoto-sheaves} that for any $\scF$ in $\sfP_{\Iw,\Iw}^{\Wak}$ and any ideal $\mathbf{Y} \subset X_*(T)$ (i.e.~any subset such that if $\lambda \in \mathbf{Y}$ and $\mu \preceq \lambda$ then $\mu \in \mathbf{Y}$) there exists a unique subobject $\scF_{\bY} \subset \scF$ which is an extension of objects $\Wak_\lambda$ with $\lambda \in \mathbf{Y}$ and such that the quotient $\scF/\scF_{\bY}$ is an extension of objects $\Wak_\mu$ with $\mu \in X_*(T) \smallsetminus \mathbf{Y}$; moreover 
the assignment $\scF \mapsto \scF_{\mathbf{Y}}$ is functorial. 

In particular, given $\lambda \in X_*(T)$ one can choose an ideal $\mathbf{Y} \subset X_*(T)$ such that $\lambda \in \mathbf{Y}$ is maximal, and consider the quotient
\[
 \mathrm{gr}_\lambda(\scF) := \scF_\mathbf{Y}/\scF_{\mathbf{Y} \smallsetminus \{\lambda\}}.
\]
By construction $\mathrm{gr}_\lambda(\scF)$ is a direct sum of copies of $\Wak_\lambda$, and we set
\[
 \mathrm{Grad}_\lambda(\scF) := \Hom_{\sfP_{\Iw,\Iw}}(\Wak_\lambda, \mathrm{gr}_\lambda(\scF)).
\]

\begin{lem}
\label{lem:gr-exact}
 The object $\mathrm{gr}_\lambda(\scF)$ does not depend on the choice of $\mathbf{Y}$. Moreover, the functor $\scF \mapsto \mathrm{Grad}_\lambda(\scF)$ is exact (for the exact structure on $\sfP_{\Iw,\Iw}^{\Wak}$ considered above).
\end{lem}

\begin{proof}
 For the first claim, see~\cite[Lemma~4.3.4]{ar-book}. For the second claim, see~\cite[Proposition~4.6.1]{ar-book}.
\end{proof}

The following result is due to Arkhipov and the first author~\cite[Theorem~4(a)]{ab} in the setting where $\bk$ is replaced by $\overline{\mathbb{Q}}_\ell$. The proof is also reproduced in~\cite[\S 7.3]{zhu-conj} and (in the present setting) in~\cite[\S 4.4]{ar-book}.

\begin{thm}[Arkhipov--Bezrukavnikov]
\label{thm:Wak-filtration}
For any $\scF$ in $\sfP_{\Loop^+G, \Loop^+G}$, the perverse sheaf $\sfZ(\scF)$ belongs to $\sfP_{\Iw,\Iw}^{\Wak}$.
\end{thm}

\begin{rmk}
\label{rmk:Wakop-2}
In~\cite{zhu-conj} and~\cite{ar-book}, Theorem~\ref{thm:Wak-filtration} is proved with the more general definition of Wakimoto sheaves alluded to in Remark~\ref{rmk:Wakop}.
In particular, for any $\scF$ in $\sfP_{\Loop^+G, \Loop^+G}$ the perverse sheaf $\sfZ(\scF)$ also admits a filtration with subquotients of the form $\Wakop_\lambda$ with $\lambda \in X_*(T)$.
\end{rmk}

To simplify notation, from now on we will use the notation
\[
\sZ := \sfZ \circ \Sat^{-1} : \Rep(G^\vee_\bk) \to \sfP_{\Iw,\Iw},
\]
where $\Sat$ is as in~\eqref{eqn:Satake-equiv},
and write $\sm_V$ for $\sm_{\Sat^{-1}(V)}$.

Thanks to Theorem~\ref{thm:Wak-filtration}, we can consider for any $V$ in $\Rep(G^\vee_\bk)$ and any $\lambda \in X_*(T)$ the vector space
$\mathrm{Grad}_\lambda(\sZ(V))$.
By~\cite[Lemma~4.8.1]{ar-book}, we have a canonical identification
\begin{equation}
\label{eqn:Grad-wt-space}
\mathrm{Grad}_\lambda(\sZ(V)) \cong V_{w_0(\lambda)},
\end{equation}
where $w_0$ is the longest element in $\Wf$. (Here, in the right-hand side we consider the weight space for the action of the canonical maximal torus $T^\vee_\bk \subset G^\vee_\bk$.)

\section{The regular quotient}
\label{sec:regular-quotient}

In this section we introduce the main player of the paper, namely the ``regular quotient'' of the category $\sfP_{\Iw,\Iw}$.

\subsection{The regular quotient}
\label{ss:regular-quotient}


We consider the Serre subcategory $\langle \IC_w : \ell(w) >0 \rangle_{\mathrm{Serre}} \subset \sfP_{\Iw,\Iw}$ generated by the objects $\IC_w$ with $w \in W$ satisfying $\ell(w)>0$, and the Serre quotient
\[
\sfP_{\Iw,\Iw}^0 := \sfP_{\Iw,\Iw} / \langle \IC_w : \ell(w) >0 \rangle_{\mathrm{Serre}}.
\]
The quotient functor $\sfP_{\Iw,\Iw} \to \sfP_{\Iw,\Iw}^0$ will be denoted $\Pi_{\Iw,\Iw}^0$. Every object in the category $\sfP_{\Iw,\Iw}^0$ has finite length, and the assignment $w \mapsto \Pi_{\Iw,\Iw}^0(\IC_w)$ induces a bijection between $\Omega$ and the set of isomorphism classes of simple objects in $\sfP_{\Iw,\Iw}^0$.

The proof of the following lemma can be copied verbatim from~\cite[\S 6.5.6]{ar-book}.

\begin{lem}\phantomsection
\label{lem:convolution-0}
\begin{enumerate}
\item
\label{it:convolution-0-1}
If $\scF$ belongs to $\langle \IC_w : \ell(w) >0 \rangle_{\mathrm{Serre}}$ and $\scG$ is any object of $\sfP_{\Iw,\Iw}$, then for any $n \in \Z$ the perverse sheaves $\pH^n(\scF \star_{\Iw} \scG)$ and $\pH^n(\scG \star_{\Iw} \scF)$ belong to $\langle \IC_w : \ell(w) >0 \rangle_{\mathrm{Serre}}$.
\item
\label{it:convolution-0-2}
For $\scF,\scG$ in $\sfP_{\Iw,\Iw}$ and any $n \in \Z \smallsetminus \{0\}$ we have $\Pi_{\Iw,\Iw}^0(\pH^n(\scF \star_{\Iw} \scG))=0$.
\end{enumerate}
\end{lem}

Lemma~\ref{lem:convolution-0} has the following consequence.

\begin{prop}\phantomsection
\label{prop:PI0-monoidal}
\begin{enumerate}
\item
The bifunctor
\[
\sfP_{\Iw,\Iw} \times \sfP_{\Iw,\Iw} \to \sfP_{\Iw,\Iw}^0
\]
sending a pair $(\scF,\scG)$ to $\Pi^0(\pH^0(\scF \star_{\Iw} \scG))$ factors through a bifunctor
\[
(-) \star_{\Iw}^0 (-) : \sfP_{\Iw,\Iw}^0 \times \sfP_{\Iw,\Iw}^0 \to \sfP_{\Iw,\Iw}^0.
\]
\item
The bifunctor $\star_{\Iw}^0$ is exact on both sides, and admits natural associativity and unit constraints.
\end{enumerate}
\end{prop}

\begin{proof}
Lemma~\ref{lem:convolution-0} implies that the bifunctor $(\scF,\scG) \mapsto \Pi_{\Iw,\Iw}^0(\pH^0(\scF \star_{\Iw} \scG))$ is exact on both sides, and vanishes on $\langle \IC_w : \ell(w) >0 \rangle_{\mathrm{Serre}} \times \sfP_{\Iw,\Iw}$ and on $\sfP_{\Iw,\Iw} \times \langle \IC_w : \ell(w) >0 \rangle_{\mathrm{Serre}}$. Therefore it factors through a bifunctor $\star^0_{\Iw}$, which is exact on both sides (see~\cite[Corollaires~2--3, p.~368--369]{gabriel}). This observation also shows that for $\scF_1, \scF_2, \scF_3$ in $\sfP_{\Iw,\Iw}^0$ we have canonical isomorphisms
\[
(\scF_1 \star^0_{\Iw} \scF_2) \star^0_{\Iw} \scF_3 \cong \Pi_{\Iw,\Iw}^0(\pH^0(\scF_1 \star_{\Iw} \scF_2 \star_{\Iw} \scF_3)) \cong \scF_1 \star^0_{\Iw} (\scF_2 \star^0_{\Iw} \scF_3),
\]
which provides the desired associativity constraint. The unit object for this product is $\delta^0:=\Pi_{\Iw,\Iw}^0(\delta_{\Fl})$.
\end{proof}

Thanks to Proposition~\ref{prop:PI0-monoidal}, we can now consider the abelian monoidal category $(\sfP_{\Iw,\Iw}^0, \star_{\Iw}^0)$. We set
\[
\sZ^0 :=\Pi_{\Iw,\Iw}^0 \circ \sZ : \Rep(G_\bk^\vee) \to \sfP_{\Iw,\Iw}^0,
\]
and for $V$ in $\Rep(G^\vee_\bk)$ we set
$\sm_V^0:=\Pi_{\Iw,\Iw}^0(\sm_V) \in \End_{\sfP_{\Iw,\Iw}^0}(\sZ^0(V))$.

The following claim is an immediate consequence of Proposition~\ref{prop:monodromy}\eqref{it:monodromy-2}.

\begin{lem}
\label{lem:monodromy-commutes}
 For any $V,V'$ in $\Rep(G_\bk^\vee)$ and any $f \in \Hom_{\sfP_{\Iw,\Iw}^0}(\sZ^0(V), \sZ^0(V'))$ we have $f \circ \sm_V^0 = \sm_{V'}^0 \circ f$.
\end{lem}

\subsection{Main result: statement and strategy of proof}
\label{ss:main}

The main result of the present paper can be stated as follows.

\begin{thm}
\label{thm:main}
Assume that the following conditions hold:
\begin{enumerate}
 \item 
 \label{it:assumption-main-1}
 either
 $X^*(T) / \Z\mathfrak{R}$ is free and
 $X_*(T)/\Z\mathfrak{R}^\vee$ has no $\ell$-torsion, or $\ell$ is very good;
 \item 
 \label{it:assumption-main-3}
 for any indecomposable factor of the root system $\mathfrak{R}$, the prime $\ell$ is strictly bigger than the corresponding value in the table of Figure~\ref{fig:bounds}.
\end{enumerate}
Then there exists a regular unipotent element $\su \in G^\vee_\bk$ and an equivalence of monoi\-dal categories
\[
\Phi^0 : (\sfP_{\Iw,\Iw}^0, \star_{\Iw}^0) \simto (\Rep(\rmZ_{G^\vee_\bk}(\su)),\otimes)
\]
such that $\Phi^0 \circ \sZ^0 \cong \For^{G^\vee_\bk}_{\rmZ_{G^\vee_\bk}(\su)}$.
\end{thm}

\begin{figure}
 \begin{tabular}{|c|c|c|c|c|c|c|c|c|}
  \hline
  $\mathbf{A}_n$ ($n \geq 1$) & $\mathbf{B}_n$ ($n \geq 2$) & $\mathbf{C}_n$ ($n \geq 3$) & $\mathbf{D}_n$ ($n \geq 4$) & $\mathbf{E}_6$ & $\mathbf{E}_7$ & $\mathbf{E}_8$ & $\mathbf{F}_4$ & $\mathbf{G}_2$ \\
  \hline
  $1$ & $n$ & $2$ & $2$ & $3$ & $19$ & $31$ & $3$ & $3$\\
  \hline
 \end{tabular}
\caption{Bounds on $\ell$}
\label{fig:bounds}
\end{figure}

\begin{rmk}
Let us briefly discuss the assumptions in Theorem~\ref{thm:main}. 
The condition that $X^*(T) / \Z\mathfrak{R}$ is free
is equivalent to requiring that the derived subgroup of $G^\vee_\bk$ is simply connected. (Alternatively, since the quotient $X^*(T)/\Z\mathfrak{R}$ identifies with the group of characters of the center $\rmZ(G)$ of $G$, this assumption is also equivalent to the condition that $\rmZ(G)$ is a torus.) 
The assumption that $X_*(T)/\Z\mathfrak{R}^\vee$ has no $\ell$-torsion guarantees that $\rmZ(G^\vee_\bk)$ and $\rmZ_{G^\vee_\bk}(\su)$ are smooth group schemes, see~\S\ref{ss:centralizer-ureg}.
(This assumption in fact follows from~\eqref{it:assumption-main-3} if $\mathfrak{R}$ has no component of type $\mathbf{A}$.)
Finally, assumption~\eqref{it:assumption-main-3} is an artefact of our proof of the statement considered in~\ref{it:central-tilting} below; it is a likely that it can be refined (at least) to the condition that $\ell$ is good for $G$. 
\end{rmk}

Our proof of Theorem~\ref{thm:main} will proceed as follows. First, using standard considerations one can reduce the proof to the case of the first possible assumption in~\eqref{it:assumption-main-1}; see~\S\ref{ss:end-proof} for details. In this case we will
use Proposition~\ref{prop:subgroup}, applied to the abelian monoi\-dal category $(\sfP_{\Iw,\Iw}^0, \star_{\Iw}^0)$ and the functor $\sZ^0 : \Rep(G^\vee_\bk) \to \sfP_{\Iw,\Iw}^0$. For this we will need to construct a central structure on this functor; this will be done in~\S\ref{ss:central-structure-Z0}, and does not require any assumption. The first real difficulty appears when trying to show that any object in $\sfP_{\Iw,\Iw}^0$ is a subquotient of an object $\sZ^0(V)$. In fact we will not show this directly; instead we will provisionally replace the category $\sfP_{\Iw,\Iw}^0$ by the full subcategory $\widetilde{\sfP}_{\Iw,\Iw}^0$ whose objects are precisely the subquotients of objects $\sZ^0(V)$. With this replacement Proposition~\ref{prop:subgroup} applies, and provides a closed subgroup $H \subset G^\vee_\bk$, see~\S\ref{ss:identification-regular-quotient}. The rest of the paper will then be devoted to the proof of the following two claims (under the assumptions of Theorem~\ref{thm:main}), which will imply the theorem:
\begin{enumerate}[label={(\alph*)}]
\item
\label{it:tP-P}
the embedding $\widetilde{\sfP}_{\Iw,\Iw}^0 \subset \sfP_{\Iw,\Iw}^0$ is an equality;
\item
\label{it:H-centralizer}
the subgroup $H$ is the centralizer of a regular unipotent element.
\end{enumerate}
These properties will be proved in Section~\ref{sec:proof-main}. The main two ingredients of these proofs are the following:
\begin{enumerate}[label={(\roman*)}]
\item
\label{it:construction-u}
the construction of the unipotent element appearing in~\ref{it:H-centralizer} (based on the use of the monodromy automorphisms $\sm^0_V$) and the proof that this element is regular, which will be explained in Section~\ref{sec:regularity};
\item
\label{it:central-tilting}
the proof that the images of the central sheaves in the Iwahori--Whittaker category $\sfP_{\IW,\Iw}$ are tilting perverse sheaves, which will be obtained in Section~\ref{sec:central-tilting} after some preliminaries (including the definition of the category $\sfP_{\IW,\Iw}$) in Section~\ref{sec:asp-IW}.
\end{enumerate}

\subsection{Central structure on \texorpdfstring{$\sZ^0$}{Z0}}
\label{ss:central-structure-Z0}

The following claim is an easy consequence of Theorem~\ref{thm:Z-central}.

\begin{lem}
\label{lem:Z0-central}
 The functor $\sZ^0$ admits a canonical structure of central functor from $\Rep(G_\bk^\vee)$ to $\sfP_{\Iw,\Iw}^0$.
\end{lem}

\begin{proof}
Let us first explain how to define the ``monoidality'' and ``centrality'' isomorphisms for the functor $\sZ^0$. Given $V,V'$ in $\Rep(G^\vee_\bk)$, we have the canonical isomorphism
\[
\phi_{\Sat^{-1}(V),\Sat^{-1}(V')} : \sfZ(\Sat^{-1}(V) \star_{\Loop^+G} \Sat^{-1}(V')) \simto \sfZ(\Sat^{-1}(V)) \star_{\Iw} \sfZ(\Sat^{-1}(V')).
\]
see Theorem~\ref{thm:gaitsgory}\eqref{it:gaitsgory-monoidal}. Using the monoidal structure on the functor $\Sat^{-1}$, one can interpret that isomorphism as an identification
\[
\sZ(V \otimes V') \simto \sZ(V) \star_{\Iw} \sZ(V').
\]
These complexes are perverse sheaves; applying the functor $\Pi_{\Iw,\Iw}^0$ we deduce the wished-for isomorphism
\[
\phi^0_{V,V'} : \sZ^0(V \otimes V') \simto \sZ^0(V) \star^0_{\Iw} \sZ^0(V').
\]
It is easily seen that these isomorphisms define a monoidal structure on $\sZ^0$.

Next, for $V$ in $\Rep(G^\vee_\bk)$ and $\scF$ in $\sfP_{\Iw,\Iw}^0$, we must define a canonical isomorphism
\[
\sigma^0_{V, \scF}: \sZ^0(V)\star^0_\Iw \scF \simto \scF\star_\Iw^0\sZ^0(V).
\]
However, the objects of $\sfP_{\Iw,\Iw}^0$ and $\sfP_{\Iw,\Iw}$ are the same; $\scF$ can therefore also be considered as an object in $\sfP_{\Iw,\Iw}$, for which we have the canonical isomorphism
\[
\sigma_{\Sat^{-1}(V),\scF} : \sfZ(\Sat^{-1}(V)) \star_\Iw \scF \simto \scF \star_\Iw \sfZ(\Sat^{-1}(V)),
\]
see Theorem~\ref{thm:gaitsgory}\eqref{it:gaitsgory-central}. Here again both sides are perverse sheaves, see Theorem~\ref{thm:gaitsgory}\eqref{it:gaitsgory-exact}, and applying the functor $\Pi_{\Iw,\Iw}^0$ we obtain the wished-for isomorphism $\sigma^0_{V, \scF}$.

We must now check that these data satisfy properties~\eqref{it:central-functor-1}--\eqref{it:central-functor-3} from \S\ref{ss:statement-reconstruction}. Properties~\eqref{it:central-functor-1} and~\eqref{it:central-functor-3} immediately follow from the corresponding properties for the functor $\sfZ$ by applying the functor $\Pi_{\Iw,\Iw}^0$, since all the involved objects are perverse. In order to check property~\eqref{it:central-functor-2}, we consider some $V$ in $\Rep(G^\vee_\bk)$ and $\scF,\scG$ in $\sfP_{\Iw,\Iw}^0$. As above $\scF$ and $\scG$ can also be considered as objects in $\sfP_{\Iw,\Iw}$, and we know that the composition
\[
\sZ(V) \star_\Iw \scF \star_\Iw \scG \xrightarrow{\sigma_{\Sat^{-1}(V),\scF} \star_\Iw \id} \scF \star_\Iw \sZ(V) \star_\Iw \scG \\
\xrightarrow{\id \star_\Iw \sigma_{\Sat^{-1}(V),\scG}} \scF \star_\Iw \scG \star_\Iw \sZ(V)
\]
coincides with $\sigma_{\Sat^{-1}(V),\scF \star_\Iw \scG}$. Let us apply the functor $\Pi_{\Iw,\Iw}^0 \circ \pH^0$. By exactness of the functors $\sZ(V) \star_\Iw (-)$ and $(-) \star_\Iw \sZ(V)$ (see Theorem~\ref{thm:gaitsgory}\eqref{it:gaitsgory-exact}), the domain, resp.~codomain of the morphism we obtain identify with
\[
\Pi_{\Iw,\Iw}^0(\sZ(V) \star_\Iw \pH^0(\scF \star_\Iw \scG)) = \sZ^0(V) \star_\Iw^0 \scF \star_\Iw^0 \scG,
\]
resp.~with
\[
\Pi_{\Iw,\Iw}^0(\pH^0(\scF \star_\Iw \scG) \star_\Iw \sZ(V)) = \scF \star_\Iw^0 \scG \star_\Iw^0 \sZ^0(V).
\]
It is clear also that under these identifications we have $\Pi_{\Iw,\Iw}^0 \circ \pH^0(\sigma_{\Sat^{-1}(V),\scF \star_\Iw \scG})=\sigma^0_{V,\scF \star_\Iw^0 \scG}$. On the other hand, since $\sZ(V) \star_\Iw \scF$ and $\scF \star_\Iw \sZ(V)$ are perverse we have
\[
\Pi_{\Iw,\Iw}^0 \circ \pH^0(\sZ(V) \star_\Iw \scF \star_\Iw \scG)=\Pi_{\Iw,\Iw}^0(\sZ(V) \star_\Iw \scF) \star^0_\Iw \scG
\]
and
\[
\Pi_{\Iw,\Iw}^0 \circ \pH^0(\scF \star_\Iw \sZ(V) \star_\Iw \scG)=\Pi_{\Iw,\Iw}^0(\scF \star_\Iw \sZ(V)) \star_\Iw^0 \scG,
\]
and under these identifications we have $\Pi_{\Iw,\Iw}^0 \circ \pH^0(\sigma_{\Sat^{-1}(V),\scF} \star_\Iw \id_\scG)=\sigma^0_{V,\scF} \star_\Iw^0 \id_{\scG}$. Similarly we have
\[
\Pi_{\Iw,\Iw}^0 \circ \pH^0(\scF \star_\Iw \sZ(V) \star_\Iw \scG)=\scF \star^0_\Iw \Pi_{\Iw,\Iw}^0(\sZ(V) \star_\Iw \scG)
\]
and
\[
\Pi_{\Iw,\Iw}^0 \circ \pH^0(\scF \star_\Iw \scG \star_\Iw \sZ(V))= \scF \star_\Iw^0 \Pi_{\Iw,\Iw}^0(\scG \star_\Iw \sZ(V)),
\]
and under these identifications we have $\Pi_{\Iw,\Iw}^0 \circ \pH^0(\id_\scF \star_\Iw \sigma_{\Sat^{-1}(V),\scG})=\id_{\scF} \star^0_\Iw \sigma^0_{V,\scG}$. Combining these descriptions we obtain that our data satisfy property~\eqref{it:central-functor-2}, which finishes the proof.
\end{proof}

\subsection{Description of a subcategory of the regular quotient as a category of representations}
\label{ss:identification-regular-quotient}


We will denote by $\widetilde{\sfP}_{\Iw,\Iw}^0$ the full abelian subcategory of $\sfP_{\Iw,\Iw}^0$ whose objects are the subquotients of objects of the form $\sZ^0(V)$ with $V$ in $\Rep(G_\bk^\vee)$. It is not difficult (using in particular Lemma~\ref{lem:socle-Deltas}) to check that this subcategory contains all the simple objects of $\sfP_{\Iw,\Iw}^0$. Obviously, the functor $\sZ^0$ factors through a functor $\Rep(G_\bk^\vee) \to \widetilde{\sfP}^0_{\Iw,\Iw}$, which we will denote similarly.

\begin{lem}
\label{lem:convolution-tilde}
If $\scF$ and $\scG$ belong to $\widetilde{\sfP}_{\Iw,\Iw}^0$, then $\scF \star^0_{\Iw,\Iw} \scG$ also belongs to $\widetilde{\sfP}_{\Iw}^0$.
\end{lem}

\begin{proof}
If $\scF$ is a subquotient of $\sZ^0(V)$ and $\scG$ a subquotient of $\sZ^0(V')$, then by exactness of the product $\star^0_{\Iw}$ (see Proposition~\ref{prop:PI0-monoidal}) the object $\scF \star^0_{\Iw} \scG$ is a subquotient of $\sZ^0(V) \star^0_{\Iw} \sZ^0(V') \cong \sZ^0(V \otimes V')$, hence belongs to $\widetilde{\sfP}_{\Iw,\Iw}^0$.
\end{proof}

As a consequence of Lemma~\ref{lem:convolution-tilde}, $\widetilde{\sfP}_{\Iw,\Iw}^0$ admits a natural structure of abelian monoi\-dal category. Moreover, Lemma~\ref{lem:Z0-central} implies that $\sZ^0 : \Rep(G_\bk^\vee) \to \widetilde{\sfP}^0_{\Iw,\Iw}$ has a canonical structure of central functor.

\begin{prop}
\label{prop:reg-quotient-RepH}
There exist
\begin{enumerate}
\item
a subgroup scheme $H \subset G_\bk^\vee$;
\item
an element $\su \in G_\bk^\vee$ such that $H \subset \rmZ_{G_\bk^\vee}(\su)$;
\item
an equivalence of monoidal categories
\[
\Phi^0 : (\widetilde{\sfP}_{\Iw,\Iw}^0, \star_{\Iw}^0) \simto (\Rep(H), \otimes);
\]
\item
and an isomorphism of functors $\eta : \Phi^0 \circ \sZ^0 \simto \For^{G_\bk^\vee}_H$
\end{enumerate}
such that
for any $V$ in $\Rep(G_\bk^\vee)$, the endomorphism $\eta(\Phi^0(\sm_V^0))$ coincides with the action of $\su$ on $V$.
\end{prop}

\begin{proof}
Using Proposition~\ref{prop:subgroup} we obtain the subgroup scheme $H \subset G^\vee_\bk$ and the equivalence of monoidal categories $\Phi^0$ such that $\Phi^0 \circ \sZ^0 \cong \For^{G^\vee_\bk}_H$. Lemma~\ref{lem:monodromy-commutes} implies in particular that the morphisms $(\sm^0_V : V \in \Rep(G_\bk^\vee))$ define an automorphism of $\sZ^0$ as a monoidal functor; composing with $\For^H \circ \Phi^0$ (where $\For^H : \Rep(H) \to \mathrm{Vect}_\bk$ is the forgetful functor) we obtain a tensor automorphism of $\For^H \circ \Phi^0 \circ \sZ^0 \cong \For^{G_\bk^\vee}$, hence an element $\su \in G_\bk^\vee$ (see~\cite[Proposition~2.8]{dm}). The fact that this automorphism is induced by an automorphism of $\For^{G_\bk^\vee}_H$ means that $\su$ commutes with $H$, i.e.~that $H \subset \rmZ_{G_\bk^\vee}(\su)$.
\end{proof}

\begin{rmk}
The monodromy automorphism with respect to loop rotation action defines, for any object $\scF$ of $\sfP_{\Iw,\Iw}^0$, an automorphism $\mathfrak{M}_\scF^0$. It follows from~\eqref{eqn:sm-monodromy} that for $V$ in $\Rep(G^\vee_\bk)$ we have $\sm^0_V = (\mathfrak{M}^0_{\sZ^0(V)})^{-1}$, see~\eqref{eqn:sm-monodromy}. This shows in particular that the automorphism of the functor $\For^{G^\vee_\bk}$ defining $\su$ is induced by an automorphism of $\For^H : \Rep(H) \to \mathrm{Vect}_\bk$. Hence $\su$ in fact belongs to the subgroup $H$. (This observation will not play any role below.)
\end{rmk}

\section{Regularity of \texorpdfstring{$\su$}{u}}
\label{sec:regularity}

In this section we assume that $X^*(T)/\Z\mathfrak{R}$ is free.

\subsection{Statement}
\label{ss:regularity-statement}

The following statement will be crucial for us.

\begin{thm}
\label{thm:u-regular}
The element $\su \in G_\bk^\vee$ is regular unipotent.
\end{thm}

In order to prove this statement we will construct a Borel subgroup $\widetilde{B}^\vee \subset G_\bk^\vee$ containing $\su$ in its unipotent radical $\widetilde{U}^\vee$, and then show that the image of $\su$ in $\widetilde{U}^\vee / (\widetilde{U}^\vee,\widetilde{U}^\vee)$ has nontrivial components in each (simple) root subgroup, see Proposition~\ref{prop:components-u} below. 

Recall (see~\S\ref{ss:filtration-central-sheaves}) that for any $V$ in $\Rep(G^\vee_\bk)$ and any ideal $\mathbf{Y} \subset \bX$ we have a canonical subobject $\sZ(V)_{\mathbf{Y}} \subset \sZ(V)$. We define an increasing $\Z$-filtration on $\sZ(V)$ by setting, for $n \in \Z$,
\[
F^n \sZ(V) := \sZ(V)_{\{\lambda \in X_*(T) \mid \langle \lambda, 2\rho \rangle \leq n\}}.
\]
(Here, $2\rho$ denotes the sum of the positive roots of $G$.)
Composing with the exact functor $\Pi_{\Iw,\Iw}^0$ we deduce a filtration on $\sZ^0(V)$, and then applying the fiber functor $\For^H \circ \Phi^0$ we obtain an exact $\otimes$-filtration on the fiber functor $\For^H \circ \Phi^0 \circ \sZ^0 \cong \For^{G_\bk^\vee}$ in the sense of~\cite[Chap. IV, \S 2.1]{sr}. (In fact the associated graded for this filtration is the functor sending $V$ to
\begin{equation}
\label{eqn:gr-fiber-functor}
\bigoplus_{\lambda \in X_*(T)} \mathrm{Grad}_\lambda(\sZ(V)) \otimes_\bk (\Phi^0 \circ \Pi_{\Iw,\Iw}^0(\Wak_\lambda)),
\end{equation}
which is exact by Lemma~\ref{lem:gr-exact}.) By~\cite[Th\'eor\`eme in Chap.~IV, \S 2.4]{sr} (a general result due to Deligne), such a filtration is automatically ``splittable.'' Hence~\cite[Chap. IV, Prop. 2.2.5]{sr} ensures that the subgroup scheme $\widetilde{B}^\vee$ of $G_\bk^\vee$ consisting of tensor automorphisms of $\For^{G_\bk^\vee}$ which are compatible with this filtration (see~\cite[Chap.~IV, \S 2.1.3]{sr}) is a parabolic subgroup, with unipotent radical the subgroup $\widetilde{U}^\vee$ of automorphisms inducing the identity on the associated graded. The Lie algebra of $\widetilde{B}^\vee$ is also described in~\cite[Chap. IV, Prop. 2.2.5]{sr}; from the dimension of this Lie algebra we see that in our case $\widetilde{B}^\vee$ is in fact a Borel subgroup in $G^\vee_\bk$. (In this computation we use~\eqref{eqn:Grad-wt-space} and the fact that in~\eqref{eqn:gr-fiber-functor} the vector space $\Phi^0 \circ \Pi_{\Iw,\Iw}^0(\Wak_\lambda)$ is $1$-dimensional, since it is invertible in the monoidal category of vector spaces.)

For $V$ in $\Rep(G^\vee_\bk)$, the automorphism $\sm^0_V$ is induced by the automorphism $\sm_V$ of $\sZ(V)$, which must respect the filtration $F^\bullet \sZ(V)$ by functoriality of this filtration. Therefore $\sm^0_V$ respects the filtration of $\sZ^0(V)$, showing that $\su$ belongs to $\widetilde{B}^\vee$. Moreover, since $\Wak_\lambda$ is equivariant under the loop-rotation action,~\eqref{eqn:sm-monodromy} implies that the action of $\sm^0_V$ on the associated graded for our filtration is trivial. Therefore, $\su$ even belongs to $\widetilde{U}^\vee$.

We now consider the quotient
\[
\overline{U}^\vee := \widetilde{U}^\vee / (\widetilde{U}^\vee,\widetilde{U}^\vee),
\]
a commutative algebraic group. This group admits a natural action (by group automorphisms) of the ``universal maximal torus'' $A^\vee$. (This group is defined as the quotient of any Borel subgroup by its unipotent radical; so in particular we have $\widetilde{B}^\vee / \widetilde{U}^\vee \cong A^\vee$ canonically.) This action ``breaks'' $\overline{U}^\vee$ as a product of subgroups $\overline{U}^\vee_{\alpha}$ isomorphic to $\Ga$ and naturally parametrized by the opposites of the simple roots of $(G_\bk^\vee,T_\bk^\vee)$, seen as characters of $T_\bk^\vee \cong B_\bk^\vee / (B_\bk^\vee,B_\bk^\vee) \cong A^\vee$.

By a well-known characterization of regular unipotent elements (see e.g.~\cite[\S 4.1]{humphreys}), Theorem~\ref{thm:u-regular} will follow from the following claim.

\begin{prop}
\label{prop:components-u}
For any simple root $\alpha$, the component $\su_{-\alpha}$ of $\su$ in $\overline{U}^\vee_{-\alpha}$ is nontrivial.
\end{prop}

The proof of this proposition will require some preparation; it will be given in~\S\ref{ss:proof-u} below.

\subsection{Computation in rank \texorpdfstring{$1$}{1}}
\label{ss:computation-rk-1}

In this subsection we fix a simple root $\alpha$ of $(G,T)$, and set $s:=s_\alpha$. We also choose an element $\varpi \in X_*^+(T)$ such that $\langle \varpi,\alpha \rangle = 1$ (or equivalently $s(\varpi)=\varpi-\alpha^\vee$). (Our running assumption in this section guarantees that such a coweight exists.)

\begin{lem}
\label{lem:Wak-varpi}
 We have
 \[
  \Wak_{\varpi} \cong \nabla^{\Iw}_{\st(\varpi)} \cong \nabla^{\Iw}_{\st(\varpi) s} \star_{\Iw} \nabla^{\Iw}_s, \qquad \Wak_{\varpi - \alpha^\vee} \cong \Delta^{\Iw}_s \star_{\Iw} \nabla^{\Iw}_{\st(\varpi) s}.
 \]
\end{lem}

\begin{proof}
 The equality $\Wak_{\varpi} = \nabla^{\Iw}_{\st(\varpi)}$ follows from the definitions (and the fact that $\varpi$ is dominant), and the isomorphism $\nabla^{\Iw}_{\st(\varpi) s} \star_{\Iw} \nabla^{\Iw}_s \cong \nabla^{\Iw}_{\st(\varpi)}$ from the fact that $\ell(\st(\varpi) s) = \ell(\st(\varpi))-1$ (which can be checked using~\eqref{eqn:formula-length}) and~\eqref{eqn:Delta-nabla-conv}. Next, we observe that $\varpi - \alpha^\vee=(2\varpi - \alpha^\vee)-\varpi$, with $2\varpi - \alpha^\vee \in X_*^+(T)$. Hence
 \[
  \Wak_{\varpi - \alpha^\vee} \cong \Delta^{\Iw}_{\st(-\varpi)} \star_{\Iw} \nabla^{\Iw}_{\st(2\varpi-\alpha^\vee)}.
 \]
Now as above we have $\ell(\st(\varpi) s) = \ell(\st(\varpi))-1$, hence $\Delta^{\Iw}_{\st(-\varpi)} \cong \Delta^{\Iw}_s \star_{\Iw} \Delta^{\Iw}_{s\st(-\varpi)}$. On the other hand we have $\ell(\st(2\varpi-\alpha^\vee)) = \ell(\st(2\varpi))-2= 2\ell(\st(\varpi))-2$, and $\st(2\varpi-\alpha^\vee) = (\st(\varpi) s)^2$. It follows that
\[
 \nabla^{\Iw}_{\st(2\varpi-\alpha^\vee)} \cong \nabla^{\Iw}_{\st(\varpi) s} \star_{\Iw} \nabla^{\Iw}_{\st(\varpi) s}.
\]
We deduce as desired that
\[
\Wak_{\varpi - \alpha^\vee} \cong \Delta^{\Iw}_s \star_{\Iw} \nabla^{\Iw}_{\st(\varpi) s},
\]
since $\Delta^{\Iw}_{s\st(-\varpi)} \star_{\Iw} \nabla^{\Iw}_{\st(\varpi) s} \cong \delta_{\Fl}$ by~\eqref{eqn:Delta-nabla-inverse}.
\end{proof}

We now set
\[
 X_1 := \Fl_{G,\st(\varpi)}, \quad X_2 := \Fl_{G,\st(\varpi-\alpha^\vee)}, \quad Y:= \Fl_{G,\st(\varpi) s}
\]
and $X:=X_1 \sqcup X_2 \sqcup Y$. Then $X_1$, $X_2$ and $Y$ are isomorphic to affine spaces, $X$ is a locally closed subvariety in $\Fl_G$, $X_1$ and $X_2$ are open in $X$, and $Y$ is closed in $X$. The corresponding embeddings will be denoted
\[
 j : X \hookrightarrow \Fl, \quad j_1 : X_1 \hookrightarrow X, \quad j_2 : X_2 \hookrightarrow X, \quad i : Y \to X.
\]

\begin{cor}
\label{cor:Wak-varpi-2}
 We have
 \[
  \Wak_{\varpi} = j_* \bigl( j_{1*} \underline{\bk}_ {X_1}[\ell(\st(\varpi))] \bigr), \quad \Wak_{\varpi-\alpha^\vee} = j_* \bigl( j_{2!} \underline{\bk}_ {X_2}[\ell(\st(\varpi))] \bigr).
 \]
\end{cor}

\begin{proof}
 The first isomorphism is clear from the first isomorphism in Lemma~\ref{lem:Wak-varpi}, since $j \circ j_1 = j_{\st(\varpi)}$. For the description of $\Wak_{\varpi-\alpha^\vee}$, we start with the isomorphism $\Wak_{\varpi - \alpha^\vee} \cong \Delta^{\Iw}_s \star_{\Iw} \nabla^{\Iw}_{\st(\varpi) s}$ from Lemma~\ref{lem:Wak-varpi}. Let us denote by $q : \Loop G \to \Fl_G$ the natural projection. By definition of convolution, we have
 \[
  \Delta^{\Iw}_s \star_{\Iw} \nabla^{\Iw}_{\st(\varpi) s} = f_*(\Delta^{\Iw}_s \, \widetilde{\boxtimes} \, \underline{\bk}_{Y} [\dim(Y)]),
 \]
where $f : q^{-1}(\overline{\Fl_{G,s}}) \times^{\Iw} Y \to \Fl_G$ is the morphism induced by multiplication in $\Loop G$. Now $f$ induces an isomorphism $q^{-1}(\overline{\Fl_{G,s}}) \times^{\Iw} Y \simto Y \sqcup X_2$, under which the open subvariety $q^{-1}(\Fl_{G,s}) \times^{\Iw} Y$ identifies with $X_2$. Hence $\Delta^{\Iw}_s \star^{\Iw} \nabla^{\Iw}_{\st(\varpi) s}$ is obtained from $\underline{\bk}_{X_2}[\ell(\st(\varpi))]$ by taking $!$-pushforward under the embedding $X_2 \hookrightarrow Y \sqcup X_2$, and then $*$-pushforward under the embedding $Y \sqcup X_2 \hookrightarrow \Fl$. The latter map can be written as the composition of the closed embedding $Y \sqcup X_2 \hookrightarrow X$ followed by $j$, and we deduce the desired identification.
\end{proof}

In the following lemma we denote by $\Perv_{\Gm^{\mathrm{rot}} \ltimes \Iw}(\Fl_G,\bk)$ the category of $\Gm \ltimes \Iw$-equivariant perverse sheaves on $\Fl_G$, where $\Gm$ acts on $\Iw$ and $\Fl_G$ by loop rotation.

\begin{prop}
\label{prop:Ext1-computation}
We have
\[
\Ext^1_{\Perv_{\Gm^{\mathrm{rot}} \ltimes \Iw}(\Fl_G,\bk)}(\Wak_{\varpi}, \Wak_{\varpi - \alpha^\vee}) =0.
\]
\end{prop}

\begin{proof}
Using Corollary~\ref{cor:Wak-varpi-2} and adjunction, we see that
\[
\Ext^1_{\Perv_{\Gm^{\mathrm{rot}} \ltimes \Iw}(\Fl_G,\bk)}(\Wak_{\varpi}, \Wak_{\varpi - \alpha^\vee}) \cong \Ext^1_{\Perv_{\Gm^{\mathrm{rot}} \ltimes \Iw}(X,\bk)}(j_{1*} \underline{\bk}_ {X_1}[d], j_{2!} \underline{\bk}_ {X_2}[d]),
\]
where $d:=\ell(\st(\varpi))$. We denote by $\overline{\jmath}_2 : \overline{X_2} = X_2 \sqcup Y \hookrightarrow X$ the embedding of the closure $\overline{X_2}$ (a smooth closed subvariety) in $X$, and consider the natural distinguished triangle
\[
j_{2!} \underline{\bk}_{X_2}[d] \to \overline{\jmath}_{2!} \underline{\bk}_{\overline{X_2}}[d] \to i_* \underline{\bk}_Y[d] \xrightarrow{[1]}.
\]
Here the second term is a simple perverse sheaf, which we will denote by $\IC_2$, and the third term is the shift by $1$ of a simple perverse sheaf which will be denoted $\IC_Y$. Hence this triangle induces a short exact sequence of $\Gm^{\mathrm{rot}} \ltimes \Iw$-equivariant perverse sheaves
\begin{equation}
\label{eqn:ses-X-2}
\IC_Y \hookrightarrow j_{2!} \underline{\bk}_{X_2}[d] \twoheadrightarrow \IC_2
\end{equation}
on $X$. Similarly, we have a short exact sequence of $\Gm^{\mathrm{rot}} \ltimes \Iw$-equivariant perverse sheaves
\begin{equation}
\label{eqn:ses-X-1}
\IC_Y \hookrightarrow j_{1!} \underline{\bk}_{X_1}[d] \twoheadrightarrow \IC_1
\end{equation}
where $\IC_1$ is the pushforward of the constant sheaf on $X_1 \sqcup Y$ shifted by $d$. Applying Verdier duality to~\eqref{eqn:ses-X-2}--\eqref{eqn:ses-X-1} we deduce short exact sequences of $\Gm^{\mathrm{rot}} \ltimes \Iw$-equivariant perverse sheaves
\begin{equation}
\label{eqn:ses-X-2'}
\IC_2 \hookrightarrow  j_{2*} \underline{\bk}_{X_2}[d] \twoheadrightarrow \IC_Y, \qquad \IC_1 \hookrightarrow  j_{1*} \underline{\bk}_{X_1}[d] \twoheadrightarrow \IC_Y.
\end{equation}

Applying the functor $\Hom_{\Perv_{\Gm^{\mathrm{rot}} \ltimes \Iw}(\Fl_G,\bk)}(-,j_{2!} \underline{\bk}_{X_2}[d])$ to the second exact sequence in~\eqref{eqn:ses-X-2'} we obtain an exact sequence
\begin{multline}
\label{eqn:Ext1-computation-long}
\Ext^1(\IC_Y, j_{2!} \underline{\bk}_{X_2}[d]) \to
\Ext^1(j_{1*} \underline{\bk}_{X_1}[d], j_{2!} \underline{\bk}_{X_2}[d]) \\ \to \Ext^1(\IC_1, j_{2!} \underline{\bk}_{X_2}[d]) \to
\Hom_{\Db_{\Gm^{\mathrm{rot}} \ltimes \Iw}(X,\bk)}(\IC_Y, j_{2!} \underline{\bk}_{X_2}[d+2]).
\end{multline}
Hence to conclude, it suffices to prove that
\begin{gather}
\label{eqn:Ext1-computation}
\Ext^1(\IC_Y, j_{2!} \underline{\bk}_{X_2}[d])=0; \\
\label{eqn:Ext1-computation-2}
\text{the third map in~\eqref{eqn:Ext1-computation-long} is injective.}
\end{gather}

For this,
consider (a portion of) the long exact sequence obtained by applying the functor $\Hom(\IC_Y,-)$ to~\eqref{eqn:ses-X-2}:
\begin{multline}
\label{eqn:proof-Ext1-computation}
\Ext^1(\IC_Y,\IC_Y) \to \Ext^1(\IC_Y, j_{2!} \underline{\bk}_{X_1}[d]) \to \Ext^1(\IC_Y,\IC_2) \\
\to \Hom_{\Db_{\Gm^{\mathrm{rot}} \ltimes \Iw}(X,\bk)}(\IC_Y,\IC_Y[2]) \to \Hom_{\Db_{\Gm^{\mathrm{rot}} \ltimes \Iw}(X,\bk)}(\IC_Y, j_{2!} \underline{\bk}_{X_2}[d+2]) \\
\to \Hom_{\Db_{\Gm^{\mathrm{rot}} \ltimes \Iw}(X,\bk)}(\IC_Y, \IC_2[2]).
\end{multline}
The first, resp.~fourth, term in this sequence identifies with $\mathsf{H}^1_{\Gm^{\mathrm{rot}} \ltimes \Iw}(Y,\bk)=0$, resp.~with $\mathsf{H}^2_{\Gm^{\mathrm{rot}} \ltimes \Iw}(Y,\bk) = \bk \otimes_\Z X^*(\Gm \times T)$. It is also easy to check using adjunction that we have
\[
\Ext^1(\IC_Y,\IC_2) \cong \mathsf{H}^0_{\Gm^{\mathrm{rot}} \ltimes \Iw}(Y,\bk)=\bk
\]
and
\[
\Hom_{\Db_{\Gm^{\mathrm{rot}} \ltimes I}(X,\bk)}(\IC_Y, \IC_2[2]) \cong \mathsf{H}^1_{\Gm^{\mathrm{rot}} \ltimes \Iw}(Y,\bk)=0.
\]

For any $\beta \in \mathfrak{R}^+$ we have a root subgroup $U_\beta \subset G$, and for $m \in \Z$ we denote by $U_{\beta,m} \subset \Loop G$ the subgroup which identifies, for any choice of isomorphism $u_\beta : \Ga \simto U_\beta$, with the image of $t \mapsto u_\beta(tz^m)$. We set
\[
 J:= \prod_{\substack{\beta \in \mathfrak{R}^+ \smallsetminus \{\alpha\} \\ 0 < n \leq \langle \varpi,\beta \rangle}} U_{\beta,n},
\]
where the product is taken in an arbitrary fixed order. Then the map $j \mapsto j z^{\varpi} \dot{s}I/I$ induces an isomorphism $J \simto Y$, and the map $(u,j) \mapsto ujz^{\varpi} \dot{s}I/I$ induces an isomorphism between $U_{\alpha} \times J$ and an open neighborhood of $Y$ in $Y \sqcup X_2$. 
Hence, arguing as in~\cite[Proof of Lemma~6.5]{bezr}, we see that the image of the map $\Ext^1(\IC_Y,\IC_2) \to \Hom_{\Db_{\Gm^{\mathrm{rot}} \ltimes I}(X,\bk)}(\IC_Y,\IC_Y[2])$ in~\eqref{eqn:proof-Ext1-computation} is the line generated by $1 \otimes \alpha$.
It follows that the third map in~\eqref{eqn:proof-Ext1-computation} is injective, proving~\eqref{eqn:Ext1-computation}, and providing moreover an isomorphism
\begin{equation}
\label{eqn:Ext2}
 \Hom_{\Db_{\Gm^{\mathrm{rot}} \rtimes I}(X,\bk)}(\IC_Y, j_{2!} \underline{\bk}_{X_1}[d][2]) \cong (\bk \otimes_\Z X^*(\Gm \times T)) / \bk \cdot (1 \otimes \alpha).
\end{equation}

Now, let us consider (a portion of) the long exact sequence obtained by applying $\Hom(\IC_1,-)$ to~\eqref{eqn:ses-X-2}:
\begin{multline*}
\Hom(\IC_1,\IC_2) \to \Ext^1(\IC_1,\IC_Y) \to \\
\Ext^1(\IC_1,j_{2!} \underline{\bk}_{X_2}[d]) \to \Ext^1(\IC_1, \IC_2).
\end{multline*}
Here the first term vanishes. The fourth term also vanishes: in fact this follows from the observation that any extension of $\IC_1$ by $\IC_2$ must be the intermediate extension of its restriction to the open subset $X_1 \sqcup X_2$ (because so are $\IC_1$ and $\IC_2$); hence this extension must split. These considerations show that the natural morphism $\Ext^1(\IC_1,\IC_Y) \to \Ext^1(\IC_1,j_{2!} \underline{\bk}_{X_2}[d])$ is an isomorphism, which allows to identify the third map in~\eqref{eqn:Ext1-computation-long} with the composition
\begin{multline*}
\Ext^1(\IC_1,\IC_Y) \to \Hom_{\Db_{\Gm^{\mathrm{rot}} \ltimes \Iw}(X,\bk)}(\IC_Y, \IC_Y[2]) \\
\to \Hom_{\Db_{\Gm^{\mathrm{rot}} \ltimes \Iw}(X,\bk)}(\IC_Y, j_{2!} \underline{\bk}_{X_2}[d][2]).
\end{multline*}
As above the first term is $1$-dimensional, the middle term identifies with $\bk \otimes_\Z X^*(\Gm \times T)$, and the first map is injective with image $\bk \cdot (1 \otimes (\alpha-\hbar))$, where $\hbar$ is the tautological character of $\Gm$. (To see this, one observes that the map $(j,u) \mapsto j z^{\varpi} u\dot{s}I/I$ induces an isomorphism between $J \times U_{-\alpha}$ and an open neighborhood of $Y$ in $Y \sqcup X_1$.) In view of~\eqref{eqn:Ext2} we deduce that our composition is indeed injective, finishing the proof of~\eqref{eqn:Ext1-computation-2}, hence of the proposition.
\end{proof}

\subsection{Proof of Proposition~\ref{prop:components-u}}
\label{ss:proof-u}

We can now give the proof of Proposition~\ref{prop:components-u}. 
We fix a simple root $\alpha$, with corresponding simple reflection $s$. 
We also fix an element $\varpi \in X_*(T)=X^*(T^\vee_\bk)$ such that $\langle \varpi,\alpha \rangle = 1$ and $\langle \varpi,\beta \rangle = 0$ for any simple root $\beta$ different from $\alpha$, and use it as our choice of cocharacter denoted similarly in~\S\ref{ss:computation-rk-1}. To simplify notation, we set $r:=\langle \varpi, 2\rho \rangle = \ell(\varpi)$.

Let $V:=\mathrm{Ind}_{\widetilde{B}^\vee}^{G_\bk^\vee}(\varpi)$ (where $\varpi$ is identified with a character of $A^\vee = \widetilde{B}^\vee / \widetilde{U}^\vee$ via the canonical isomorphism $T_\bk^\vee \cong A^\vee$, see~\S\ref{ss:regularity-statement}). For the filtration on $V$ constructed in ~\S\ref{ss:regularity-statement} we have $V=F^r(V)$, $V/F^{r-2}(V)$ is one-dimensional (with weight $\varpi$ for the action of $A^\vee$), and so is $F^{r-2}(V)/F^{r-4}(V)$ (with weight $\varpi-\alpha^\vee$ for the action of $A^\vee$). This shows that the perverse sheaf
\[
\scF := \sZ(V)/F^{r-4} \sZ(V)
\]
fits in an exact sequence of $\Iw$-equivariant perverse sheaves
\begin{equation}
\label{eqn:extension-Wak-varpi}
\Wak_{\varpi-\alpha^\vee} \hookrightarrow \scF \twoheadrightarrow \Wak_\varpi.
\end{equation}

\begin{lem}
\label{lem:extension-non-split}
The extension~\eqref{eqn:extension-Wak-varpi} is nonsplit.
\end{lem}

\begin{proof}
Let $\wt(\varpi)$ be the set of $T_\bk^\vee$-weights of $V$. Then since $\sZ(V)$ admits a filtration whose subquotients are of the form $\Wak_\lambda$ with $\lambda \in \wt(\varpi)$ (see~\S\ref{ss:filtration-central-sheaves}), in view of Property~\eqref{it:Wak-support} in~\S\ref{ss:Wakimoto-sheaves} its support must be
\[
\bigcup_{\lambda \in \wt(\varpi)} \overline{\Fl_{G,\st(\lambda)}}.
\]
It is clear that the subvariety $X$ of~\S\ref{ss:computation-rk-1} is contained (and open) in the closed subvariety $\overline{\Fl_{G,\st(\varpi)}} \cup \overline{\Fl_{G,\st(\varpi-\alpha^\vee)}}$. Moreover, $X$ does not intersect any subvariety $\overline{\Fl_{G,\st(\lambda)}}$ with $\lambda \in \wt(\varpi) \smallsetminus \{\varpi, \varpi-\alpha^\vee\}$; in fact this is clear for dimension reasons if $\lambda \notin \Wf \cdot \varpi$, and if $\lambda \in \Wf \cdot \varpi \smallsetminus \{\varpi, \varpi-\alpha^\vee\}$ then $\st(\varpi)s$ is not smaller than $\st(\lambda)$ in the Bruhat order.

From this observation we deduce that the restriction of $\sZ(V)$ to $X$ coincides with that of $\scF$, and is a perverse sheaf. Restricting~\eqref{eqn:extension-Wak-varpi} to $X$ and using Corollary~\ref{cor:Wak-varpi-2} we see moreover that we have an exact sequence
\begin{equation}
\label{eqn:extension-Wak-varpi-rest}
(j_2)_! \underline{\bk}_{X_2}[r] \hookrightarrow \sZ(V)_{|X} \twoheadrightarrow (j_1)_* \underline{\bk}_{X_1}[r].
\end{equation}

Now, similar considerations using the ``opposite'' Wakimoto sheaves (see Remark~\ref{rmk:Wakop} and Remark~\ref{rmk:Wakop-2}) show that $\sZ(V)$ has a canonical subobject $\scF'$ such that $\sZ(V)_{|X} = \scF'_{|X}$ and which fits in a short exact sequence $\Wakop_{\varpi} \hookrightarrow \scF' \twoheadrightarrow \Wakop_{\varpi-\alpha^\vee}$. Restricting to $X$ and using the ``opposite'' version of Corollary~\ref{cor:Wak-varpi-2} we obtain a short exact sequence
\[
(j_1)_! \underline{\bk}_{X_1}[r] \hookrightarrow \sZ(V)_{|X} \twoheadrightarrow (j_2)_* \underline{\bk}_{X_2}[r].
\]
The existence of this exact sequence forces~\eqref{eqn:extension-Wak-varpi-rest} (and hence~\eqref{eqn:extension-Wak-varpi}) to be nonsplit. (For instance, if~\eqref{eqn:extension-Wak-varpi-rest} were split there would exist no injective map $(j_1)_! \underline{\bk}_{X_1}[r] \hookrightarrow \sZ(V)_{|X}$.)
\end{proof}

Now we consider the action of $\su$ on the $2$-dimensional $\widetilde{B}^\vee$-module $V/F^{r-4}(V)$. (Since $\su$ belongs to $\widetilde{B}^\vee$, it stabilizes $F^{r-4}(V)$, hence indeed induces an endomorphism of this quotient.) It is clear from the description of this module that this action coincides with the action of the component $\su_{-\alpha}$. Hence to conclude the proof it suffices to prove that this action is nontrivial, or in other words that the endomorphism of $\Pi_{\Iw,\Iw}^0(\scF)$ induced by $\sm_V^0$ is nontrivial. In view of~\eqref{eqn:sm-monodromy}, for this it suffices to prove that $\Pi_{\Iw,\Iw}^0(\mathfrak{M}_\scF)$ is nontrivial.

First, we claim that $\mathfrak{M}_\scF$ is nontrivial. In fact, in view of~\cite[Lemma~2.6]{bezr} (see also~\cite[\S 10.1]{bezr} for a discussion of the \'etale version) this is equivalent to saying that $\scF$ is not equivariant under the action of $\Gm^{\mathrm{rot}}$. However, by construction $\scF$ is $\Iw$-equivariant. If it were $\Gm^{\mathrm{rot}}$-equivariant, then it would be $\Gm^{\mathrm{rot}} \ltimes \Iw$-equivariant, which is impossible by Proposition~\ref{prop:Ext1-computation} and Lemma~\ref{lem:extension-non-split}.

Now that this claim is established, we know that the nilpotent endomorphism $\mathfrak{M}_\scF - \id_\scF$ is nonzero. This endomorphism stabilizes the filtration~\eqref{eqn:extension-Wak-varpi}, and induces the zero endomorphism on the associated graded. It therefore factors as a composition
\[
\scF \twoheadrightarrow \Wak_\varpi \xrightarrow{f} \Wak_{\varpi-\alpha^\vee} \hookrightarrow \scF,
\]
where $f \neq 0$. Now we have $\Wak_\varpi = \nabla^{\Iw}_\varpi$, whose top is isomorphic to $\IC_\omega$ for some $\omega \in \Omega$ by Lemma~\ref{lem:socle-Deltas}. Hence the image of $f$ admits $\IC_\omega$ as a composition factor. Since $\IC_\omega$ is not killed by $\Pi_{\Iw,\Iw}^0$ it follows that $\Pi_{\Iw,\Iw}^0(f) \neq 0$, so that $\Pi_{\Iw,\Iw}^0(\mathfrak{M}_\scF) \neq \id_{\Pi_{\Iw,\Iw}^0(\scF)}$, which finishes the proof.

\section{Antispherical and Iwahori--Whittaker categories}
\label{sec:asp-IW}



In this section we drop the assumption that $X^*(T)/\Z\mathfrak{R}$ is free. Instead, from~\S\ref{ss:asp-IW-statement} on we will assume that either $X^*(T)/\Z\mathfrak{R}$ is free or $\ell$ is very good for $G$. (The definitions and statement in~\S\ref{ss:IW-Perv} will later be used without this assumption.)

\subsection{Iwahori--Whittaker perverse sheaves}
\label{ss:IW-Perv}

We denote by $B^+ \subset G$ the Borel subgroup containing $T$ and opposite to $B$, and by $\Iw^+ \subset \Loop^+ G$ the corresponding Iwahori subgroup. Let also $\Iwu^+$ be the pro-unipotent radical of $\Iw^+$, i.e.~the inverse image of the unipotent radical $U^+$ of $B^+$ under the natural map $\Iw^+ \to B^+$. We choose a nondegenerate character $\chi' : U^+ \to \Ga$ (i.e.~a group morphism nontrivial on each simple root subgroup), and denote by $\chi$ its composition with the projection $\Iwu^+ \to U^+$. Finally we choose an Artin--Schreier $\bk$-local system $\scL_{\mathrm{AS}}$ on $\Ga$ (see e.g.~\cite{bezr}), and consider the
derived category $\sfD_{\IW,\Iw} = \Db_{\IW}(\Fl_G,\bk)$ of Iwahori--Whit\-taker sheaves on $\Fl_G$, i.e.~ the $(\Iwu^+, \chi^*(\scL_{\mathrm{AS}}))$-equivariant derived category of $\bk$-sheaves on $\Fl_G$. (As for the $\Loop^+ G$-equivariant derived category of $\Gr_G$ or the $\Iw$-equivariant derived category of $\Fl_G$, the ``true'' definition of this category requires restricting to finite unions of orbits, and considering a finite-type quotient of $\Iwu^+$ through which the action factors; this procedure is standard, and will not be discussed.) This category is endowed with the perverse t-structure, whose heart will be denoted $\sfP_{\IW, \Iw}$.

The $\Iwu^+$-orbits on $\Fl_G$ are parametrized in the standard way by $W$. Those which support a nonzero $(\Iwu^+, \chi^*(\scL_{\mathrm{AS}}))$-equivariant local system are those corresponding to elements in the subset $\fW \subset W$ of elements $w$ which are minimal in $\Wf w$. For any $\lambda \in X_*(T)$ we will denote by $w_\lambda$ the minimal length element in $\Wf \cdot \st(\lambda)$, by $\Fl_{G,\lambda}^{\IW}$ the corresponding orbit, and by $\scL_{\chi,\lambda}$ the unique rank-$1$ $(\Iwu^+, \chi^*(\scL_{\mathrm{AS}}))$-equivariant local system on $\Fl^\IW_{G,\lambda}$. We also denote by $j^\IW_\lambda : \Fl_{G,\lambda}^\IW \to \Fl_G$ the embedding, and set
\[
 \Delta_\lambda^\IW := (j_{\lambda}^\IW)_! \scL_{\chi,\lambda} [\dim(\Fl^\IW_{G,\lambda})], \quad \nabla_\lambda^\IW := (j_{\lambda}^\IW)_* \scL_{\chi,\lambda} [\dim(\Fl^\IW_{G,\lambda})].
\]
(These objects are perverse sheaves because $j^\IW_\lambda$ is an affine embedding.) As usual we have
\begin{equation}
\label{eqn:Hom-DN-IW}
\Hom_{\Db_\IW(\Fl_G,\bk)}(\Delta_\lambda^\IW, \nabla_\mu^\IW[n]) \cong 
\begin{cases}
\bk & \text{if $\lambda=\mu$ and $n=0$;} \\
0 & \text{otherwise,}
\end{cases}
\end{equation}
and the image of any nonzero morphism $\Delta^\IW_\lambda \to \nabla^\IW_\lambda$ is simple; this simple object will be denoted $\IC_\lambda^\IW$. It is well known also that the category $\sfP_{\IW,\Iw}$ has a natural structure of highest weight category, with weight poset $X_*(T)$ (for the order $\leq^\IW$ defined by $\lambda \leq^\IW \mu$ iff $\Fl^\IW_{G,\lambda} \subset \overline{\Fl^\IW_{G,\mu}}$) and standard, resp.~costandard, objects $(\Delta^\IW_\lambda : \lambda \in X_*(T))$, resp.~$(\nabla^\IW_\lambda : \lambda \in X_*(T))$.

In Section~\ref{sec:central-tilting} we will need some properties of the associated order on $X_*(T)$. Here, for $\mu \in X_*(T)$ we denote by $\mu^+$ the unique dominant $\Wf$-conjugate of $\mu$, and for $\lambda \in X_*(T)$ we set
\[
 \wt(\lambda)=\{\mu \in X_*(T) \mid  \mu^+\preceq\lambda^+\}, 
\]
where $\preceq$ is as in~\S\ref{ss:Wakimoto-sheaves}.
In other words, $\wt(\lambda)$ is the set of $T_\bk^\vee$-weights of the Weyl module for $G^\vee_\bk$ of highest weight $\lambda^+$ (so that this notation extends the notation used already in the proof of Lemma~\ref{lem:extension-non-split}). The other standard description of these weights shows that $\wt(\lambda)$ coincides with the set denoted $\mathsf{conv}(\lambda)$ in~\cite{prinblock}, i.e.~the intersection of the convex hull of $\Wf\lambda$ with $\lambda + \Z \mathfrak{R}^\vee$.

\begin{lem}\phantomsection
\label{lem:inclusion-orbits}
\begin{enumerate}
\item
\label{it:inclusion-orbits-1}
If $\lambda,\mu \in X_*(T)$ and $\Fl^\IW_{G,\lambda} \subset \overline{\Fl^\IW_{G,\mu}}$, then $\lambda \in \wt(\mu)$.
\item
\label{it:inclusion-orbits-2}
If $\lambda \in X_*^+(T)$ and $\mu \in \Wf(\lambda)$, then $\Fl^{\IW}_{G,\mu} \subset \overline{\Fl^\IW_{G,\lambda}}$.
\end{enumerate}
\end{lem}

\begin{proof}
It is not difficult to see that for any $\nu \in X_*(T)$ we have 
\[
\overline{\Fl^\IW_{G,\nu}}=\overline{\Fl_{G,w_0 w_\nu}},
\]
where (as in~\S\ref{ss:filtration-central-sheaves}) $w_0$ is the longest element in $\Wf$.
Hence $\Fl^\IW_{G,\lambda} \subset \overline{\Fl^\IW_{G,\mu}}$ iff $w_0 w_\lambda \leq w_0 w_\mu$ in the Bruhat order. In view of~\cite[Lemma~2.2]{douglass}, this is equivalent to the property that $w_\lambda \leq w_\mu$. Hence~\eqref{it:inclusion-orbits-1} follows from~\cite[Lemma~9.12(3)]{prinblock}, and~\eqref{it:inclusion-orbits-2} from the fact that $w_\mu \leq w_\lambda=\st(\lambda)$ in this case (see e.g.~\cite[Lemma~2.4]{mr}).
\end{proof}

\subsection{Statement}
\label{ss:asp-IW-statement}

As explained at the beginning of the section, from now on in this section we assume that either $X^*(T)/\Z\mathfrak{R}$ is free or $\ell$ is very good for $G$.

Note that since the closure $\overline{\Fl^\IW_{G,0}}$ does not contain any orbit $\Fl^\IW_{G,\lambda}$ with $\lambda \neq 0$, the natural morphism $\Delta^\IW_0 \to \nabla^\IW_0$ is an isomorphism (and this complex coincides with $\IC^\IW_0$). We consider the functor
\begin{equation}
\label{eqn:def-AvIW}
 \Av_{\IW} : \sfD_{\Iw,\Iw} \to \sfD_{\IW,\Iw}
\end{equation}
defined by
\[
 \Av_\IW(\scF)=\Delta_0^\IW \star_{\Iw} \scF.
\]

We will denote by $\sfP_{\Iw,\Iw}^{\asp}$ the Serre quotient of the category $\sfP_{\Iw,\Iw}$ by the Serre subcategory generated by the simple objects $\IC_w$ with $w \notin \fW$. The associated quotient functor will be denoted $\Pi_{\Iw,\Iw}^{\asp} : \sfP_{\Iw,\Iw} \to \sfP_{\Iw,\Iw}^\asp$. (Here, ``$\asp$'' stands for ``anti-spherical.'')

The main result of the present section is the following statement. 

\begin{thm}\phantomsection
\label{thm:Av-ff}
\begin{enumerate}
 \item
 \label{it:Av-t-exact}
 The functor $\Av_{\IW}$ is t-exact for the perverse t-structures on $\sfD_{\Iw,\Iw}$ and $\sfD_{\IW,\Iw}$.
 \item 
 \label{it:Av-ff}
 The functor
 \[
 \Av_{\IW} : \sfP_{\Iw,\Iw} \to \sfP_{\IW,\Iw}
 \]
 factors through a fully-faithful functor
 \[
 \Av_\IW^\asp : \sfP_{\Iw,\Iw}^{\asp} \to \sfP_{\IW,\Iw}.
 \]
 Moreover, the essential image of this functor is closed under subquotients.
\end{enumerate}
\end{thm}

\begin{rmk}
We will see later (see~\S\ref{ss:comparison}) that, under appropriate assumptions,
the functor in Theorem~\ref{thm:Av-ff}\eqref{it:Av-ff} is an equivalence of categories.
\end{rmk}

The proof of 
Theorem~\ref{thm:Av-ff} is very close to the corresponding proof in~\cite[\S 2]{ab}, which is repeated with more details in~\cite[Chap.~6]{ar-book}. We will not repeat the proofs that can simply be copied from these references.

\subsection{Exactness}

We start by proving Theorem~\ref{thm:Av-ff}\eqref{it:Av-t-exact}.


The following lemma can be proved as in~\cite[Lemma~6.4.4 \& Lemma~6.4.5]{ar-book}.

\begin{lem}
\phantomsection
\label{lem:Av-D-N}
\begin{enumerate}
 \item 
 \label{it:Av-kills-IC}
 For $w \in W$ we have $\Av_{\IW}(\IC_w)=0$ unless $w \in \fW$.
 \item
 \label{it:Av-D-N}
 For any $w \in W$ we have
\[
\Av_{\IW}(\Delta^{\Iw}_w) \cong \Delta^{\IW}_\lambda \quad \text{and} \quad \Av_{\IW}(\nabla^{\Iw}_w) \cong \nabla^{\IW}_\lambda,
\]
where $\lambda \in X_*(T)$ is the unique element such that $\Wf \cdot w = \Wf \cdot w_\lambda$.
\end{enumerate}
\end{lem}

As in~\cite[Corollary~6.4.7]{ar-book} one deduces
the following corollary, which in particular establishes Theorem~\ref{thm:Av-ff}\eqref{it:Av-t-exact}.

\begin{cor}
\label{cor:Av-exact}
The functor $\Av_\IW$ is t-exact. Moreover, for any $\lambda \in \bX$ we have
\[
\Av_{\IW}(\IC_{w_\lambda}) \cong \IC_\lambda^\IW.
\]
\end{cor}

\subsection{An adjoint functor}
\label{ss:adjoint-IW}

The functor $\Av_\IW$ admits another description. Namely, we will denote by $\Iw_0$ the intersection of $\Iw$ and $\Iwu^+$, or in other words the kernel of the morphism $\mathrm{ev} : \Loop^+ G \to G$ from~\S\ref{ss:Fl}. We can then consider the associated equivariant derived category $\Db_{\Iw_0}(\Fl_G,\bk)$, and the functor
\[
 \lInd_{\Iw_0}^{\IW} : \Db_{\Iw_0}(\Fl_G,\bk) \to \Db_{\IW}(\Fl_G,\bk)
\]
sending a complex $\mathscr{F}$ to $a_!(\chi^*(\scL_{\mathrm{AS}}) \, \widetilde{\boxtimes} \, \mathscr{F})[\dim(U^+)]$, where $a : \Iwu^+ \times^{\Iw_0} \Fl_G \to \Fl_G$ is the action morphism. It is not difficult to see that we have
\[
 \Av_\IW = \lInd_{\Iw_0}^{\IW} \circ \For^{\Iw}_{\Iw_0},
\]
where $\For^{\Iw}_{\Iw_0} : \sfD_{\Iw,\Iw} \to \Db_{\Iw_0}(\Fl_G,\bk)$ is the forgetful functor, and moreover that the functor $\lInd_{\Iw_0}^{\IW}$ is left adjoint to $\For^{\IW}_{\Iw_0} [\dim(U^+)]$, where $\For^{\IW}_{\Iw_0} : \Db_\IW(\Fl_G,\bk) \to \Db_{\Iw_0}(\Fl_G,\bk)$ is the forgetful functor (see e.g.~\cite[Lemma~A.3]{modrap1} for similar considerations).

Now we also consider
the induction functor
\[
\rInd_{\Iw_0}^{\Iw} : \Db_{\Iw_0}(\Fl_G,\bk) \to \sfD_{\Iw,\Iw}
\]
sending a complex $\scF$ to $b_*(\underline{\bk} \, \widetilde{\boxtimes} \, \scF)[\dim(B)]$, where $b : \Iw \times^{\Iw_0} \Fl_G \to \Fl_G$ is the action morphism.
This functor is right adjoint to the functor $\For^{\Iw}_{\Iw_0}[-\dim(B)]$. Therefore, the functor
\[
 \rInd_{\Iw_0}^{\Iw} \circ \For^{\IW}_{\Iw_0} [-\dim(T)]
\]
is right adjoint to the functor $\Av_\IW$.

As seen in Corollary~\ref{cor:Av-exact} the functor $\Av_\IW$ is t-exact; therefore its right adjoint $\rInd_{\Iw_0}^{\Iw} \circ \For^{\IW}_{\Iw_0} [-\dim(T)]$ is left t-exact, and the functor
\begin{equation}
\label{eqn:def-F}
 \mathsf{F} : \sfP_{\IW,\Iw} \to \sfP_{\Iw,\Iw}
\end{equation}
defined by
\[
\mathsf{F}(\mathscr{F}) = \pH^0 \bigl( \rInd_{\Iw_0}^{\Iw} \circ \For^{\IW}_{\Iw_0}(\mathscr{F}) [-\dim(T)] \bigr)
\]
is right adjoint to the functor $\Av_\IW : \sfP_{\Iw,\Iw} \to \sfP_{\IW,\Iw}$.

\subsection{Some properties of tilting objects}
\label{ss:properties-tilting}

As explained in~\S\ref{ss:IW-Perv} the category $\sfP_{\IW,\Iw}$ has a natural highest weight structure, with weight poset $X_*(T)$ (which we have identified with $\fW$). In such a category we can consider the tilting objects, and the general theory of highest weight categories (see e.g.~\cite[\S 7]{riche-hab}) shows that the set of isomorphism classes of indecomposable tilting objects in $\sfP_{\IW,\Iw}$ is in a natural bijection with $X_*(T)$; the object associated with $\lambda$ will be denoted $\mathscr{T}^\IW_\lambda$. In this subsection we prove some properties of these objects which will not play any role in the proof of Theorem~\ref{thm:Av-ff}, but which will be needed in Section~\ref{sec:proof-main}
below.

We begin with the following claim. Here we denote by $\sfP_{\IW,\Iw}^0$ the quotient of the abelian category $\sfP_{\IW,\Iw}$ by the Serre subcategory generated by the objects of the form $\IC^\IW_\lambda$ with $\ell(w_\lambda)>0$; the associated quotient functor will be denoted $\Pi^0_{\IW,\Iw} : \sfP_{\IW,\Iw} \to \sfP_{\IW,\Iw}^0$.

\begin{lem}
\label{lem:Hom-tiltings-IW}
For any tilting objects $\mathscr{T}$ and $\mathscr{T}'$ in $\sfP_{\IW,\Iw}$, the quotient functor $\Pi^0_{\IW,\Iw}$ induces an isomorphism
\[
\Hom_{\sfP_{\IW,\Iw}}(\mathscr{T},\mathscr{T}') \simto \Hom_{\sfP_{\IW,\Iw}^0}(\Pi^0_{\IW,\Iw}(\mathscr{T}), \Pi^0_{\IW,\Iw}(\mathscr{T}')).
\]
\end{lem}

\begin{proof}
The idea of this proof is copied from~\cite[\S 2.1]{bbm}. From the definition of morphism spaces in a quotient category, the lemma will follow if we prove that the socle, resp.~top, of each object $\Delta^\IW_\lambda$, resp.~$\nabla^\IW_\lambda$, is a sum of objects $\IC^\IW_\mu$ with $\ell(w_\mu)=0$. Here it suffices to prove the claim for $\nabla^\IW_\lambda$; the other one will follow by Verdier duality. If we choose $w \in W$ such that $\Wf \cdot w = \Wf \cdot w_\lambda$, then for any $\mu \in X_*(T)$ we have
\begin{multline*}
\Hom_{\sfP_{\IW,\Iw}}(\nabla^\IW_\lambda, \IC^\IW_\mu) \cong \Hom_{\sfP_{\IW,\Iw}}(\Av_\IW(\nabla^{\Iw}_w), \IC^\IW_\mu) \\
\cong \Hom_{\sfP_{\Iw,\Iw}}(\nabla^{\Iw}_w, \mathsf{F}(\IC^\IW_\mu))
\end{multline*}
by Lemma~\ref{lem:Av-D-N}\eqref{it:Av-D-N} and adjunction. By Lemma~\ref{lem:socle-Deltas}, the top of $\nabla^{\Iw}_w$ is $\IC^{\Iw}_\omega$, where $\omega \in \Omega$ is the unique element such that $\omega W^{\mathrm{Cox}} = wW^{\mathrm{Cox}}$. If $\ell(w_\mu)>0$, and if $s \in S$ is such that $\ell(w_\mu s)<\ell(w_\mu)$, then $w_\mu s \in \fW$, and $\IC^\IW_\mu$ is the shifted pullback of a simple perverse sheaf on the partial affine flag variety $\Fl_G^s$ associated with $s$. Then $\mathsf{F}(\IC^\IW_\mu)$ also the shifted pullback of a perverse sheaf on $\Fl_G^s$, which shows that 
\[
\Hom_{\sfP_{\Iw,\Iw}}(\nabla^{\Iw}_w, \mathsf{F}(\IC^\IW_\mu))=0
\]
and finishes the proof.
\end{proof}

The other result we will need is the following.

\begin{lem}
\label{lem:tilting-PIW}
 For any $\lambda \in X_*(T)$, the object $\mathscr{T}^\IW_\lambda$ is isomorphic to a subobject of a direct sum of objects of the form $\mathscr{T}^\IW_\mu$ with $\mu \in X_*^+(T)$.
\end{lem}

\begin{proof}
We denote by $\Iwu$ the pro-unipotent radical of $\Iw$, by $\Perv_{\Iwu}(\Fl_G,\bk)$ the category of $\Iwu$-equivariant perverse sheaves on $\Fl_G$ (or equivalently the category of perverse sheaves which are constructible with respect to the stratification by $\Iw$-orbits), and set for $w \in W$
\[
 \Delta^{\Iwu}_w:=\For^{\Iw}_{\Iwu}(\Delta^{\Iw}_w), \quad \nabla^{\Iwu}_w:=\For^{\Iw}_{\Iwu}(\nabla^{\Iw}_w),
\]
where $\For^{\Iw}_{\Iwu} : \sfD_{\Iw,\Iw} \to \Db_{\Iwu}(\Fl_G,\bk)$ is the forgetful functor. Then the category $\Perv_{\Iwu}(\Fl_G,\bk)$ admits a structure of highest weight category with weight poset $W$ (for the Bruhat order) and standard, resp.~costandard, objects the objects $\Delta^{\Iwu}_w$ with $w \in W$, resp.~$\nabla^{\Iwu}_w$ with $w \in W$.

The setup of~\S\ref{ss:adjoint-IW} allows us to construct a functor $\Av_\IW' : \Db_{\Iwu}(\Fl_G,\bk) \to \Db_\IW(\Fl_G,\bk)$ such that $\Av_\IW = \Av_{\IW}' \circ \For^{\Iw}_{\Iwu}$. Then by Lemma~\ref{lem:Av-D-N}\eqref{it:Av-D-N},
for any $w \in W$ we have
\[
\Av_{\IW}'(\Delta^{\Iwu}_w) \cong \Delta^{\IW}_\lambda \quad \text{and} \quad \Av_{\IW}'(\nabla^{\Iwu}_w) \cong \nabla^{\IW}_\lambda,
\]
where $\lambda \in X_*(T)$ is the unique element such that $\Wf \cdot w = \Wf \cdot w_\lambda$. These formulas show that $\Av_\IW'$ is t-exact, and sends tilting objects in $\Perv_{\Iwu}(\Fl_G,\bk)$ to tilting objects in $\sfP_{\IW,\Iw}$. Moreover, if we denote by $\mathscr{T}^{\Iwu}_w$ the indecomposable tilting object in $\Perv_{\Iwu}(\Fl_G,\bk)$ associated with $w \in W$, then one can check that the indecomposable tilting object in $\sfP_{\IW,\Iw}$ associated with $\lambda \in X_*(T)$ is a direct summand in $\Av_\IW'(\mathscr{T}^{\Iwu}_{w_\lambda})$. Therefore, to conclude we only have to prove that any indecomposable tilting object in $\Perv_{\Iwu}(\Fl_G,\bk)$ is isomorphic to a subobject of a tilting object whose image under $\Av_\IW'$ is a direct sum of tilting objects associated with dominant weights.

Let us fix an indecomposable tilting object $\mathscr{T}$ in $\Perv_{\Iwu}(\Fl_G,\bk)$.
Denote by $\tFl_G$ the ``extended'' affine flag variety, i.e.~the ind-scheme representing the fppf quotient $(\Loop G/\Iwu)_{\mathrm{fppf}}$; then we have a $T$-bundle $\pi : \tFl_G \to \Fl_G$.
Recall the ``free pro-unipotently monodromic'' tilting perverse sheaves on $\tFl_G$, see~\cite{by} or~\cite[\S 5.4]{br}. By~\cite[Proposition~5.12]{br} (see also~\cite[Lemma~A.7.3]{by} for a closely related statement) there exists an indecomposable free pro-unipotently monodromic tilting perverse sheaf $\widehat{\mathscr{T}}$ such that $\mathscr{T} = \pi_! \widehat{\mathscr{T}}[\dim(T)]$. Now we have a convolution action $\hatstar$ of free pro-unipotently monodromic tilting perverse sheaves on the category $\Db_{\Iwu}(\Fl_G,\bk)$ (see~\cite[\S 7.3]{br} for a discussion of a similar construction in the context of usual flag varieties), and the considerations in the proof of~\cite[Lemma~7.8]{br} show that for any free pro-unipotently monodromic tilting perverse sheaf $\widehat{\mathscr{T}}'$ the functor $\widehat{\mathscr{T}}' \hatstar (-)$ is t-exact and sends tilting perverse sheaves to tilting perverse sheaves. In particular, since $\Delta^{\Iwu}_e$ is a subobject of $\mathscr{T}^{\Iwu}_{w_0}$ by Lemma~\ref{lem:socle-Deltas}, the tilting object $\mathscr{T} = \pi_! \widehat{\mathscr{T}}[\dim(T)]= \widehat{\mathscr{T}} \hatstar \Delta^{\Iwu}_e$ is a subobject of the tilting object $\mathscr{T}' = \widehat{\mathscr{T}} \hatstar \mathscr{T}^{\Iwu}_{w_0}$. 

What remains to be proved is that $\Av_\IW'(\mathscr{T}')$ is a direct sum of tilting objects associated with dominant weights. However if $\lambda \in X_*(T)$ is not dominant then there exists $s \in S \cap \Wf$ such that $w_\lambda s > w_\lambda$ and $w_\lambda s \in \fW$ (see e.g.~\cite[Lemma~2.4 and its proof]{mr1}). Then if $q_s$ is the projection from $\Fl_G$ to the partial affine flag variety $\Fl_G^s$ associated with $s$, the indecomposable tilting object in $\sfP_{\IW,\Iw}$ associated with $\lambda$ is not killed by the functor $(q_s)_!$ (because the restriction of its image to the $\Iwu^+$-orbit associated with $w_\lambda$ is nonzero). Hence to prove the claim we only have to observe that $(q_s)_! \mathscr{T}'=0$ for any $s \in S \cap \Wf$. This fact follows from the properties that $(q_s)_! \mathscr{T}^{\Iwu}_{w_0}=0$ (which is standard, and can e.g.~be deduced from the fact that $\mathscr{T}^{\Iwu}_{w_0}$ is the projective cover of $\Delta^{\Iwu}_e$ in $\Perv_{\Iwu}(\Loop^+ G/I)$, see~\cite[Proposition~5.26]{modrap1}), and that the functor $\widehat{\mathscr{T}} \hatstar (-)$ commutes with $(q_s)_!$ in the appropriate sense.
\end{proof}

Recall that in a highest weight category, any object is isomorphic to a subquotient of a tilting object (see e.g.~\cite[Proposition~7.17]{riche-hab}).
Therefore, Lemma~\ref{lem:tilting-PIW} has the following consequence.

\begin{cor}
\label{cor:subquot-PIW}
 Any object in $\sfP_{\IW,\Iw}$ is isomorphic to a subquotient of a direct sum of objects of the form $\mathscr{T}^\IW_\mu$ with $\mu \in X_*^+(T)$.
\end{cor}

\subsection{Proof of Theorem~\texorpdfstring{\ref{thm:Av-ff}\eqref{it:Av-ff}}{}}

We have now introduced all the tools needed for the proof of Theorem~\ref{thm:Av-ff}\eqref{it:Av-ff}.
We begin with a general lemma.

\begin{lem}
\label{lem:adjunction-quotients}
 Let $\mathsf{A}$ and $\mathsf{B}$ be abelian categories, let $\mathsf{M}$ be a Serre subcategory of $\mathsf{A}$, let $\mathsf{N}$ be a Serre subcategory of $\mathsf{B}$, and consider the associated quotient functors $\Pi_\mathsf{A} : \mathsf{A} \to \mathsf{A}/\mathsf{M}$ and $\Pi_\mathsf{B} : \mathsf{B} \to \mathsf{B}/\mathsf{N}$. Assume we are given functors $F : \mathsf{A} \to \mathsf{B}$, $G:\mathsf{B} \to \mathsf{A}$ and $F^0 : \mathsf{A}/\mathsf{M} \to \mathsf{B}/\mathsf{N}$, $G^0 : \mathsf{B}/\mathsf{N} \to \mathsf{A}/\mathsf{M}$ such that the following diagram commutes:
 \[
  \xymatrix@C=1.5cm{
  \mathsf{A} \ar[d]_-{\Pi_{\mathsf{A}}} \ar@<0.5ex>[r]^-{F} & \mathsf{B} \ar@<0.5ex>[l]^-{G} \ar[d]^-{\Pi_{\mathsf{B}}} \\
  \mathsf{A}/\mathsf{M} \ar@<0.5ex>[r]^-{F^0} & \mathsf{B}/\mathsf{N}. \ar@<0.5ex>[l]^-{G^0}
  }
 \]
If $F$ is left adjoint to $G$, then $F^0$ is left adjoint to $G^0$.
\end{lem}

\begin{proof}
 We have adjunction morphisms $F \circ G \to \id$ and $\id \to G \circ F$. Composing with $\Pi_{\mathsf{B}}$ and $\Pi_{\mathsf{A}}$ respectively on the left, we deduce morphisms of functors $F^0 \circ G^0 \circ \Pi_{\mathsf{B}} \to \Pi_{\mathsf{B}}$ and $\Pi_{\mathsf{A}} \to G^0 \circ F^0 \circ \Pi_{\mathsf{A}}$. Such morphisms must be induced by morphisms $F^0 \circ G^0 \to \id$ and $\id \to G^0 \circ F^0$. Since our original morphisms satisfy the zig-zag relations, so do these morphisms, which implies the desired claim.
\end{proof}

By Lemma~\ref{lem:Av-D-N}\eqref{it:Av-kills-IC} and~\cite[Corollaires~2--3, p.~368--369]{gabriel}, the restriction of $\Av_\IW$ to the hearts of the perverse t-structures factors through an exact functor
\[
\Av_\IW^{\asp} : \sfP_{\Iw,\Iw}^\asp \to \sfP_{\IW,\Iw};
\]
i.e.~we have $\Av_\IW = \Av_\IW^{\asp} \circ \Pi_{\Iw,\Iw}^{\asp}$. Moreover, by Lemma~\ref{lem:adjunction-quotients},
the functor $\mathsf{F}^\asp := \Pi_{\Iw,\Iw}^{\asp} \circ \mathsf{F}$ is right adjoint to $\Av_\IW^{\asp}$.


The following lemma can be proved using the same arguments as in~\cite[Proof of Lemma~6.4.10]{ar-book}.

\begin{lem}
\label{lem:vanishing-convolution}
Let $\mathscr{F}$ in $\sfP_{\Iw,\Iw}$.
 If $\mathscr{G} \in \sfD_{\Iw,\Iw}$ belongs to the full subcategory generated under extensions by the objects of the form $\IC_w[n]$ with $(w,n) \in \Wf \times \Z_{\leq 0}$ and either $n<0$ or $n=0$ and $w \neq e$, then we have
\[
\Pi_{\Iw,\Iw}^{\asp} \circ \pH^{-1}(\mathscr{G} \star_{\Iw} \mathscr{F})=\Pi_{\Iw,\Iw}^{\asp} \circ \pH^0(\mathscr{G} \star_{\Iw} \mathscr{F})=0.
\]
\end{lem}

The fact that the functor $\Av_\IW^{\asp}$ is fully faithful
is now equivalent to the following lemma.

\begin{lem}
\label{lem:AvIW-ff}
 The adjunction morphism $\id \to \mathsf{F}^{\asp} \circ \Av_\IW^\asp$ is an isomorphism.
\end{lem}

\begin{proof}
 We have to prove that for any $\mathscr{F}$ in $\sfP_{\Iw,\Iw}^{\asp}$ the morphism
 \[
  \mathscr{F} \to \mathsf{F}^{\asp} \circ \Av_\IW^\asp(\mathscr{F})
 \]
is an isomorphism, or in other words that for any $\mathscr{F}$ in $\sfP_{\Iw,\Iw}$ the morphism
\[
 \Pi_{\Iw,\Iw}^\asp(\mathscr{F}) \to \Pi_{\Iw,\Iw}^\asp \circ \mathsf{F} \circ \Av_\IW(\mathscr{F})
\]
induced by the adjunction morphism $\id \to \mathsf{F} \circ \Av_\IW$ is an isomorphism. Now we have
\begin{multline*}
 \mathsf{F} \circ \Av_\IW(\mathscr{F}) = \pH^0 \bigl( \rInd_{\Iw_0}^\Iw (\Delta^\IW_0 \star_{\Iw} \mathscr{F}) [-\dim(T)]\bigr) \\
 \cong \pH^0 \bigl( \rInd_{\Iw_0}^\Iw (\Delta^\IW_0) \star_{\Iw} \mathscr{F} [-\dim(T)]\bigr),
\end{multline*}
and the adjunction morphism considered here is induced by the adjunction morphism
\begin{equation}
\label{eqn:adjunction-Deltae}
 \Delta^{\Iw}_e \to \rInd_{\Iw_0}^\Iw(\Av_\IW(\Delta^{\Iw}_e))[-\dim(T)] = \rInd_{\Iw_0}^\Iw(\Delta^\IW_0)[-\dim(T)].
\end{equation}
Hence to conclude we must show that for any $\mathscr{F}$ in $\sfP_{\Iw,\Iw}$ the latter morphism induces an isomorphism
\[
 \Pi_{\Iw,\Iw}^\asp(\mathscr{F}) \to \Pi_{\Iw,\Iw}^{\asp} \Bigl( \pH^0 \bigl( \rInd_{\Iw_0}^\Iw (\Delta^\IW_0) \star_\Iw \mathscr{F} [-\dim(T)] \bigr) \Bigr).
\]

We start by analyzing the morphism~\eqref{eqn:adjunction-Deltae}. 
The $G$-action on the base point of $\Fl_G$ induces an isomorphism between $\overline{\Fl_{G,w_0}}=\Loop^+ G/\Iw$ and the flag variety $G/B$, from which we obtain equivalences of categories
\begin{gather*}
\Db_{\Iw}(\overline{\Fl_{w_0}},\bk) \cong \Db_B(G/B,\bk), \\
\Db_{(\Iwu^+, \chi^*(\scL_{\mathrm{AS}}))}(\overline{\Fl_{G,w_0}},\bk) \cong \Db_{(U^+,(\chi')^*(\scL_{\mathrm{AS}}))}(G/B,\bk).
\end{gather*}
Using these equivalences, we can apply~\cite[Lemma~12.1]{bezr},\footnote{In that lemma it is assumed that $X^*(T)/\Z\mathfrak{R}$ is free; however this statement also holds in general if $\ell$ is very good, see~\cite[Remark~12.2]{bezr}.}
which guarantees that the complex on the right-hand side of~\eqref{eqn:adjunction-Deltae} is concentrated in non-negative (perverse) degrees, and that moreover we have
\[
\pH^0(\rInd_{\Iw_0}^\Iw(\Delta^\IW_0)[-\dim(T)]) = \pH^{-\dim(T)} \left( \rInd_{\Iw_0}^\Iw(\Delta^\IW_0) \right) \cong \Delta^\Iw_{w_0}.
\]
It follows that~\eqref{eqn:adjunction-Deltae} is induced by a nonzero morphism $\Delta^\Iw_e \to \Delta^\Iw_{w_0}$. Such a morphism is automatically injective, and the composition factors of its cokernel are of the form $\IC_w$ with $w \in \Wf \smallsetminus \{e\}$, see Lemma~\ref{lem:socle-Deltas}. In particular, the cone $\mathscr{C}$ of~\eqref{eqn:adjunction-Deltae} belongs to the full subcategory of $\sfD_{\Iw,\Iw}$ generated under extensions by the objects of the form $\IC_w[n]$ with $(w,n) \in \Wf \times \Z_{\leq 0}$ and either $n<0$ or $n=0$ and $w \neq e$. By Lemma~\ref{lem:vanishing-convolution} it follows that we have
\[
 \Pi_{\Iw,\Iw}^{\asp} \circ \pH^{-1}(\mathscr{C} \star_\Iw \mathscr{F})=\Pi_{\Iw,\Iw}^{\asp} \circ \pH^0(\mathscr{C} \star_\Iw \mathscr{F})=0.
\]
The desired claim follows, using the long exact sequence in cohomology associated with the distinguished triangle
\[
 \mathscr{F} \to \rInd_{\Iw_0}^\Iw(\Delta^\IW_0) \star_\Iw \mathscr{F} [-\dim(T)] \to \mathscr{C} \star_\Iw \mathscr{F} \xrightarrow{[1]}
\]
induced by~\eqref{eqn:adjunction-Deltae} and the exactness of $\Pi_{\Iw,\Iw}^\asp$.
\end{proof}

To complete the proof of Theorem~\ref{thm:Av-ff}\eqref{it:Av-ff}, it only remains to prove that the essential image of $\Av_\IW^\asp$ is closed under subquotients. This property follows from Corollary~\ref{cor:Av-exact} and the following general lemma.

\begin{lem}
\label{lem:essential-image}
Let $\mathsf{A}$ and $\mathsf{B}$ be abelian categories, and let $F : \mathsf{A} \to \mathsf{B}$ be an exact fully faithful functor. Assume that every object in $\mathsf{A}$ has finite length, and that $F(M)$ is simple for any simple object $M$ in $\mathsf{A}$. Then the essential image of $F$ is closed under subquotients.
\end{lem}

\begin{proof}
Our assumption implies that every object in the essential image of $F$ has finite length, and that each of its composition factors is of the form $F(M)$ with $M$ simple in $\mathsf{A}$. Any subquotient of such an object can be obtained by repeatedly taking a cokernel of a morphism from an object of the form $F(M)$ with $M$ simple in $\mathsf{A}$, or a kernel of a morphism to an object of this form. To prove the lemma it therefore suffices to show that if $M,N$ are objects of $\mathsf{A}$ with $M$ simple and if $f : F(M) \to F(N)$, resp.~$f : F(N) \to F(M)$ is an embedding, resp.~a surjection, then $\mathrm{cok}(f)$, resp.~$\ker(f)$, belongs to the essential image of $F$. However, since $F$ is fully faithful we can write $f=F(g)$ for some morphism $g : M \to N$, resp.~$g : N \to M$, and by exactness we then have $\mathrm{cok}(f)=F(\mathrm{cok}(g))$, resp.~$\ker(f)=F(\ker(g))$.
\end{proof}

\begin{rmk}
Lemma~\ref{lem:essential-image} is stated in~\cite[Lemme~4.2.6.1]{bbd} under the extra assumption that $F$ admits left and right adjoints. Our arguments above show that this condition is in fact not necessary.
\end{rmk}

\section{Central and tilting perverse sheaves}
\label{sec:central-tilting}



The considerations in this section do \emph{not} require the assumptions of Sections~\ref{sec:regularity}--\ref{sec:asp-IW} (unless this is explicitly notified).

\subsection{Statement and strategy of proof}

To simplify notation, we now set
\[
\sZ^\IW := \Av_\IW \circ \sZ : \Rep(G^\vee_\bk) \to \sfP_{\IW,\Iw}.
\]
Our goal in this section is to prove the following result.

\begin{thm}
\label{thm:central-tilting}
Assume that 
for any indecomposable factor of the root system $\mathfrak{R}$, the prime $\ell$ is strictly bigger than the corresponding value in the table of Figure~\ref{fig:bounds}.
 Then for any tilting $G^\vee_\bk$-module $V$, the perverse sheaf $\sZ^\IW(V)$ is tilting.
\end{thm}



%

\begin{rmk}
 It is likely that the assumption on $\ell$ in Theorem~\ref{thm:central-tilting} can be weakened; however, in view of~\cite[\S 3.8]{jmw2}, this extension will require a different strategy.
\end{rmk}


The first step towards Theorem~\ref{thm:central-tilting} is the following claim, whose proof can be copied from~\cite[Lemma~25]{ab} or~\cite[Proposition~6.5.7]{ar-book}. (This statement does not require any assumption on $G$ or $\ell$.)

\begin{prop}
\label{prop:central-tilting-1}
 If $V,V'$ are tilting $G^\vee_\bk$-modules such that $\sZ^\IW(V)$ and $\sZ^\IW(V')$ are tilting, then $\sZ^\IW(V \otimes V')$ is tilting.
\end{prop}

Assume now for a moment that $G$ is simple (i.e.~quasi-simple of adjoint type), so that $G^\vee_\bk$ is quasi-simple and simply-connected.
Recall that in this context a dominant coweight $\lambda \in X_*^+(T)$ is called \emph{minuscule} if for any root $\alpha$ we have $\langle \lambda,\alpha \rangle \in \{0,\pm 1\}$, and moreover that in this case we have $\wt(\lambda)=\Wf\lambda$ (where the left-hand side is as in~\S\ref{ss:IW-Perv}). We will also denote by $\alpha_0$ the maximal root of $(G,T)$ (for our given choice of positive roots), and by $\alpha_0^\vee$ the associated coroot.\footnote{Note the slight conflict with the notation of~\S\ref{ss:fixed-points-qmin}. But here the group under consideration is $\bG=G^\vee_\bk$, whose highest \emph{short} root is $\alpha_0^\vee$.} Then $\alpha_0^\vee$ is dominant, and we have $\wt(\alpha_0^\vee)=\Wf\alpha_0^\vee \cup \{0\}$.

The following result is proved in~\cite[\S\S 3.6--3.7]{jmw2}.

\begin{prop}
\label{prop:tilting-jmw}
 Assume that $G$ is simple and that $\ell$ is strictly bigger than the value corresponding to its type in the table of Figure~{\rm \ref{fig:bounds}}. Then any indecomposable tilting $G^\vee_\bk$-module appears as a direct summand of a tensor product of indecomposable tilting modules whose highest weights are either minuscule or equal to $\alpha_0^\vee$. Moreover, if $G$ is of type $\mathbf{A}$ then only minuscule coweights are needed.
\end{prop}

This proposition shows the importance of understanding the objects $\sZ(V)$ (or $\sZ^\IW(V)$) when $V$ is an indecomposable tilting module whose highest weight is either minuscule or equal to $\alpha_0^\vee$. These objects will be studied in the rest of the section. In particular, in~\S\ref{ss:proof-prop-central-tilting-2} we will prove the following claim.

\begin{prop}
\label{prop:central-tilting-2}
Assume that $G$ is simple, and let 
$V$ be an indecomposable tilting $G^\vee_\bk$-module with highest weight $\lambda$.
We assume that
\begin{enumerate}
 \item 
 \label{it:minuscule}
 either $\lambda$ is minuscule;
 \item 
  \label{it:qminuscule}
 or $\lambda=\alpha_0^\vee$, $G$ is not of type $\mathbf{A}$, and $\ell$ satisfies the conditions in Figure~\ref{fig:conditions} with respect to the type of $G^\vee_\bk$.
\end{enumerate}
Then $\sZ^\IW(V)$ is tilting.
\end{prop}

For now,
let us explain why these results imply Theorem~\ref{thm:central-tilting}.

\begin{proof}[Proof of Theorem~\ref{thm:central-tilting}]
Consider the quotient group $G/\rmZ(G)$ (which is semisimple of adjoint type), and the associated affine Grassmannian $\Gr_{G/\rmZ(G)}$. The quotient morphism $G \to G/\rmZ(G)$ induces a morphism
\begin{equation}
\label{eqn:morphism-GrG-quotient}
\Gr_G \to \Gr_{G/\rmZ(G)}.
\end{equation}
We claim that this morphism restricts, on the reduced ind-subscheme associated with each connected component of $\Gr_G$, to a universal homeomorphism onto the reduced ind-scheme associated with a connected component of $\Gr_{G/\rmZ(G)}$. In fact, since the connected components of $\Gr_G$ are permuted by the action of $\Loop G$ it suffices to check this claim for the reduced ind-scheme $(\Gr_G)^0_{\mathrm{red}}$ associated with the connected component containing the base point. Now by~\cite[Proposition~6.6]{pr}, if $\mathscr{D}G$ is the derived subgroup of $G$, the embedding $\mathscr{D}G \hookrightarrow G$ induces an isomorphism
\[
(\Gr_{\mathscr{D}G})^0_{\mathrm{red}} \simto (\Gr_G)^0_{\mathrm{red}}.
\]
We have $\rmZ(\mathscr{D}G)=\rmZ(G) \cap \mathscr{D} G$, and the morphism $\mathscr{D}G/\rmZ(\mathscr{D}G) \to G/\rmZ(G)$ induced by the embedding $\mathscr{D}G \hookrightarrow G$ is an isomorphism, which reduces the proof of our claim to the case $G$ is semisimple. In this case the claim follows from the results of~\cite[\S 6]{pr}; see~\cite[Footnote on p.~3]{br}. (In case $p$ does not divide the order of the fundamental group of $G/\rmZ(G)$, then it is proved in~\cite{pr} that $\Gr_G$ and $\Gr_{G/\rmZ(G)}$ are reduced, and that the above morphism restricts to an \emph{isomorphism} on each connected component of $\Gr_G$.)

Our claim implies that
pushforward along the morphism~\eqref{eqn:morphism-GrG-quotient} induces a fully faithful exact functor
\[
\Perv_{\Loop^+ G}(\Gr_G,\bk) \to \Perv_{\Loop^+(G/\rmZ(G))}(\Gr_{G/\rmZ(G)},\bk).
\]
(To justify the fact that this functor takes values in $\Loop^+(G/\rmZ(G))$-equivariant perverse sheaves, we use the fact that equivariance follows from constructibility with respect to the stratification by $\Loop^+ (G/\rmZ(G))$-orbits; see~\cite[Proposition~2.1]{mv} or~\cite[Proposition~1.10.8]{br}.) 
It can be easily seen
that this functor is monoidal; hence in view of Theorem~\ref{thm:Satake} it provides a monoidal functor $\Rep(G^\vee_\bk) \to \Rep((G/\rmZ(G))^\vee_\bk)$. It is well known that this functor is induced by a natural morphism $(G/\rmZ(G))^\vee_\bk \to G^\vee_\bk$ which identifies $(G/\rmZ(G))^\vee_\bk$ with the simply-connected cover of the derived subgroup of $G^\vee_\bk$.

Since the morphism $\Fl_G \to \Gr_G$ is a (Zariski) locally trivial fibration with fibers $G/B$, and similarly for the morphism $\Fl_{G/\rmZ(G)} \to \Gr_{G/\rmZ(G)}$, we also have a similar claim as above for the natural morphism $\Fl_G \to \Fl_{G/\rmZ(G)}$, which allows to identify, for any connected component of $\Fl_G$, the subcategory of $\sfP_{\IW,\Iw}$ consisting of perverse sheaves supported on this connected component with a full category of the corresponding category for $G/\rmZ(G)$. Since the construction of the functor $\sfZ$ is compatible with these identifications in the obvious way, since $\sZ(V)$ is supported on a single connected component if $V$ is indecomposable, and since the pullback to $(G/\rmZ(G))^\vee_\bk$ of a tilting $G^\vee_\bk$-module is a tilting $(G/\rmZ(G))^\vee_\bk$-module, this reduces the proof of the theorem to the case 
%
$G$ is semisimple of adjoint type. Then, since such a group is a product of simple groups, we can further assume that $G$ is simple.

Recall that (by the general theory of highest weight categories) a direct summand of a tilting object in $\sfP_{\IW,\Iw}$ is again tilting.
If $G$ is of type $\mathbf{A}$, then by Proposition~\ref{prop:tilting-jmw} any indecomposable tilting $G^\vee_\bk$-module is a direct summand of a tensor product of tilting modules with minuscule highest weights; hence in this case the theorem follows from Proposition~\ref{prop:central-tilting-1} and case~\eqref{it:minuscule} in Proposition~\ref{prop:central-tilting-2}. If $G$ is not of type $\mathbf{A}$, the 
claim follows similarly from Proposition~\ref{prop:central-tilting-1}, Proposition~\ref{prop:tilting-jmw}, and Proposition~\ref{prop:central-tilting-2} (cases~\eqref{it:minuscule} and~\eqref{it:qminuscule}).
\end{proof}

\subsection{Stalks and costalks of objects \texorpdfstring{$\sZ^{\IW}(V)$}{Z(V)}}
\label{ss:stalks-costalks}

The proof of the following proposition can be copied from~\cite[Lemma~27]{ab} or~\cite[Proposition~6.5.4]{ar-book}.

\begin{prop}
\label{prop:Euler-central-sheaves}
For any $G^\vee_\bk$-module $V$ and any $\mu \in X_*(T)$ we have
\begin{gather*}
\sum_{i \geq 0} (-1)^i \cdot \dim \bigl( \Hom_{\Db_{\IW}(\Fl_G,\bk)}(\Delta^\IW_\mu, \sZ^\IW(V)[i]) \bigr) = \dim(V_\mu), \\
\sum_{i \geq 0} (-1)^i \cdot \dim \bigl( \Hom_{\Db_{\IW}(\Fl_G,\bk)}(\sZ^\IW(V), \nabla^\IW_\mu[i]) \bigr) = \dim(V_\mu).
\end{gather*}
\end{prop}

\begin{rmk}
\label{rmk:multiplicities-ZIW}
Once Theorem~\ref{thm:central-tilting} will be proved, in view of~\eqref{eqn:Hom-DN-IW} we will know that the morphism spaces
\[
\Hom_{\Db_{\IW}(\Fl_G,\bk)}(\Delta^\IW_\mu, \sZ^\IW(V)[i]) \text{ and } \Hom_{\Db_{\IW}(\Fl_G,\bk)}(\sZ^\IW(V), \nabla^\IW_\mu[i])
\]
vanish for $i\neq0$. Proposition~\ref{prop:Euler-central-sheaves} will then imply that for any $V$ tilting and $\lambda \in X_*(T)$ we have
\[
(\sZ^\IW(V) : \Delta^\IW_\lambda) = (\sZ^\IW(V) : \nabla^\IW_\lambda) = \dim(V_\lambda),
\]
where the first, resp.~second, expression means the multiplicity in a standard, resp.~costandard, filtration.
From this one can in particular deduce that if $V$ has a highest weight $\lambda \in X_*^+(T)$, then the indecomposable tilting object $\mathscr{T}^\IW_\lambda$
is a direct summand of $\sZ^{\IW}(V)$.
\end{rmk}

\begin{cor}
\label{cor:support-ZIW}
Let $V$ in $\Rep(G^\vee_\bk)$ and $\lambda \in X_*^+(T)$. If the $T_\bk^\vee$-weights of $V$ all belong to $\wt(\lambda)$, then $\sZ^\IW(V)$ is supported on the closed subvariety $\bigcup_{\mu \in \wt(\lambda) \cap X_*^+(T)} \overline{\Fl^\IW_{G,\mu}}$.
\end{cor}

\begin{proof}
Let $U$ be the open complement of $\bigcup_{\mu \in \wt(\lambda) \cap X_*^+(T)} \overline{\Fl^\IW_{G,\mu}}$. Then Lemma~\ref{lem:inclusion-orbits} shows that $\Fl^\IW_{G,\nu} \subset U$ iff $\nu \in X_*(T) \smallsetminus \wt(\lambda)$. Our assumption and Proposition~\ref{prop:Euler-central-sheaves} show that the Euler characteristic of the restriction to any such orbit of the perverse sheaf $\sZ^\IW(V)_{|U}$ vanishes. Hence this perverse sheaf is $0$ (because the restriction of a perverse sheaf to a stratum which is open in its support is nonzero, and concentrated in one degree).
\end{proof}

The proof of the following lemma is similar to that of~\cite[Lemma~28]{ab} or~\cite[Lemma~6.5.11]{ar-book}.

\begin{lem}
\label{lem:stalks-ZIW-Wf}
For any $V$ in $\Rep(G^\vee_\bk)$, any $\lambda \in X_*(T)$, any $x \in \Wf$ and any $n \in \Z$ we have
\[
\Hom_{\Db_\IW(\Fl_G,\bk)} \bigl( \Delta^\IW_\lambda, \sZ^\IW(V)[n] \bigr) \cong \Hom_{\Db_\IW(\Fl_G,\bk)} \bigl( \Delta^\IW_{x(\lambda)}, \sZ^\IW(V)[n] \bigr)
\]
and
\[
\Hom_{\Db_\IW(\Fl_G,\bk)} \bigl( \sZ^\IW(V), \nabla^\IW_\lambda[n] \bigr) \cong \Hom_{\Db_\IW(\Fl_G,\bk)} \bigl( \sZ^\IW(V), \nabla^\IW_{x(\lambda)}[n] \bigr).
\]
\end{lem}

The next statement involves the element $\su$ from Proposition~\ref{prop:reg-quotient-RepH}.

\begin{lem}
\label{lem:tilting-0weight}
Assume that either $\ell$ is very good, or $X^*(T)/\Z\mathfrak{R}$ is free.
For any $G^\vee_\bk$-module $V$ which admits a good filtration, we have
 \[
  \dim \Bigl( \Hom_{\sfP_{\IW,\Iw}} \bigl( \Delta_0^{\IW}, \Av_{\IW}(\sZ(V)) \bigr) \Bigr) \leq \dim( V^{\su}).
 \]
For any $G^\vee_\bk$-module $V$ which admits a Weyl filtration, we have
 \[
  \dim \Bigl( \Hom_{\sfP_{\IW,\Iw}} \bigl( \Av_{\IW}(\sZ(V)), \nabla_0^{\IW} \bigr) \Bigr) \leq \dim( V^{\su}).
 \]
\end{lem}

\begin{proof}
We prove the first claim only; the second one can be treated similarly (using the fact that the invariants and coinvariants of $\su$ acting on $V$ have the same dimension).

 We have $\Delta_0^\IW = \Av_\IW(\IC_e)$; hence Theorem~\ref{thm:Av-ff} provides an isomorphism
 \[
  \Hom_{\sfP_{\IW,\Iw}} \bigl( \Delta_0^{\IW}, \Av_{\IW}(\sZ(V)) \bigr) \cong \Hom_{\sfP_{\Iw,\Iw}^{\asp}} \bigl( \Pi_{\Iw,\Iw}^\asp(\IC_e), \Pi_{\Iw,\Iw}^\asp(\sZ(V)) \bigr).
 \]
 Recall now the quotient functor $\Pi_{\Iw,\Iw}^0 : \sfP_{\Iw,\Iw} \to \sfP_{\Iw,\Iw}^0$. It is clear that this functor factors through an exact functor $\sfP_{\Iw,\Iw}^{\asp} \to \sfP_{\Iw,\Iw}^0$; one can easily check that the latter functor identifies $\sfP_{\Iw,\Iw}^0$ with the Serre quotient of $\sfP_{\Iw,\Iw}^{\asp}$ by the Serre subcategory generated by the simple objects $\Pi_{\Iw,\Iw}^{\asp}(\IC_w)$ with $w \in \fW \smallsetminus \Omega$.
 Since $\Pi_{\Iw,\Iw}^\asp(\IC_e)$ is simple,
any nonzero morphism $\Pi_{\Iw,\Iw}^\asp(\IC_e) \to \Pi_{\Iw,\Iw}^\asp(\sZ(V))$ is injective, and its image is not killed by the quotient functor $\sfP_{\Iw,\Iw}^{\asp} \to \sfP_{\Iw,\Iw}^0$; we deduce that this functor induces an injective morphism
\[
 \Hom_{\sfP_{\Iw,\Iw}^{\asp}} \bigl( \Pi_{\Iw,\Iw}^\asp(\IC_e), \Pi_{\Iw,\Iw}^\asp(\sZ(V)) \bigr) \hookrightarrow \Hom_{\sfP_{\Iw,\Iw}^0} \bigl( \sZ^0(\bk), \sZ^0(V) \bigr).
\]
Next, since $\sm_\bk^0=\id$, using Lemma~\ref{lem:monodromy-commutes} we see that
\[
 \Hom_{\sfP_{\Iw,\Iw}^0} \bigl( \sZ^0(\bk), \sZ^0(V) \bigr) = \Hom_{\sfP_{\Iw,\Iw}^0} \bigl( \sZ^0(\bk), \ker(\sm_{V}^0 - \id) \bigr).
\]
And then using Proposition~\ref{prop:reg-quotient-RepH} we obtain an injection
\[
\Hom_{\sfP_{\Iw,\Iw}^0} \bigl( \sZ^0(\bk), \ker(\sm_{V}^0 - \id) \bigr) \hookrightarrow \Hom_\bk(\bk,V^{\su})=V^\su.
\]
The desired inequality follows.
\end{proof}

\subsection{Proof of Proposition~\ref{prop:central-tilting-2}}
\label{ss:proof-prop-central-tilting-2}

%
In this subsection we assume that $G$ is simple, and give the proof of Proposition~\ref{prop:central-tilting-2}. We will use the easy observation that a perverse sheaf in $\Perv_\IW(\Fl_G,\bk)$ is tilting iff its restriction and corestriction to each orbit $\Fl^\IW_{G,\lambda}$ contained in its support are perverse.

Let us first consider case~\eqref{it:minuscule}.
Corollary~\ref{cor:support-ZIW} shows that $\sZ^\IW(V)$ is supported on $\overline{\Fl^\IW_{G,\lambda}}$, which contains the orbit $\Fl^\IW_{G,\mu}$ iff $\mu \in \Wf \cdot \lambda$. The orbit $\Fl^\IW_{G,\lambda}$ is open in this subvariety; hence the restriction and corestriction of the perverse sheaf $\sZ^\IW(V)$ to this stratum are perverse. Using Lemma~\ref{lem:stalks-ZIW-Wf}, we deduce the same property for each orbit $\Fl^\IW_{G,\mu}$ with $\mu \in \Wf \cdot \lambda$, which completes the proof in this case.

Now we consider case~\eqref{it:qminuscule}. In particular we assume that $\lambda=\alpha_0^\vee$, that $G$ is not of type $\mathbf{A}$, and that $\ell$ satisfies the conditions in Figure~\ref{fig:conditions} with respect to the type of $G^\vee_\bk$. The perverse sheaf $\sZ^\IW(V)$ is supported on $\overline{\Fl_{G,\lambda}^\IW} \cup \overline{\Fl_{G,0}^\IW}$, which contains the orbit $\Fl^\IW_{G,\mu}$ iff $\mu \in \Wf \cdot \lambda \cup \{0\}$. As in case~\eqref{it:minuscule}, the restriction and corestriction of $\sZ^{\IW}(V)$ to any orbit $\Fl^\IW_{G,\mu}$ with $\mu \in \Wf \cdot \lambda$ are perverse, so to conclude it suffices to prove the similar claim for $\Fl^\IW_{G,0}$. Let $i : \overline{\Fl^\IW_{G,0}} \to \Fl_G$ be the embedding, and let $j$ be the embedding of the open complement. Then we have distinguished triangles
\begin{gather*}
j_! j^* \sZ^\IW(V) \to \sZ^\IW(V) \to i_* i^* \sZ^\IW(V) \xrightarrow{[1]}, \\
i_* i^! \sZ^\IW(V) \to \sZ^\IW(V) \to j_* j^* \sZ^\IW(V) \xrightarrow{[1]}.
\end{gather*}
What we have proved so far implies that $j^* \sZ^\IW(V)$ admits a filtration with standard subquotients, and a filtration with costandard subquotients; hence both $j_! j^* \sZ^\IW(V)$ and $j_* j^* \sZ^\IW(V)$ are perverse sheaves; it follows that $i^* \sZ^\IW(V)$ is concentrated in perverse degrees $0$ and $-1$, while $i^! \sZ^\IW(V)$ is concentrated in perverse degrees $0$ and $1$. 

By Lemma~\ref{lem:tilting-0weight} 
we have
\[
\dim \Hom_{\sfP_{\IW,\Iw}} \bigl( \Delta_0^{\IW}, \sZ^\IW(V) \bigr) \leq \dim(V^{\su}).
\]
Now by Theorem~\ref{thm:u-regular} the element $\su$ is regular unipotent, so that by Lemma~\ref{lem:fixed-points-qmin} we have
\[
\dim(V^{\su})=\dim(V_0).
\]
Finally, by Proposition~\ref{prop:Euler-central-sheaves} we have
\[
\dim(V_0) = \sum_{n \geq 0} (-1)^n \cdot \dim \Hom_{\Db_{\IW}(\Fl_G,\bk)}(\Delta^\IW_0, \sZ^\IW(V)[n]).
\]
Combining all these remarks, we obtain that
\begin{multline*}
\dim \Hom_{\sfP_{\IW,\Iw}} \bigl( \Delta_0^{\IW}, \sZ^\IW(V) \bigr) \\
\leq \sum_{n \geq 0} (-1)^n \cdot \dim \Hom_{\Db_{\IW}(\Fl_G,\bk)}(\Delta^\IW_0, \sZ^\IW(V)[n]).
\end{multline*}
The remarks in the preceding paragraph imply that the right-hand side equals
\[
\dim \Hom_{\sfP_{\IW,\Iw}} \bigl( \Delta_0^{\IW}, \sZ^\IW(V) \bigr) - \dim \Hom_{\Db_{\IW}(\Fl_G,\bk)} \bigl( \Delta_0^{\IW}, \sZ^\IW(V)[1] \bigr),
\]
so that we necessarily have $\dim \Hom_{\Db_{\IW}(\Fl_G,\bk)} \bigl( \Delta_0^{\IW},\sZ^\IW(V)[1] \bigr)=0$; in other words, the corestriction of $\sZ^\IW(V)$ to $\Fl_{G,0}^\IW$ vanishes in perverse degree $1$, hence is perverse. Similar arguments show that the restriction of this perverse sheaf to $\Fl^\IW_{G,0}$ is perverse, which completes the proof.

%

\section{Proof of Theorem~\texorpdfstring{\ref{thm:main}}{}}
\label{sec:proof-main}


In this section we finally give the proof of Theorem~\ref{thm:main}.


\subsection{Comparison of \texorpdfstring{$\widetilde{\sfP}_{\Iw,\Iw}^0$}{PI0} and \texorpdfstring{$\sfP_{\Iw,\Iw}^0$}{PI0}}
\label{ss:comparison}

In this subsection we assume that either $X^*(T)/\Z\mathfrak{R}$ is free or $\ell$ is very good for $G$. In addition, we assume that for any indecomposable factor of the root system $\mathfrak{R}$, the prime $\ell$ is strictly bigger than the corresponding value in the table of Figure~\ref{fig:bounds}. We start with the following consequence of Theorem~\ref{thm:central-tilting}.

\begin{prop}
\label{prop:subquotients-PIW}
Any object in $\sfP_{\IW,\Iw}$ is isomorphic to a subquotient of an object of the form $\sZ^\IW(V)$ with $V$ in $\Rep(G^\vee_\bk)$.
\end{prop}

\begin{proof}
In view of Corollary~\ref{cor:subquot-PIW}, to prove the proposition it suffices to prove that each object $\mathscr{T}^\IW_\mu$ with $\mu \in X_*^+(T)$ is a subquotient of an object $\sZ^\IW(V)$. This property follows from Theorem~\ref{thm:central-tilting}, see Remark~\ref{rmk:multiplicities-ZIW}.
\end{proof}

Using this proposition we obtain the following.

\begin{cor}
\label{cor:PIasp-PIW}
The functor $\Av_{\IW}^{\asp}$ of Theorem~\ref{thm:Av-ff}\eqref{it:Av-ff} is an equivalence of categories.
\end{cor}

\begin{proof}
By Theorem~\ref{thm:Av-ff}\eqref{it:Av-ff} we already know that $\Av_{\IW}^{\asp}$ is fully faithful, and that its essential image is stable under subquotients. Now this essential image contains the objects $\sZ^{\IW}(V)$ with $V$ in $\Rep(G^\vee_\bk)$; hence Proposition~\ref{prop:subquotients-PIW} implies that it contains all objects, or in other words that $\Av_{\IW}^{\asp}$ is essentially surjective.
\end{proof}

Corollary~\ref{cor:PIasp-PIW} allows us to restate Proposition~\ref{prop:subquotients-PIW} as saying that any object in $\sfP_{\Iw,\Iw}^{\asp}$ is a subquotient of an object $\Pi_{\Iw,\Iw}^{\asp}(\sZ(V))$ with $V$ in $\Rep(G^\vee_\bk)$. 
Now, recall the quotient functor
$\sfP_{\Iw,\Iw}^{\asp} \to \sfP_{\Iw,\Iw}^0$ considered in the proof of Lemma~\ref{lem:tilting-0weight}.
Applying this functor we deduce the following, which establishes Property~\ref{it:tP-P} in~\S\ref{ss:main}.

\begin{cor}
\label{cor:tPI0}
The inclusion $\widetilde{\sfP}_{\Iw,\Iw}^0 \to \sfP_{\Iw,\Iw}^0$ is an equivalence of categories.
\end{cor}

\subsection{Description of \texorpdfstring{$H$}{H}}

In this subsection we assume that $X^*(T)/\Z\mathfrak{R}$ is free, that the $\Z$-module $X_*(T) / \Z\mathfrak{R}^\vee$ has no $\ell$-torsion, and finally that for any indecomposable factor of the root system $\mathfrak{R}$, the prime $\ell$ is strictly bigger than the corresponding value in the table of Figure~\ref{fig:bounds}.

Recall the subgroup $H \subset \rmZ_{G^\vee_\bk}(\su)$ introduced in Proposition~\ref{prop:reg-quotient-RepH}. Corollary~\ref{cor:tPI0} allows us to restate this proposition as providing an equivalence of monoidal categories $\sfP_{\Iw,\Iw}^0 \cong \Rep(H)$. In view of the comments in~\S\ref{ss:main}, to complete the proof of Theorem~\ref{thm:main}, it therefore only remains to prove the following.

\begin{prop}
\label{prop:description-H}
 The embedding $H \subset \rmZ_{G^\vee_\bk}(\su)$ is an equality.
\end{prop}

\begin{proof}
The proof will be based on the results of~\S\S\ref{ss:proof-criterion}--\ref{ss:fixed-points-Zu}, which will be applied to $\mathbf{G}=G^\vee_\bk$. Since the conditions in Figure~\ref{fig:bounds} imply in particular that $\ell$ is good for $G^\vee_\bk$, our assumptions indeed guarantee that these results are applicable.

Recall that the connected components of $\Gr_G$ are in a canonical bijection with the quotient of $X_*(T)$ by the coroot lattice of $(G,T)$, which itself indentifies with $X^*(\rmZ(G^\vee_\bk))$. In particular, every object of $\sfP_{\Loop^+G, \Loop^+G}$ has a canonical decomposition as a direct sum of subobjects parametrized by $X^*(\rmZ(G^\vee_\bk))$; under the geometric Satake equivalence (see~\S\ref{ss:Satake}), this decomposition corresponds to the canonical decomposition of a $G^\vee_\bk$-module according to the (diagonalizable) action of $\rmZ(G^\vee_\bk)$.

The projection $\Fl_G \to \Gr_G$ induces a bijection between the sets of connected components of $\Fl_G$ and $\Gr_G$. Hence any object of $\sfP_{\Iw,\Iw}^0$ also admits a canonical decomposition as a direct sum of subobjects parametrized by $X^*(\rmZ(G^\vee_\bk))$, and from this we see that the embedding $\rmZ(G^\vee_\bk) \hookrightarrow G^\vee_\bk$ factors through $H$; in other words $H$ contains $\rmZ(G^\vee_\bk)$.

We have proved in Theorem~\ref{thm:u-regular} that $\su$ is regular unipotent. Hence the proposition will follow from Lemma~\ref{lem:criterion-equality} provided we prove that for any finite-dimensional tilting $G^\vee_\bk$-module $V$ the embedding
\[
V^{\rmZ_{G^\vee_\bk}(\su)}\hookrightarrow V^{H}
\]
is an equality, or in other words that $\dim(V^{\rmZ_{G^\vee_\bk}(\su)})=\dim(V^H)$. Now, by Proposition~\ref{prop:Vu-V0} the left-hand side is equal to the dimension of the $0$-weight space $V_0$ of $V$. On the other hand, in view of Proposition ~\ref{prop:reg-quotient-RepH} we have
\[
V^H \cong \Hom_H(\bk,V) \cong \Hom_{\sfP_{\Iw,\Iw}^0}(\delta^0, \sZ^0(V)).
\]
Corollary~\ref{cor:PIasp-PIW} allows to identify the category $\sfP_{\Iw,\Iw}^0$ with the category $\sfP_{\IW,\Iw}^0$ of~\S\ref{ss:properties-tilting}, and using
Lemma~\ref{lem:Hom-tiltings-IW} we deduce an isomorphism
\[
\Hom_{\sfP_{\Iw,\Iw}^0}(\delta^0, \sZ^0(V)) \cong \Hom_{\sfP_{\IW,\Iw}}(\Delta^\IW_0, \sZ^\IW(V)).
\]
Theorem~\ref{thm:central-tilting} implies that the right-hand side is equal to the multiplicity $(\sZ^\IW(V):\Delta^\IW_0)$. Finally, by Remark~\ref{rmk:multiplicities-ZIW} this multiplicity is also equal to $\dim(V_0)$, which shows that $\dim(V^{\rmZ_{G^\vee_\bk}(\su)})=\dim(V^H)$ and therefore finishes the proof.
\end{proof}

\subsection{End of the proof}
\label{ss:end-proof}

We have now proved Theorem~\ref{thm:main} under the first possible assumption in~\eqref{it:assumption-main-1} (and assumption~\eqref{it:assumption-main-3}). So, it only remains to explain how to deduce the version of the theorem under the second possible assumption in~\eqref{it:assumption-main-1}, namely that $\ell$ is very good for $G$ (and, of course, assumption~\eqref{it:assumption-main-3}).

In this case one can proceed as follows. The center $\rmZ(G)$ of $G$ can be written (non canonically) as a product of a torus and a finite diagonalizable group. Choose such a decomposition, and let $G'$ be the quotient of $G$ by this finite part. Then the quotient morphism $G \to G'$ induces a morphism of ind-schemes $\Fl_G \to \Fl_{G'}$, which as in the proof of Theorem~\ref{thm:central-tilting} identifies (up to some universal homeomorphism) $\Fl_G$ with a union of connected components of $\Fl_{G'}$. We also have a ``dual'' morphism $(G')^{\vee}_\bk \to G^\vee_\bk$ which is a finite \'etale isogeny, and an obvious commutative diagram involving the geometric Satake equivalences for $G$ and $G'$. If we denote by $\Iw'$ the Iwahori group associated with $G'$ (and the Borel subgroup given by the image of $B$), then we can consider the categories $\sfP_{\Iw,\Iw}^0$ (associated with $G$) and $\sfP_{\Iw',\Iw'}^0$ (associated with $G'$).

The group $G'$ satisfies the conditions considered above, so Theorem~\ref{thm:main} applies to it: there is a regular unipotent element $\su' \in (G')^\vee_\bk$ and an equivalence of monoidal categories
\begin{equation}
\label{eqn:equivalence-G'}
 (\sfP_{\Iw',\Iw'}^0, \star^0) \cong (\Rep(\rmZ_{(G')^\vee_\bk}(\su')), \otimes).
\end{equation}
Using Corollary~\ref{cor:PIasp-PIW} to identify these categories with those considered in~\S\ref{ss:properties-tilting}, one sees that $\sfP_{\Iw,\Iw}^0$ identifies with a monoidal subcategory of $(\sfP_{\Iw',\Iw'}^0, \star^0)$, in fact the direct summand consisting of perverse sheaves supported on the components in the image of $\Fl_G$. On the other hand, denoting by $\su$ the image of $\su'$ in $G^\vee_\bk$, by Proposition~\ref{prop:Z-ureg} the centralizer $\rmZ_{G^\vee_\bk}(\su)$ identifies with the quotient of $\rmZ_{(G')^\vee_\bk}(\su')$ by the kernel of the isogeny $(G')^{\vee}_\bk \to G^\vee_\bk$. As a consequence, the category $\Rep(\rmZ_{G^\vee_\bk}(\su))$ embeds as a direct summand in $\Rep(\rmZ_{(G')^\vee_\bk}(\su'))$ stable under tensor products. One can easily check that the equivalence~\eqref{eqn:equivalence-G'} restricts to an equivalence between the direct summands $\sfP_{\Iw,\Iw}^0$ and $\Rep(\rmZ_{G^\vee_\bk}(\su))$, which finishes the proof.

\subsection{Complement: indecomposability of the objects \texorpdfstring{$\sZ^{\IW}(V)$}{ZIW(V)}}

We finish this section by 
proving the following complement to Theorem~\ref{thm:central-tilting}.

\begin{prop}
 If $V$ is an indecomposable tilting $G^\vee_\bk$-module, $\sZ^{\IW}(V)$ is indecomposable.
\end{prop}

\begin{proof}
By considerations similar to those in~\S\ref{ss:end-proof}, the proof reduces to the case where the conditions in Proposition~\ref{prop:description-H} are satisfied.

 Recall the quotient functor $\Pi^0_{\IW,\Iw} : \sfP_{\IW,\Iw} \to \sfP_{\IW,\Iw}^0$, and consider an indecomposable tilting $G^\vee_\bk$-module $V$. In view of Lemma~\ref{lem:Hom-tiltings-IW}, to prove that $\sZ^{\IW}(V)$ is indecomposable it suffices to prove that $\Pi^0_{\IW,\Iw}(\sZ^{\IW}(V))$ is indecomposable. As observed in the course of the proof of Proposition~\ref{prop:description-H}, Corollary~\ref{cor:PIasp-PIW} allows us to identify $\sfP_{\Iw,\Iw}^0$ with $\sfP_{\IW,\Iw}^0$, so that to conclude it suffices to prove that
 $\sZ^0(V)$ is indecomposable. 
 Now Proposition~\ref{prop:reg-quotient-RepH} and Proposition~\ref{prop:description-H} imply that we have an algebra isomorphism
 \[
  \End_{\sfP_{\Iw,\Iw}^0}(\sZ^0(V)) \cong \End_{\rmZ_{G^\vee_\bk}(\su)}(V);
 \]
hence to conclude we have to show that the algebra $\End_{\rmZ_{G^\vee_\bk}(\su)}(V)$ is local. By Frobenius reciprocity we have an algebra isomorphism
\[
 \End_{\rmZ_{G^\vee_\bk}(\su)}(V) \cong \Hom_{G^\vee_\bk}(V,V \otimes \scO(G^\vee_\bk/\rmZ_{G^\vee_\bk}(\su)))
\]
(where the description of the algebra structure on the right-hand side involves the multiplication in the ring $\scO(G^\vee_\bk/\rmZ_{G^\vee_\bk}(\su))$). 
As seen in the course of the proof of Lemma~\ref{lem:Spr-resolution-isom}, the nilpotent cone $\mathcal{N}_{G^\vee_\bk}$ of $G^\vee_\bk$ is isomorphic to its unipotent cone; 
using Lemma~\ref{lem:U-sc} we deduce a $G^\vee_\bk$-equivariant isomorphism of algebras
\[
\scO(G^\vee_\bk/\rmZ_{G^\vee_\bk}(\su)) \cong \scO(\mathcal{N}_{G^\vee_\bk}).
\]
The dilation action of $\Gm$ on $\mathcal{N}_{G^\vee_\bk}$ induces a grading on $\scO(\mathcal{N}_{G^\vee_\bk})$, such that $G^\vee_\bk$ acts by graded algebra automorphisms, and whose degree-$0$ part is $\bk$. This structure induces a grading on $\End_{\rmZ_{G^\vee_\bk}(\su)}(V)$, whose degree-$0$ part is the algebra $\End_{G^\vee_\bk}(V)$, which is local since $V$ is indecomposable. We deduce that $\End_{\rmZ_{G^\vee_\bk}(\su)}(V)$ is local, see~\cite[Theorem~3.1]{gg}.
\end{proof}


\end{document}